%% file: MIYATAKEv2.tex
\documentclass[11pt,a4paper]{article}

\DeclareMathAlphabet{\mathcal}{OMS}{cmsy}{m}{n}
\usepackage[top=25mm,bottom=35mm,left=20mm,right=20mm]{geometry} 

\usepackage{url}
\usepackage{graphicx}
\usepackage{mathptmx}     
\usepackage{mathtools}
\mathtoolsset{showonlyrefs=true}
\DeclarePairedDelimiter\paren{\lparen}{\rparen}
\usepackage{color}
\usepackage{doi}
\usepackage{amsmath}
\usepackage{amsfonts}
\usepackage{amssymb}
\usepackage{amsthm}
\usepackage{bm}

\usepackage[linesnumbered,ruled]{algorithm2e}

\newcommand{\bbR}{\mathbb{R}}
\newcommand{\bbT}{\mathbb{T}}
\newcommand{\px}{\partial_x}
\newcommand{\rmd}{\mathrm{d}}

\newcommand{\calH}{\mathcal{H}}

\newcommand{\fracdel}[2]{\cfrac{\delta #1}{\delta #2}}

\newcommand{\ba}{\bm a}
\newcommand{\bc}{\bm c}
\newcommand{\bd}{\bm d}
\newcommand{\be}{\bm e}
\newcommand{\bff}{\bm f}
\newcommand{\bg}{\bm g}
\newcommand{\br}{\bm r}
\newcommand{\bs}{\bm s}
\newcommand{\bu}{\bm u}
\newcommand{\bv}{\bm v}

\newcommand{\by}{{\bm y}}
\newcommand{\bz}{{\bm z}}

\DeclareMathOperator{\veco}{vec}
\DeclareMathOperator{\sech}{sech}
\DeclareMathOperator{\diag}{diag}

\DeclareMathOperator{\rank}{rank}
\DeclareMathOperator*{\mtmax}{\textsf{max}}

\newcommand{\tn}{{\tilde{n}}}
\newcommand{\tr}{{\tilde{r}}}

\theoremstyle{definition}
\newtheorem{proposition}{Proposition}

\newtheorem{remark}{Remark}

\usepackage{tikz}
\usepackage{pgfplotstable}
\usetikzlibrary{arrows,positioning,plotmarks,external,patterns,angles,
decorations.pathmorphing,backgrounds,fit,shapes,graphs,calc}

\begin{document}
\title{Structure-preserving model reduction for dynamical systems \\ with a first integral}
\author{
Yuto Miyatake\thanks{Cybermedia Center, Osaka University, 
			1-32 Machikaneyama, Toyonaka, Osaka 560-0043, Japan,
\href{mailto:miyatake@cas.cmc.osaka-u.ac.jp}{miyatake@cas.cmc.osaka-u.ac.jp}}
}

\maketitle

\begin{abstract}
Since the expense of the numerical integration of large scale dynamical systems
 is often computationally prohibitive,
model reduction methods, which approximate such systems by simpler and 
much lower order ones,
are often employed to reduce the computational effort.
In this paper, for dynamical systems with a first integral,
new structure-preserving model reduction approaches
are presented 
that yield reduced-order systems while preserving the first integral.
We apply energy-preserving integrators to the reduced-order systems and show some numerical experiments
that demonstrate the favourable behaviour of the proposed approaches.
\end{abstract}

\section{Introduction}
\label{sec1}

Since the expense of the numerical integration of large scale dynamical systems
is often computationally prohibitive,
model reduction methods, which approximate high dimensional systems 
by simpler and much lower order ones, are often employed 
to reduce the computational effort~\cite{as01}.
The proper orthogonal decomposition  (POD) method with Galerkin projection,
which was first introduced by Moore~\cite{mo81},
is one of standard data-driven model reduction methods.
This method extracts a few basis vectors that fit the empirical solution data 
with a good accuracy,
and project the high dimensional system to the subspace spanned by the basis vectors.
The POD-Galerkin approach can often provide an efficient surrogate system,
and has found applications in a wide range of areas such as
structural dynamics~\cite{as03}, fluid mechanics~\cite{hl12,il00,rc04}, 
and time-dependent partial differential equations~\cite{kv01,sk04}.
However, when the vector field of the original system is nonlinear,
the complexity of evaluating the nonlinear term of the reduced-order system
remains as expensive as that of the original problem.
To resolve this issue,
Chaturantabut and Sorensen proposed the 
discrete empirical interpolation method (DEIM)
based on the POD-Galerkin method and an interpolatory projection~\cite{cs10,cs12}.

Though the aforementioned data-driven approaches 
work preferably for many applications,
they rarely inherit underlying mathematical structures of the original system,
such as symmetry, symplecticity and energy-preservation.
For dynamical systems with some mathematical structures,
numerical integrators that inherit such properties,
referred to as geometric numerical integrators or 
structure-preserving integrators, are often preferred,
since they usually produce qualitatively better numerical solutions than
standard general-purpose integrators such as the famous fourth-order 
explicit Runge--Kutta method (see, e.g.~\cite{hlw06}).
Therefore, model reduction while preserving such properties would be 
preferred: for example, 
if the reduced-order system inherits the mathematical structures, 
one could easily choose an appropriate 
numerical integrator for the reduced-order system.
Structure-preserving model reduction methods have received attention 
in recent years (see~\cite{ah17,cb16,gw17,pm16} and references therein).

In this paper, we are concerned with a dynamical system with a first integral, 
i.e. a dynamical system with a conservation law.
Such a system can always be formulated as a skew-gradient system of the form
\begin{align} \label{sys1}
\frac{\rmd}{\rmd t}\by = S(\by) \nabla_\by H(\by), \quad \by(0) = \by_0 \in \bbR^n,
\end{align}
where $S(\by)\in\bbR^{n\times n}$ is a skew-symmetric matrix, 
and the function $H:\bbR^n\to \bbR$ is assumed to be 
sufficiently differentiable~\cite{qc96}.
Indeed, the function $H$ is constant along the solution:
\begin{align*}
\frac{\rmd}{\rmd t}H(\by) 
= \nabla _\by H(\by)^\top \dot{\by} 
= \nabla_\by H(\by) ^\top S(\by) \nabla_\by H(\by) = 0
\end{align*}
due to the skew-symmetry of $S(\by)$, where the dot stands for the differentiation
with respect to $t$.

When $S(\by)$ is a constant skew-symmetric matrix, that is, 
it is independent of $\by$,
several structure-preserving model reduction methods have been studied.
If $S$ is of the form
\begin{align} \label{eq:SJ}
S =  J_{2\tn}^{-1} =
\begin{bmatrix}
0_\tn & I_\tn \\
-I_\tn & 0_\tn
\end{bmatrix},
\quad n = 2\tn,
\end{align}
where $0_n, I_n\in\bbR^{n\times n}$ denote 
the zero and identity matrices, respectively,
the corresponding system is called a Hamiltonian system.
Peng and Mohseni~\cite{pm16} proposed model reduction techniques 
that find a lower order Hamiltonian system, to which
any structure-preserving integrators developed for Hamiltonian systems
can be applicable.
For the case $S(\by)$ is a constant skew-symmetric matrix 
but is not necessarily of the form \eqref{eq:SJ},
Gong et al.~\cite{gw17}  proposed a model reduction approach that yields
a lower-order skew-gradient system with a constant skew-symmetric matrix.
These structure-reserving model reduction methods are briefly reviewed in
Section~\ref{sec2}.
For other structure-preserving model reduction methods, 
see, for example,~\cite{ah17,cb16} and references therein.

In the line of these research, we are concerned with the case
$S(\by)$ depends on $\by$.
This situation often arises, for example, as a Hamiltonian system with some constraints
and from discretizing a Hamiltonian partial differential equation (PDE).
The simple application of the approach~\cite{gw17} gives 
a lower order skew-gradient system;
however, the computational complexity for evaluating the vector field 
may still depend on $n$ (the size of the original problem).
In this paper, we study structure-preserving model reduction techniques
so that the vector field of the reduced-order system can be evaluated efficiently. 
We classify target systems into two types.
First class is the case $S(\by)$ depends linearly on $\by$, 
and has a specific structure such as $S(\by) = YD+DY$, 
where $D\in\bbR^{n\times n}$ is a constant skew-symmetric matrix and 
$Y = \diag (\by) \in \bbR^{n\times n}$.
In this case,
we show that the computational complexity for the reduced-order system 
based on the approach~\cite{gw17} 
is already independent of $n$.
We shall develop a new approach for more general cases, as a second class,
based on the approach~\cite{gw17} and DEIM.

This paper is organized as follows.
In Section~\ref{sec2},
the proper orthogonal decomposition method with Galerkin-projection,
the discrete empirical interpolation method and 
some structure-preserving model reduction methods are briefly reviewed.
Structure-preserving model reduction methods for \eqref{sys1} are discussed in
Sections~\ref{sec:linear} and ~\ref{sec:nonlinear}.
Section~\ref{sec:linear} considers the first class,
and Section~\ref{sec:nonlinear} treats the second, i.e. general cases.
We demonstrate the effect of the methods by some numerical results in 
Section~\ref{sec:num}.
Finally, concluding remarks are given in Section~\ref{sec:conclusion}.

\section{Preliminaries: POD, DEIM and some structure-preserving model reduction methods}
\label{sec2}

In this section,
we briefly review the proper orthogonal decomposition method 
with Galerkin-projection,
the discrete empirical interpolation method and 
some structure-preserving model reduction methods.

\subsection{Model reduction with Galerkin-projection}

Model reduction methods considered in this paper, except for
those in Section~\ref{subsubsec:smr}, are based on the Galerkin projection.
The basic procedure to construct a reduced-order system is summarized below.

Let us consider a system of ordinary differential equations of the form
\begin{align} \label{eq:odef}
\frac{\rmd}{\rmd t} \by = \bff (\by), \quad \by(0) = \by_0\in\bbR^n,
\end{align}
as a full-order model,
where $\bff:\bbR^n \to \bbR^n$ is supposed to be sufficiently smooth.
A standard way of constructing a reduced-order system is to
project the solution of \eqref{eq:odef}  onto an appropriate subspace of $\bbR^n$.
Assume that the flow $\by(t)$ can be well approximated 
in a lower dimensional subspace, i.e.
a linear combination of some basis vectors $\bv_i\in\bbR^n$ ($i=1,\dots,r$):
\begin{align} \label{rl1}
\by(t) \approx \sum_{i=1}^r z_i(t) \bv_i,
\end{align}
where $r\ll n$. 
Without loss of generality, the basis vectors are chosen such that they are orthonormal.
Let $V := [\bv_1,\dots,\bv_r]\in\bbR^{n\times r}$. Then $V^\top V = I_r$.
By using this notation and $\bz(t) := [z_1(t),\dots,z_r(t)]^\top$, 
the relation \eqref{rl1} can be written as
\begin{align}
\by(t) \approx V \bz(t).
\end{align}
Substituting $V \bz$ into $\by$ in \eqref{eq:odef} yields the overdetermined system
\begin{align*}
V\frac{\rmd}{\rmd t} \bz = \bff(V\bz).
\end{align*} 
Applying the Galerkin method by multiplying $V^\top$ from the left 
leads to the reduced-order system
\begin{align} \label{eq:ros1}
\frac{\rmd}{\rmd t}\bz = V^\top \bff(V\bz), \quad \bz (0) = V^\top \by_0.
\end{align}

\subsection{The proper orthogonal decomposition method}

The proper orthogonal decomposition (POD) method is a popular
approach of finding an appropriate matrix $V$ based on empirical solution data~\cite{mo81}.

The POD method seeks to extract important information from empirical solution data,
called snapshots, of the full-order system.
A snapshot matrix $Y$ consists of either numerical solutions or observed data at 
some time instances $t=t_1,t_2,\dots,t_s$.
Let $Y := [\by_1,\dots, \by_s]\in\bbR^{n\times s}$,
where $\by _i \approx \by (t_i)$.
We then consider the following optimization problem
\begin{align} \label{opt1}
\min_{\rank (V) = r} \sum_{j=1}^s \| \by_j - VV^\top \by_j \| ^2
\quad \text{such that} \quad V^\top V = I_r,
\end{align}
where $\| \cdot \| $ denotes the $2$-norm in the Euclidean space. 
The optimal solution to this problem is given 
by the singular value decomposition (SVD) for $Y$.
Let $\bv_1,\dots,\bv_r$ be the left singular vectors of $Y$ 
corresponding to the first $r$ leading nonzero singular values.
Then the POD matrix $V := [\bv_1,\dots, \bv_r]$ solves 
the above optimization problem \eqref{opt1}.
If $\rank(Y) = d$ and $\sigma_1 \geq \sigma_2\geq \cdots \geq \sigma_d > 0$,
it follows  that
\begin{align}
\sum_{j=1}^s \| \by_j - VV^\top \by_j \| ^2 = \sum_{j=r+1}^d \sigma_j^2,
\end{align}
and thus the dimension $r$ of the reduced-order system 
is usually set such that it satisfies
\begin{align*}
\frac{\sum_{j=r+1}^d \sigma_j^2}{\sum_{j=1}^d \sigma_j^2} < \epsilon
\end{align*}
for a small constant $\epsilon$ ($0<\epsilon \ll 1$).

If the vector field $\bff$ is linear, that is, 
$\bff(\by) = A \by$ for some constant matrix $A\in\bbR^{n\times n}$,
the reduced-order system \eqref{eq:ros1} becomes
\begin{align} \label{eq:ros2}
\frac{\rmd}{\rmd t} \bz = \hat{A}\bz,
\end{align}
where $\hat{A} = V^\top A V \in \bbR^{r\times r}$.
Since the matrix $\hat{A}$ can be computed in the off-line stage,
the computational complexity of evaluating the vector field $\hat{A}\bz$ depends
only on $r$ and is independent of $n$ (the size of the original problem).

\subsection{The discrete empirical interpolation method}
\label{subsec:deim}

In general, the vector field $\bff$ is often nonlinear.
Let 
\begin{align}
\bff (\by) = A \by + \bg (\by),
\end{align}
where $A\in\bbR^{n\times n}$ is a constant matrix and 
$\bg:\bbR^n \to \bbR^n$ denotes a nonlinear part.
In this case, the reduced-order system \eqref{eq:ros1} becomes
\begin{align}
\frac{\rmd}{\rmd t}\bz = \hat{A}\bz + V^\top \bg ( V\bz),
\end{align}
and notice that the computational complexity for the second term $V^\top \bg ( V\bz)$
may still depend on $n$ due to the nonlinearlity:
one first needs to compute the state variable $\by := V\bz$ in the original
coordinate system, next evaluate the nonlinear vector field  $\bg (\by)$, 
and then project $\bg (\by)$ back onto the column space of $V$. 
This could make solving the reduced-order system more expensive than solving the
original full-order system.

The discrete empirical interpolation method (DEIM) was proposed by 
Chaturantabut and Sorensen~\cite{cs10} to reduce the computational complexity of evaluating the nonlinear term.
Let $\bg(t) := \bg(V\bz(t))$ to simplify the notation.
We consider the approximation to $\bg(t)$ by means of 
a constant matrix 
$U \in \bbR^{n\times m}$ ($m\ll n$) and
a time-dependent vector $\bc(t) \in \bbR^m$:
\begin{align}
\bg(t) \approx U \bc(t).
\end{align}
The DEIM tells us how to construct appropriate  $U$ and $\bc(t)$.
We first explain the construction of $\bc (t)$ assuming we already have $U$.
We require that $\bg(t)$ and $U\bc(t)$ are equal for  $m$ variables out of $n$ variables,
i.e.
\begin{align} \label{deim:req1}
\bg_{\varrho_i} (t) = U_{\varrho_i} \bc(t), \quad i = 1,\dots, m,
\end{align}
where $U_{\varrho_i}$ denotes the $\varrho_i$th row of $U$.
By using $P := [\be_{\varrho_1},\dots,\be_{\varrho_m}]\in\bbR^{n\times m}$
where $\be_{\varrho_i}$ denotes the $\varrho_i$th column of the identity matrix 
of size $n$-by-$n$, the condition \eqref{deim:req1} can be rewritten as
\begin{align}
P^\top \bg(t) = P^\top U \bc(t).
\end{align}
Now, let us assume that $P^\top U \in \bbR^{m\times m}$ is nonsingular.
Then, $\bc(t)$ is given by
\begin{align}
\bc(t) = (P^\top U)^{-1} P^\top \bg(t),
\end{align}
and thus $V^\top \bg(t)$ can be approximated by
\begin{align}
V^\top \bg(t) \approx \underbrace{V^\top U (P^\top U)^{-1}}_{r\times m} 
\underbrace{P^\top \bg(t)}_{m\times 1}.
\end{align}
Note that $V^\top U (P^\top U)^{-1}\in\bbR^{r\times m}$ can be computed 
in the off-line stage, and thus
if the computational complexity for $P^\top \bg (t)$ is independent of $n$,
the complexity for the approximation of $V^\top \bg(t)$ is also independent of $n$. 
Note that the computational complexity for $P^\top \bg (t)$ 
varies from problem to problem.
It is independent of $n$ for many applications, though there are some exceptions
(see~\cite{cs10} for more details).

The procedure for constructing the matrices $U$ and $P$
is summarized in Algorithm~\ref{algo:DEIM},
where $[|\rho|, \varrho] = \mtmax \{ | \bg | \}$
implies that $|\rho| = \max \{ | \bg | \}$ and 
$\varrho$ is the first index of the maximum value(s).
Let $\bg_1, \bg_2, \dots, \bg_s$ be snapshot data for 
$\bg(\by)$ at some time instances,
and let $G := [\bg_1, \bg_2, \dots, \bg_s]$.
Applying the SVD to this matrix gives the POD basis vectors
$ \bu_1,\bu_2,\dots,\bu_m$.
The matrix $P$ can be constructed by a greedy algorithm.
Initially, the first interpolation index $\rho_1\in\{1,2,\dots,n\}$
is selected such that it is 
corresponding to the largest magnitude of the first basis function $\bu_1$.
The remaining indices $\rho_i$ ($i=2,3,\dots,m$)
are selected such that they correspond  to the largest magnitude
of the residual (defined in line 5).
Note 
that $P^\top U$ is nonsingular if $\rho \neq 0$~\cite{cs10}.

\begin{algorithm}
\label{algo:DEIM}
\caption{DEIM}
\SetKwInOut{input}{Input}
\SetKwInOut{output}{Output}
\input{ $\{ \bu_l\}_{l=1}^m \in \bbR^n$ linearly independent}\output{$U\in\bbR^{n\times m}$ and $P\in\bbR^{n\times m}$}
$[|\rho|, \varrho_1] = \mtmax \{ | \bu_1 | \}$ \\
$U = [\bu_1] $, $P = [\be_{\varrho_1}]$ \\
\For{$l=2$ to $m$}{
Solve $(P^\top U) \bc = P^\top \bu_l$ for $\bc$\\
$\br = \bu_l - U\bc$\\
$[|\rho|, \varrho_l] = \mtmax \{ | r | \}$ \\
$U\leftarrow [ U \ \bu_l]$, $P = [P \ \be_{\varrho_l} ]$
}
\end{algorithm}

\subsection{Structure-preserving model reduction}

We here review two structure-preserving model reduction methods.

\subsubsection{Hamiltonian systems}
\label{subsubsec:smr}
Consider Hamiltonian systems
\begin{align} \label{eq:hs}
\frac{\rmd}{\rmd t}\by =  J_{2\tn}^{-1} \nabla_\by H(\by), \quad 
\by(0) = \by_0 \in \bbR^n,
\end{align}
where $J_{2\tn}^{-1}$ is defined in \eqref{eq:SJ}.
The reduced-order system based on the standard Galerkin projection reads
\begin{align}
\frac{\rmd}{\rmd t} \bz = V^\top  J_{2\tn}^{-1} \nabla_\by H(V\bz).
\end{align}
However, this reduced-order system is not always a Hamiltonian system.

Peng and Mohseni~\cite{pm16} proposed 
structure-preserving model reduction methods
that find a reduced-order Hamiltonian system~\cite{pm16},
and their idea is briefly summarized below.
We assume that the matrix $V\in\bbR^{n\times r}$ is a symplectic matrix:
\begin{align} \label{cond:symmat}
V^\top J_{2\tn}^{-1} V = J_{2\tr}^{-1},
\end{align}
and define the symplectic inverse of the matrix $V$, denoted by $V^+$, by
\begin{align}
V^+ := J_{2\tr}^\top V^\top J_{2\tn}. 
\end{align}
Here the constraint $V^\top V = I_r$ is not required.
It follows that
\begin{align}\label{smp1}
V^+ V = I_{2\tr}, \quad V^+ J_{2\tn}^{-1} = J_{2\tr}^{-1} V^\top.
\end{align}
For the overdetermined system
\begin{align}
V \frac{\rmd}{\rmd t} \bz = J_{2\tn} \nabla_\by H(V\bz),
\end{align}
which is obtained by substituting $V\bz$ into $\by$ in \eqref{eq:hs},
multiplying $V^+$ from the left and using the properties \eqref{smp1} yield
\begin{align}
\frac{\rmd}{\rmd t} \bz = V^+ J_{2\tn}^{-1} \nabla_{\by} H(V\bz)
=  J_{2\tr}^{-1} V^\top  \nabla_{\by} H(V\bz)
= J_{2\tr}^{-1} \nabla _{\bz} H(V\bz).
\end{align}
The last equality follows due to the chain rule
$\nabla _\bz H(V\bz) = V^\top \nabla _\by H(V\bz)$,
which will be frequently used in this paper.
Observe that this is a Hamiltonian system 
for the Hamiltonian $\tilde{H}(\bz) := H(V\bz)$.
This approach is referred to as the symplectic model reduction.

The standard POD matrix does not always satisfy \eqref{cond:symmat}, and 
several approaches to construct an appropriate symplectic matrix $V$ 
were proposed in~\cite{pm16}.
For nonlinear Hamiltonian problems, the authors also proposed the so-called 
symplectic DEIM to reduce the computational complexity.

\subsubsection{General constant skew gradient systems}
\label{subsec:g}
If $S$ in \eqref{sys1} is a constant skew-symmetric matrix but is not of the form \eqref{eq:SJ},
the original system may not be a Hamiltonian system.
In this case, the reduced-order system based on the standard Galerkin projection
\begin{align} \label{eq:g1}
\frac{\rmd}{\rmd t}\bz = V^\top S \nabla_\by H(V\bz)
\end{align}
is not a skew-gradient system in general.
Further, the symplectic model reduction, 
which makes use of the structure of $J_{2\tn}^{-1}$,
is not applicable.

Below we review the approach by Gong et al.~\cite{gw17}.
The key is that formally inserting 
$VV^\top \in\bbR^{n\times n}$ between
$S$ and $\nabla H(V\bz)$ in \eqref{eq:g1}
yields a small skew-gradient system
\begin{align} \label{eq:g2}
\frac{\rmd}{\rmd t}\bz = V^\top S V V^\top \nabla_\by H(V\bz)  = S_r \nabla_{\bz} \tilde{H}(\bz),
\end{align}
where $S_r := V^\top S V$ and $\tilde{H}(\bz) := H(V\bz) $.
Since $VV^\top \neq I_n$ in general, 
the system \eqref{eq:g2} differs from \eqref{eq:g1},
and thus we need to carefully consider the relation 
between \eqref{eq:g2} and the original system.
As long as $V$ is the POD matrix 
generated from the standard snapshot solution data, the relation might be subtle.
Then, the authors proposed to device the snapshot matrix $Y$:
\begin{align} \label{gmat}
Y = [\by_1,\dots, \by_s, \mu \nabla_\by  H (\by_1) ,\dots, 
\mu \nabla _\by H (\by_s)]\in\bbR^{n\times s}
\end{align}
for some constant $\mu>0$.
The left singular vectors corresponding to the leading nonzero singular values
extract the information of the gradient $\nabla _\by H(\by)$ as well as $\by$.
Thus, for the POD matrix for this snapshot matrix, 
$VV^\top \nabla_\by H(\by)$ could be a good approximation
to $\nabla_\by H(\by)$.
The error analysis was also given in~\cite{gw17}.

\section{Structure-preserving model reduction for particular skew-gradient systems}
\label{sec:linear}

We now consider the case that $S(\by)$ in \eqref{sys1} may depend in $\by$.
The approach by Gong et al.~\cite{gw17}, 
which was summarized in Section~\ref{subsec:g},
is also applicable to the general cases to find the reduced-order skew-gradient system
\begin{align} \label{eq:rs1}
\frac{\rmd}{\rmd t}\bz = V^\top S(V\bz ) V V^\top \nabla_\by H(V\bz)
 =S_r (\bz) \nabla_{\bz} \tilde{H}(\bz),
\end{align}
where $S_r(\bz) := V^\top S(V\bz) V$ and $\tilde{H}(\bz) := H(V\bz) $.
But the computational complexity of evaluating the matrix $S_r(\bz)$ is
not always independent of the size of the full-order system.
In this section,
we show that if $S(\by)$ is of the form 
\begin{align}
S(\by) = YD + DY,
\end{align}
where $D\in\bbR^{n\times n}$ is a constant skew-symmetric matrix
and $Y = \diag (\by) \in \bbR^{n\times n}$,
 the computational complexity of 
evaluating $S_r(\bz)$ is independent of $n$.

Below, we use the Kronecker product, which is defined by
\begin{align}
A \otimes B := 
\begin{bmatrix}
a_{11} B & a_{12}B & \cdots & a_{1n} B \\ 
a_{21} B & a_{22}B & \cdots & a_{2n} B \\
\vdots & \vdots & \ddots & \vdots \\
a_{m1} B & a_{m2}B & \cdots & a_{mn} B
\end{bmatrix}
\in \bbR^{mp \times nq}
\end{align}
for $A = [a_{ij}] \in \bbR^{m\times n}$ and $B \in \bbR^{p\times q}$.
To simplify the notation, 
we shall consider the computational complexity for $S_r (\bz)$ after vectorizing it. 
For this aim,
for $A=[\ba_1, \ba_2, \dots , \ba_n]\in\bbR^{m\times n}$, 
we define the vec operator,
$\veco: \bbR^{m\times n} \to \bbR^{mn}$,
by
\begin{align}
\veco (A) := \begin{bmatrix}
\ba_1 \\ \ba_2 \\ \vdots \\ \ba_n
\end{bmatrix},
\end{align}
and the inverse vec operator, $\veco^{-1}: \bbR^{mn} \to \bbR^{m\times n}$,
by
$\veco^{-1} (\veco(A)) = A$.
We frequently use the following property: 
for $A\in\bbR^{m\times n}$ and $B\in\bbR^{n\times p}$,
it follows that
\begin{align} \label{vecpro1}
\veco (AB) =( I_p \otimes A) \veco (B) = (B^\top \otimes I_m) \veco (A),
\end{align}
where $I_n\in\bbR^{n\times n}$ is the identity matrix (see~e.g.~\cite[p.~275]{de97}).

Let us consider the computational complexity for $\veco (V^\top S(V\bz) V)$,
which is equivalent to discuss the complexity for $V^\top S(V\bz) V$.
By using \eqref{vecpro1}, it follows that
\begin{align}
\veco (V^\top (YD + DY) V)
&= (I_r \otimes V^\top) \veco ((YD + DY) V) 
= (I_r \otimes V^\top) (V^\top \otimes I_n) \veco (YD + DY) \\
&= \underbrace{(V\otimes V)^\top}_{r^2 \times n^2} 
\underbrace{\paren*{(D^\top \otimes I_n) + (I_n \otimes D)}}_{n^2\times n^2}
\underbrace{\veco (Y)}_{n^2 \times 1}.
\end{align}
Note that $\veco (Y) = [y_1, 0,\dots, 0 ,y_2, 0,\dots , y_n]^\top \in \bbR^{n^2 \times 1}$,
where only $(nk+1)$th elements ($k=0,\dots ,n-1$) are nonzero.
We define $\tilde{D}\in\bbR^{n^2\times n}$ by collecting the $(nk+1)$th columns ($k=0,\dots ,n-1$) 
of the matrix $\paren*{(D^\top \otimes I_n) + (I_n \otimes D)}$.
Note that $\tilde{D}$ is explicitly given by
\begin{align}
\tilde{D} = 
-
\begin{bmatrix}
\diag (d_{11},\dots,d_{1n} ) \\
\diag (d_{21},\dots,d_{2n} ) \\
\vdots \\
\diag (d_{n1},\dots,d_{nn} )
\end{bmatrix}
+ 
\begin{bmatrix}
\bd_1 & 0 & \cdots & 0 \\
0 & \bd_2 & \cdots & 0 \\
\vdots & \vdots & \ddots & \vdots \\
0 & \cdots & 0 & \bd_n
\end{bmatrix},
\end{align}
where $D = [\bd_1, \bd_2, \dots, \bd_n]$.
By using this notation, $\veco (V^\top S(\by) V)$ can be simplified as follows:
\begin{align}
\veco (V^\top S(\by)V)
&=  \underbrace{(V\otimes V)^\top}_{r^2 \times n^2} 
\underbrace{\tilde{D}}_{n^2 \times n}
\underbrace{\by}_{n\times 1}.
\end{align}
Then, substituting $V\bz$ into $\by$ in $\veco (V^\top S(\by)V)$ yields
\begin{align}
\veco (V^\top S(V\bz)V)
&=  
\underbrace{(V\otimes V)^\top}_{r^2 \times n^2} 
\underbrace{\tilde{D}V}_{n^2 \times r}
\underbrace{\bz}_{r\times 1},
\end{align}
and thus
\begin{align}
S_r (\bz) = V^\top S(V\bz)V = \veco^{-1} \paren*{(V\otimes V)^\top \tilde{D}V \bz}.
\end{align}
Since $(V\otimes V)^\top \tilde{D}V \in \bbR^{r^2 \times r}$ 
can be computed in the off-line stage,
the computational complexity of evaluating $S_r(\bz)$ is independent of $n$.

\begin{remark}
The above discussion can be applicable to a bit more general cases.
For example, the approach applies, in a similar manner,
to $S(\by) = YD+DY + D_\text{c}$, where $D_\text{c} $ is a constant skew-symmetric matrix, 
and $S(\by)$ with a different ordering of $Y$ such as $Y =\diag(y_n, y_1,y_2,\dots,y_{n-1})$.
\end{remark}

\section{Structure-preserving model reduction for general skew-gradient systems}
\label{sec:nonlinear}

We consider general cases, such as the case that
$S(\by)$ nonlinearly depends on $\by$.
In such cases, the computational complexity of evaluating 
$S_r (\bz ) = V^\top S(V\bz) V$
may depend on $n$.
In this section, 
we show that the complexity can be reduced by utilizing the idea of the DEIM
while preserving the skew-gradient structure.

Let $S(t) := S(V\bz (t))$ to simplify the notation.
Using constant skew-symmetric matrices $U_j \in \bbR^{n\times n}$ ($j=1,\dots, m$),
and a time-dependent vector $\bc (t) \in \bbR^m$,
we approximate $S(t)$ by
\begin{align}
S(t) \approx \sum_{j=1}^m U_j c_j (t) .
\end{align}
This relation can be written as
\begin{align}
\bs (t) \approx U \bc (t),
\end{align}
where
$\bs (t) = \veco (S(t))$ and
 $U = [\veco (U_1), \veco (U_2), \dots, \veco (U_m)] \in \bbR^{n^2 \times m}$.
 
Following the discussion in Section~\ref{subsec:deim},
we require that $\bs (t)$ and $U\bc(t)$ are equal for $m$ variables
out of $n^2$ variables,
i.e. we require that
\begin{align}
P^\top \bs (t) = P^\top U \bc (t)
\end{align}
with $P = [\be_{\varrho_1},\dots,\be_{\varrho_m}]\in\bbR^{n^2\times m}$.
Let us assume that $P^\top U \in \bbR^{m\times m}$ is nonsingular.
Then, $\bc (t)$ is given by
\begin{align}
\bc (t) = (P^\top U)^{-1} P^\top \bs (t),
\end{align}
and thus $V^\top S(t) V$ can be approximated by
\begin{align} 
V^\top S(t)V
= \veco^{-1} \paren*{(V\otimes V )^\top \veco (S(t))}  
 \approx
 \veco^{-1} \paren*{(V\otimes V )^\top U \bc(t)} 
 =  
\veco^{-1} \paren*{
 \underbrace{(V\otimes V )^\top U (P^\top U)^{-1}}_{r^2 \times m}
 \underbrace{ P^\top \bs (t)}_{m\times 1}
} .  \label{ap:Sdeim}
\end{align}
Since $(V\otimes V )^\top U (P^\top U)^{-1}$ can be computed in the off-line stage,
the computational complexity of evaluating the most right hand side of \eqref{ap:Sdeim}
is independent of $n$
as long as the complexity for $P^\top \bs (t)$ is independent of $n$.

The matrices $U$ and $P$ are constructed as follows.
Let $S_1,\dots,S_s$ and $\bs_1=\veco(S_1),\dots,\bs_s =\veco (S_s)$
be snapshot data for $S(t)$.
Applying the SVD yields the POD basis vectors $\bu_1,\dots,\bu_m$,
and then the construction of $U$ and $P$ follows the DEIM procedure
summarized in Algorithm~\ref{algo:DEIM}.

We note that if the snapshot matrices $S_1,\dots S_s$ are skew-symmetric,
the approximation \eqref{ap:Sdeim} is also skew-symmetric.
This is readily seen from the following properties.
First,
since each $U_i$ is a linear combination of the snapshot matrices, all of which 
are skew-symmetric, the following property holds.
\begin{proposition}
If $S_1,\dots,S_s$ are skew-symmetric, $U_1,\dots,U_m$ are also skew-symmetric.
\end{proposition}
Next, since $S(t)$ is approximated by a linear combination of $U_i$, 
the following property holds.
\begin{proposition}
If $U_1,\dots,U_m \in \bbR^{n\times n}$ are skew-symmetric.
Then 
$\veco^{-1} \paren*{
(V\otimes V )^\top U (P^\top U)^{-1}
P^\top \bs (t)
}$
is also skew-symmetric.
\end{proposition}

\section{Numerical experiments}
\label{sec:num}

In this section, we check the performance of the proposed structure-preserving
model reduction methods.
The main aim of this section is to check the preservation of a first integral and 
the stability for the reduced-order system.
We employ the KdV equation and the modified KdV equation
as our toy problems.
For both equations, an appropriate finite difference semi-discretization yields 
skew-gradient systems.

All the computations are performed in a computation environment: 
3.5 GHz Intel Core i5, 8GB memory, OS X 10.13. 
We use MATLAB (R2015a). Nonlinear equations are solved by 
the matlab function $\mathsf{fsolve}$ with tolerance $10^{-16}$. 
The singular value decomposition is performed by $\mathsf{svd}$. 
Note that this function
computes all singular values and the corresponding singular vectors, 
which is not always necessary
for practical applications.
We employ this function just to observe the behaviour of the singular values for the test 
problems.

\subsection{KdV equation}

As an illustrative example, we consider the KdV equation
\begin{align} \label{eq:kdv}
y_t + 6yy_x + y_{xxx} = 0, \quad y(t_0,\cdot) = y_0, 
\quad x \in \bbT,
\end{align}
where $\bbT$ denotes the torus of length $L$.
This equation is completely integrable, and thus has infinitely many
conservation laws.
Among them, we here consider the $L^2$-norm preservation:
\begin{align} \label{eq:l2}
\frac{\rmd}{\rmd t} \calH [y] = 0, \quad \calH [y] = 
\int_\bbT \frac{y^2}{2}\,\rmd x .
\end{align}
We shall call this quantity and its discrete version the energy.
The KdV equation can be formulated as
\begin{align} \label{vs:kdv}
y_t = -\paren*{2\paren*{y\px + \px y} + \px^3} \fracdel{\calH}{y},
\end{align}
where the variational derivative of $\calH$ is given by $\delta \calH / \delta y = y$.

We here discretize the KdV equation \eqref{S:kdv} in space as follows.
Let $\by (t) := [ y_1(t),\dots , y_n(t)]^\top$,
where $\Delta x = L / n$ and $y_i (t)$ denotes the approximation to $ y(t, i\Delta x )$.
We use the central difference operators in a matrix form:
\begin{align*}
D_1 := \frac{1}{2\Delta x }
\begin{bmatrix}
0  & 1 & &  & -1 \\
-1 & 0  & 1 & &   \\
 & \ddots & \ddots & \ddots &  \\
 1 & & & -1 & 0
\end{bmatrix} \in \bbR^{n\times n}, \quad 
D_2 := \frac{1}{2\Delta x }
\begin{bmatrix}
-2  & 1 & &  & 1 \\
1 & -2  & 1 & &   \\
 & \ddots & \ddots & \ddots &  \\
 1 & & & 1 & -2
\end{bmatrix} \in \bbR^{n\times n}
\end{align*}
and $D_3:=D_1 D_2\in \bbR^{n\times n}$.
An appropriate discretization of
the variational form \eqref{vs:kdv}  yields the skew-gradient system
\begin{align}  \label{eq:sdkdv}
\dot{\by} = S(\by) \nabla_\by H(\by), \quad H(\by):= \frac{1}{2}\by^\top \by
\end{align}
with
\begin{align} \label{S:kdv}
S(\by) =-\paren*{2\paren*{YD_1 + D_1 Y} + D_3},
\end{align}
where $Y = \diag (\by)$.
Note that $\nabla H(\by) = \by$ is linear.
Since $S(\by )$ is of the form~\eqref{S:kdv},
the approach discussed in Section~\ref{sec:linear} is applicable.
For additional details on the spatial discretization based on variational structure,
see, e.g.~\cite{cgm12,f99,fm11} and references therein.

For the temporal discretization,
applying the standard mid-point rule to \eqref{eq:sdkdv}
yields an energy-preserving integrator:
for the solution to 
\begin{align}
\frac{\by_{n+1}-\by_{n}}{\Delta t} 
= S\paren*{\frac{\by_{n+1}+\by_{n}}{2}} 
\nabla_\by H \paren*{\frac{\by_{n+1}+\by_{n}}{2}},
\end{align}
where $\Delta t$ denotes the time stepsize and 
$\by_{n} \approx \by(n\Delta t)$ ($n=0,1,2,\dots$),
it follows that $H(\by_{n+1}) = H(\by_{n})$.
Note that any Runge--Kutta methods with the property
$b_i a_{ij} + b_j a_{ji} = b_i b_j$ ($i,j=1,\dots ,s$),
where $a_{ij}$ and $b_i$ are Runge--Kutta coefficients
and 
$s$ denotes the number of the stages,
preserve any linear and quadratic invariants~\cite{hlw06}.
The simplest example is the mid-point rule, which is the second order method.

\begin{remark}
Runge--Kutta methods cannot be energy-preserving in general.
For more general forms of the energy function,
applying the discrete gradient method yields an energy-preserving integrator~\cite{ch11,go96,ha10,mq99,mi16,qm08}.
\end{remark}

In the full order simulation,
we set $L=20$, $n=500$, $\Delta x = L/n = 0.04$.
The initial vale is set to $u(0,x) =2 (1.5^2)  \sech^2(3x/2)$.
The corresponding solution is a solitary wave 
if the spatial domain is $(-\infty,\infty)$, but 
in the bounded domain, the solution behaves almost periodically 
if $L$ is sufficiently large. 
We set $T = 3$ and $\Delta t =T/600$, which means we collect $601$ snapshot data:
$Y = [\by_0, \by_1, \dots, \by_{600}]\in\bbR^{500 \times 601}$
(note that $\nabla_\by H(\by) = \by$ in this case, cf.~\eqref{gmat}).

Fig.~\ref{fig:kdv1:svd} plots the singular values of the snapshot matrix $Y$.
A fast decay of the singular values indicates that 
a few modes can express the data with a good accuracy.
Fig.~\ref{fig:kdv1:EnergyError} shows the error growth of the energy $H(V\bz)$.
The energy is well preserved with considerable accuracy.
We plot the solution in Fig.~\ref{fig:kdv1:sol}.
Due to the structure-preservation, 
the numerical solution seems stable even for  small $r$.
We observe that the solution becomes smooth as $r$ gets large, and note that
the solution of the full-order model, which is not displayed, 
is almost identical to the result for $r=60$.
Global errors measured by the discrete version of the $L^2$-norm 
are plotted in Fig.~\ref{fig:kdv1:error},
where the solution to the reduced-order system is compared with that to 
the full-order system $\by^{\text{full}}$.
We observe that the error gets small as $r$ gets large.
When $r=60$,
the error remains small for $t>T=3$.

In this example, 
it should be noted that
the DEIM (or other techniques to reduce the complexity for the nonlinear term)
was not used, and
that only the standard POD matrix 
was used to find the reduced-order system whose vector field can be efficiently
evaluated.
The idea of inheriting the skew-gradient structure made this possible.

\input{figkdv1svd.tex}

\input{figkdv1Eerror.tex}

\input{figkdv1sol.tex}

\input{figkdv1error.tex}

\subsection{Modified KdV equation}
We next consider the modified  KdV (mKdV) equation
\begin{align} \label{eq:mkdv}
y_t + 6y^2y_x + y_{xxx} = 0, \quad y(t_0,\cdot) = y_0, 
\quad x \in \bbT.
\end{align}
The $L^2$-norm preservation \eqref{eq:l2} also holds for this equation,
which can be easily checked based on the variational structure
\begin{align} \label{vs:mkdv}
y_t = -\paren*{\frac{3}{2}\paren*{y^2\px + \px y^2} + \px^3} \fracdel{\calH}{y}.
\end{align}
As is the case with the previous subsection,
we discretize the variational form \eqref{vs:mkdv} as
\begin{align} 
\dot{\by} = S(\by) \nabla H(\by), \quad H(\by):= \frac{1}{2}\by^\top \by
\end{align}
with
\begin{align} \label{S:mkdv}
S(\by) =-\paren*{\frac{3}{2}\paren*{Y^2D_1 + D_1 Y^2} + D_3},
\end{align}
where $Y = \diag (\by)$.
In this case, since $S(\by)$ depends nonlinearly on $\by$,
we apply the approach presented in Section~\ref{sec:nonlinear}.

Note that the matrix \eqref{S:kdv} is sparse and only $4n$ entries are nonzero.
In our numerical experiments, instead of applying the SVD 
to the $n^2$-by-$s$ matrix, which is quite time-consuming,
we apply the SVD to the $4n$-by-$s$ matrix by simply ignoring the zero entries.

In the full-order simulation,
we set $L=10$, $n=500$, $\Delta x = L/n = 0.02$.
The initial vale is set to $u(0,x) =\sqrt{c}  \sech(\sqrt{c}x)$ with $c=4$.
The corresponding solution is a solitary wave 
in the unbounded spatial domain $(-\infty,\infty)$, but 
in the bounded domain, the solution behaves almost periodically 
if $L$ is sufficiently large. 
We set $T = 3$ and $\Delta t =T/750$, which means we collect $751$ snapshot data:
$Y = [\by_0, \by_1, \dots, \by_{750}]\in\bbR^{500 \times 751}$.

Fig.~\ref{fig:mkdv1:svd} plots the singular values of the snapshot matrix $Y$ and $[\veco (S(\by_0)), \veco (S(\by_1)),\dots,\veco (S(\by_{750}))]$.
A fast decay of the singular values indicates that a few modes can express the data with a good accuracy.
Fig.~\ref{fig:mkdv1:EnergyError} shows the error growth of the energy $H(V\bz)$.
The energy is well preserved with considerable accuracy.
We plot the solution in Fig.~\ref{fig:mkdv1:sol}.
The result with $r=50$ is almost identical to the solution of the full-order model, 
which is not displayed here.
But the solution differs substantially when $r=40$, which indicates
the number of basis matrices for the skew-symmetric matrix is 
sensitive to the qualitative behaviour.
Global errors measured by the discrete version of the $L^2$-norm 
are plotted in Fig.~\ref{fig:mkdv1:error}.
For $r=50$, the global error remains small even for $t>T=3$.

\input{figmkdv1svd.tex}
\input{figmkdv1Eerror.tex}
\input{figmkdv1sol.tex}
\input{figmkdv1error.tex}

\section{Concluding remarks}
\label{sec:conclusion}

In this paper, structure-preserving model reduction methods
for skew-gradient systems
of the form \eqref{sys1} have been studied.
We have shown that if $S(\by)$ has a specific structure, the previous approach 
proposed by Gong et al.~\cite{gw17} can efficiently reduce the size of the full-order system,
and also proposed a new approach for general cases that is based on the approach by Gong et al.~\cite{gw17}
and the discrete empirical interpolation method.
Since the reduced-order systems keep the skew-gradient structure and thus have the energy-conservation law,
energy-preserving integrators can be easily applied and one could expect the good qualitative behaviour.

\bibliographystyle{plain}
\bibliography{references}
\end{document}

%% file: figkdv1svd.tex
\begin{figure}[htbp]
\centering

\begin{tikzpicture}
\tikzstyle{every node}=[]
\begin{semilogyaxis}[width=6cm,
xmax=500,xmin=0,
xlabel={order of modes},ylabel={singlura values},
ylabel near ticks,
	]
\addplot[thick,
] table {
1	1.95E+02
2	1.82E+02
3	1.63E+02
4	1.44E+02
5	1.40E+02
6	1.34E+02
7	1.31E+02
8	1.20E+02
9	1.19E+02
10	1.05E+02
11	1.03E+02
12	8.90E+01
13	8.84E+01
14	7.43E+01
15	7.35E+01
16	6.10E+01
17	6.02E+01
18	4.91E+01
19	4.89E+01
20	3.93E+01
21	3.90E+01
22	3.11E+01
23	3.09E+01
24	2.44E+01
25	2.43E+01
26	1.90E+01
27	1.89E+01
28	1.48E+01
29	1.47E+01
30	1.14E+01
31	1.13E+01
32	8.75E+00
33	8.73E+00
34	6.71E+00
35	6.67E+00
36	5.12E+00
37	5.09E+00
38	3.88E+00
39	3.88E+00
40	2.95E+00
41	2.94E+00
42	2.24E+00
43	2.23E+00
44	1.69E+00
45	1.68E+00
46	1.27E+00
47	1.27E+00
48	9.58E-01
49	9.54E-01
50	7.20E-01
51	7.18E-01
52	5.40E-01
53	5.39E-01
54	4.05E-01
55	4.03E-01
56	3.03E-01
57	3.02E-01
58	2.26E-01
59	2.26E-01
60	1.69E-01
61	1.69E-01
62	1.26E-01
63	1.26E-01
64	9.39E-02
65	9.38E-02
66	7.00E-02
67	6.98E-02
68	5.21E-02
69	5.19E-02
70	3.87E-02
71	3.87E-02
72	2.88E-02
73	2.87E-02
74	2.14E-02
75	2.13E-02
76	1.59E-02
77	1.58E-02
78	1.18E-02
79	1.18E-02
80	8.73E-03
81	8.71E-03
82	6.47E-03
83	6.45E-03
84	4.79E-03
85	4.79E-03
86	3.55E-03
87	3.54E-03
88	2.63E-03
89	2.62E-03
90	1.94E-03
91	1.94E-03
92	1.44E-03
93	1.44E-03
94	1.06E-03
95	1.06E-03
96	7.86E-04
97	7.85E-04
98	5.81E-04
99	5.80E-04
100	4.29E-04
101	4.28E-04
102	3.17E-04
103	3.16E-04
104	2.34E-04
105	2.34E-04
106	1.73E-04
107	1.72E-04
108	1.28E-04
109	1.27E-04
110	9.42E-05
111	9.41E-05
112	6.95E-05
113	6.94E-05
114	5.13E-05
115	5.12E-05
116	3.79E-05
117	3.78E-05
118	2.79E-05
119	2.79E-05
120	2.06E-05
121	2.06E-05
122	1.52E-05
123	1.52E-05
124	1.12E-05
125	1.12E-05
126	8.28E-06
127	8.27E-06
128	6.11E-06
129	6.09E-06
130	4.50E-06
131	4.50E-06
132	3.32E-06
133	3.31E-06
134	2.45E-06
135	2.44E-06
136	1.81E-06
137	1.80E-06
138	1.33E-06
139	1.33E-06
140	9.82E-07
141	9.80E-07
142	7.25E-07
143	7.23E-07
144	5.34E-07
145	5.33E-07
146	3.94E-07
147	3.93E-07
148	2.91E-07
149	2.90E-07
150	2.14E-07
151	2.14E-07
152	1.58E-07
153	1.58E-07
154	1.17E-07
155	1.16E-07
156	8.60E-08
157	8.58E-08
158	6.34E-08
159	6.34E-08
160	4.68E-08
161	4.67E-08
162	3.45E-08
163	3.45E-08
164	2.55E-08
165	2.54E-08
166	1.88E-08
167	1.87E-08
168	1.39E-08
169	1.38E-08
170	1.02E-08
171	1.02E-08
172	7.54E-09
173	7.53E-09
174	5.58E-09
175	5.56E-09
176	4.11E-09
177	4.10E-09
178	3.04E-09
179	3.02E-09
180	2.24E-09
181	2.23E-09
182	1.65E-09
183	1.65E-09
184	1.22E-09
185	1.22E-09
186	9.01E-10
187	8.97E-10
188	6.63E-10
189	6.62E-10
190	4.90E-10
191	4.89E-10
192	3.61E-10
193	3.61E-10
194	2.67E-10
195	2.66E-10
196	1.97E-10
197	1.96E-10
198	1.45E-10
199	1.45E-10
200	1.07E-10
201	1.07E-10
202	7.91E-11
203	7.90E-11
204	5.85E-11
205	5.83E-11
206	4.31E-11
207	4.30E-11
208	3.19E-11
209	3.18E-11
210	2.36E-11
211	2.34E-11
212	1.74E-11
213	1.73E-11
214	1.28E-11
215	1.28E-11
216	9.49E-12
217	9.48E-12
218	6.98E-12
219	6.97E-12
220	5.20E-12
221	5.19E-12
222	3.81E-12
223	3.80E-12
224	2.85E-12
225	2.83E-12
226	2.12E-12
227	2.12E-12
228	1.68E-12
229	1.65E-12
230	1.38E-12
231	1.38E-12
232	1.12E-12
233	1.10E-12
234	8.70E-13
235	8.44E-13
236	6.59E-13
237	6.34E-13
238	5.32E-13
239	5.26E-13
240	4.13E-13
241	4.01E-13
242	3.10E-13
243	3.08E-13
244	2.80E-13
245	2.79E-13
246	2.57E-13
247	2.54E-13
248	2.37E-13
249	2.36E-13
250	2.34E-13
251	2.31E-13
252	2.23E-13
253	2.21E-13
254	2.15E-13
255	2.11E-13
256	2.06E-13
257	2.00E-13
258	1.96E-13
259	1.89E-13
260	1.88E-13
261	1.83E-13
262	1.80E-13
263	1.77E-13
264	1.69E-13
265	1.67E-13
266	1.65E-13
267	1.61E-13
268	1.60E-13
269	1.58E-13
270	1.53E-13
271	1.51E-13
272	1.43E-13
273	1.40E-13
274	1.34E-13
275	1.33E-13
276	1.30E-13
277	1.28E-13
278	1.25E-13
279	1.22E-13
280	1.11E-13
281	1.06E-13
282	9.71E-14
283	9.62E-14
284	9.38E-14
285	9.24E-14
286	8.87E-14
287	7.66E-14
288	7.02E-14
289	6.86E-14
290	6.08E-14
291	5.92E-14
292	5.57E-14
293	5.31E-14
294	5.10E-14
295	4.84E-14
296	4.59E-14
297	4.52E-14
298	4.19E-14
299	4.07E-14
300	3.92E-14
301	3.76E-14
302	3.63E-14
303	3.58E-14
304	3.42E-14
305	3.16E-14
306	3.08E-14
307	2.98E-14
308	2.87E-14
309	2.74E-14
310	2.69E-14
311	2.62E-14
312	2.47E-14
313	2.44E-14
314	2.37E-14
315	2.17E-14
316	2.12E-14
317	1.97E-14
318	1.84E-14
319	1.74E-14
320	1.66E-14
321	1.46E-14
322	1.36E-14
323	1.36E-14
324	1.20E-14
325	1.17E-14
326	1.17E-14
327	1.17E-14
328	1.17E-14
329	1.17E-14
330	1.17E-14
331	1.17E-14
332	1.17E-14
333	1.17E-14
334	1.17E-14
335	1.17E-14
336	1.17E-14
337	1.17E-14
338	1.17E-14
339	1.17E-14
340	1.17E-14
341	1.17E-14
342	1.17E-14
343	1.17E-14
344	1.17E-14
345	1.17E-14
346	1.17E-14
347	1.17E-14
348	1.17E-14
349	1.17E-14
350	1.17E-14
351	1.17E-14
352	1.17E-14
353	1.17E-14
354	1.17E-14
355	1.17E-14
356	1.17E-14
357	1.17E-14
358	1.17E-14
359	1.17E-14
360	1.17E-14
361	1.17E-14
362	1.17E-14
363	1.17E-14
364	1.17E-14
365	1.17E-14
366	1.17E-14
367	1.17E-14
368	1.17E-14
369	1.17E-14
370	1.17E-14
371	1.17E-14
372	1.17E-14
373	1.17E-14
374	1.17E-14
375	1.17E-14
376	1.17E-14
377	1.17E-14
378	1.17E-14
379	1.17E-14
380	1.17E-14
381	1.17E-14
382	1.17E-14
383	1.17E-14
384	1.17E-14
385	1.17E-14
386	1.17E-14
387	1.17E-14
388	1.17E-14
389	1.17E-14
390	1.17E-14
391	1.17E-14
392	1.17E-14
393	1.17E-14
394	1.17E-14
395	1.17E-14
396	1.17E-14
397	1.17E-14
398	1.17E-14
399	1.17E-14
400	1.17E-14
401	1.17E-14
402	1.17E-14
403	1.17E-14
404	1.17E-14
405	1.17E-14
406	1.17E-14
407	1.17E-14
408	1.17E-14
409	1.17E-14
410	1.17E-14
411	1.17E-14
412	1.17E-14
413	1.17E-14
414	1.17E-14
415	1.17E-14
416	1.17E-14
417	1.17E-14
418	1.17E-14
419	1.17E-14
420	1.17E-14
421	1.17E-14
422	1.17E-14
423	1.17E-14
424	1.17E-14
425	1.17E-14
426	1.17E-14
427	1.17E-14
428	1.17E-14
429	1.17E-14
430	1.17E-14
431	1.17E-14
432	1.17E-14
433	1.17E-14
434	1.17E-14
435	1.17E-14
436	1.17E-14
437	1.17E-14
438	1.17E-14
439	1.17E-14
440	1.17E-14
441	1.17E-14
442	1.17E-14
443	1.17E-14
444	1.17E-14
445	1.17E-14
446	1.17E-14
447	1.17E-14
448	1.17E-14
449	1.17E-14
450	1.17E-14
451	1.17E-14
452	1.17E-14
453	1.17E-14
454	1.17E-14
455	1.17E-14
456	1.17E-14
457	1.17E-14
458	1.17E-14
459	1.17E-14
460	1.17E-14
461	1.17E-14
462	1.17E-14
463	1.17E-14
464	1.17E-14
465	1.17E-14
466	1.17E-14
467	1.17E-14
468	1.17E-14
469	1.17E-14
470	1.17E-14
471	1.17E-14
472	1.17E-14
473	1.17E-14
474	1.17E-14
475	1.17E-14
476	1.17E-14
477	1.17E-14
478	1.17E-14
479	1.17E-14
480	1.17E-14
481	1.17E-14
482	1.17E-14
483	1.17E-14
484	1.17E-14
485	1.17E-14
486	1.17E-14
487	1.17E-14
488	1.17E-14
489	1.17E-14
490	1.17E-14
491	1.17E-14
492	1.17E-14
493	1.09E-14
494	9.82E-15
495	8.24E-15
496	6.81E-15
497	5.75E-15
498	4.04E-15
499	2.70E-15
500	1.49E-15
};
\end{semilogyaxis}
\end{tikzpicture}

\caption{Singular values corresponding to the POD modes for the matrix $Y$ for the KdV equation.
}
\label{fig:kdv1:svd}
\end{figure}
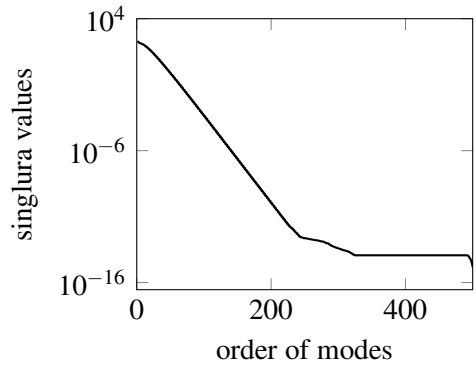

%% file: figkdv1Eerror.tex
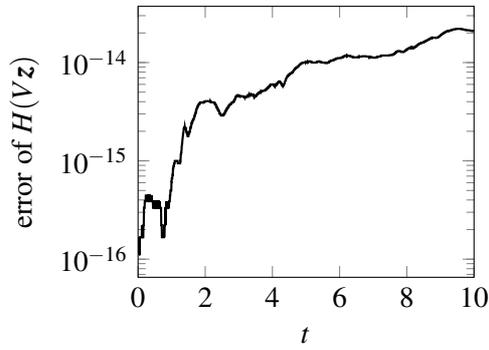
\begin{figure}[htbp]
\centering

\begin{tikzpicture}
\tikzstyle{every node}=[]
\begin{semilogyaxis}[width=6cm,
xmax=10,xmin=0,
xlabel={$t$},ylabel={error of $H(V\bz)$},
ylabel near ticks,
	]
\addplot[thick,
] table {
0	0
0.005	0
0.01	0
0.015	0
0.02	0
0.025	1.11E-16
0.03	1.11E-16
0.035	1.11E-16
0.04	1.67E-16
0.045	1.11E-16
0.05	1.67E-16
0.055	1.67E-16
0.06	1.67E-16
0.065	1.67E-16
0.07	1.67E-16
0.075	1.67E-16
0.08	1.67E-16
0.085	1.67E-16
0.09	1.67E-16
0.095	1.67E-16
0.1	1.67E-16
0.105	1.67E-16
0.11	1.67E-16
0.115	2.22E-16
0.12	1.67E-16
0.125	1.67E-16
0.13	1.67E-16
0.135	1.67E-16
0.14	1.67E-16
0.145	1.67E-16
0.15	1.67E-16
0.155	1.67E-16
0.16	2.22E-16
0.165	2.22E-16
0.17	2.22E-16
0.175	3.33E-16
0.18	3.33E-16
0.185	3.33E-16
0.19	3.33E-16
0.195	3.33E-16
0.2	3.89E-16
0.205	3.89E-16
0.21	3.89E-16
0.215	3.89E-16
0.22	3.89E-16
0.225	3.89E-16
0.23	4.44E-16
0.235	4.44E-16
0.24	4.44E-16
0.245	4.44E-16
0.25	4.44E-16
0.255	4.44E-16
0.26	4.44E-16
0.265	4.44E-16
0.27	4.44E-16
0.275	4.44E-16
0.28	4.44E-16
0.285	3.89E-16
0.29	3.89E-16
0.295	3.89E-16
0.3	3.89E-16
0.305	3.89E-16
0.31	4.44E-16
0.315	4.44E-16
0.32	3.89E-16
0.325	3.89E-16
0.33	4.44E-16
0.335	4.44E-16
0.34	4.44E-16
0.345	4.44E-16
0.35	3.89E-16
0.355	4.44E-16
0.36	4.44E-16
0.365	3.89E-16
0.37	3.89E-16
0.375	3.89E-16
0.38	3.89E-16
0.385	3.89E-16
0.39	3.89E-16
0.395	3.89E-16
0.4	3.89E-16
0.405	3.89E-16
0.41	4.44E-16
0.415	3.89E-16
0.42	3.89E-16
0.425	3.89E-16
0.43	3.89E-16
0.435	3.89E-16
0.44	3.89E-16
0.445	3.33E-16
0.45	3.33E-16
0.455	3.33E-16
0.46	3.89E-16
0.465	3.89E-16
0.47	3.89E-16
0.475	3.89E-16
0.48	3.89E-16
0.485	3.89E-16
0.49	3.89E-16
0.495	3.89E-16
0.5	3.33E-16
0.505	3.33E-16
0.51	3.89E-16
0.515	3.89E-16
0.52	3.89E-16
0.525	3.33E-16
0.53	3.33E-16
0.535	3.33E-16
0.54	3.33E-16
0.545	3.33E-16
0.55	3.89E-16
0.555	3.89E-16
0.56	3.89E-16
0.565	3.33E-16
0.57	3.89E-16
0.575	3.33E-16
0.58	3.33E-16
0.585	3.89E-16
0.59	3.89E-16
0.595	3.89E-16
0.6	3.89E-16
0.605	3.89E-16
0.61	3.33E-16
0.615	3.89E-16
0.62	3.89E-16
0.625	3.89E-16
0.63	3.33E-16
0.635	3.33E-16
0.64	3.89E-16
0.645	3.33E-16
0.65	3.33E-16
0.655	3.33E-16
0.66	3.33E-16
0.665	3.33E-16
0.67	3.33E-16
0.675	3.33E-16
0.68	3.33E-16
0.685	3.33E-16
0.69	2.22E-16
0.695	2.22E-16
0.7	2.22E-16
0.705	2.22E-16
0.71	1.67E-16
0.715	1.67E-16
0.72	1.67E-16
0.725	1.67E-16
0.73	1.67E-16
0.735	1.67E-16
0.74	1.67E-16
0.745	1.67E-16
0.75	1.67E-16
0.755	1.67E-16
0.76	1.67E-16
0.765	1.67E-16
0.77	1.67E-16
0.775	2.22E-16
0.78	1.67E-16
0.785	2.22E-16
0.79	2.22E-16
0.795	2.22E-16
0.8	2.22E-16
0.805	2.22E-16
0.81	3.33E-16
0.815	2.22E-16
0.82	3.33E-16
0.825	3.33E-16
0.83	3.33E-16
0.835	3.33E-16
0.84	3.89E-16
0.845	3.33E-16
0.85	3.33E-16
0.855	3.33E-16
0.86	3.33E-16
0.865	3.33E-16
0.87	3.33E-16
0.875	3.33E-16
0.88	3.33E-16
0.885	3.33E-16
0.89	3.89E-16
0.895	3.33E-16
0.9	3.33E-16
0.905	3.33E-16
0.91	3.33E-16
0.915	3.89E-16
0.92	3.33E-16
0.925	3.89E-16
0.93	3.89E-16
0.935	3.89E-16
0.94	3.89E-16
0.945	3.89E-16
0.95	4.44E-16
0.955	4.44E-16
0.96	4.44E-16
0.965	5.00E-16
0.97	5.00E-16
0.975	5.00E-16
0.98	5.00E-16
0.985	5.00E-16
0.99	5.00E-16
0.995	6.11E-16
1	6.11E-16
1.005	6.66E-16
1.01	6.66E-16
1.015	6.66E-16
1.02	6.66E-16
1.025	7.22E-16
1.03	7.22E-16
1.035	7.22E-16
1.04	8.33E-16
1.045	8.33E-16
1.05	8.33E-16
1.055	8.88E-16
1.06	8.88E-16
1.065	8.88E-16
1.07	8.88E-16
1.075	9.44E-16
1.08	9.44E-16
1.085	9.44E-16
1.09	9.44E-16
1.095	9.99E-16
1.1	9.99E-16
1.105	9.99E-16
1.11	9.99E-16
1.115	9.99E-16
1.12	9.99E-16
1.125	9.99E-16
1.13	9.99E-16
1.135	9.99E-16
1.14	9.99E-16
1.145	9.99E-16
1.15	9.99E-16
1.155	9.99E-16
1.16	9.99E-16
1.165	9.99E-16
1.17	9.99E-16
1.175	9.99E-16
1.18	9.99E-16
1.185	9.44E-16
1.19	9.44E-16
1.195	9.44E-16
1.2	9.44E-16
1.205	9.44E-16
1.21	9.44E-16
1.215	9.44E-16
1.22	9.44E-16
1.225	9.44E-16
1.23	9.44E-16
1.235	9.44E-16
1.24	9.44E-16
1.245	9.44E-16
1.25	9.44E-16
1.255	9.44E-16
1.26	9.99E-16
1.265	9.99E-16
1.27	1.11E-15
1.275	1.17E-15
1.28	1.17E-15
1.285	1.22E-15
1.29	1.22E-15
1.295	1.28E-15
1.3	1.39E-15
1.305	1.44E-15
1.31	1.50E-15
1.315	1.50E-15
1.32	1.61E-15
1.325	1.67E-15
1.33	1.72E-15
1.335	1.78E-15
1.34	1.89E-15
1.345	1.89E-15
1.35	1.94E-15
1.355	2.00E-15
1.36	2.11E-15
1.365	2.11E-15
1.37	2.16E-15
1.375	2.16E-15
1.38	2.16E-15
1.385	2.16E-15
1.39	2.22E-15
1.395	2.16E-15
1.4	2.16E-15
1.405	2.16E-15
1.41	2.16E-15
1.415	2.16E-15
1.42	2.16E-15
1.425	2.16E-15
1.43	2.11E-15
1.435	2.00E-15
1.44	2.00E-15
1.445	1.94E-15
1.45	1.94E-15
1.455	1.94E-15
1.46	1.89E-15
1.465	1.89E-15
1.47	1.89E-15
1.475	1.78E-15
1.48	1.78E-15
1.485	1.78E-15
1.49	1.78E-15
1.495	1.78E-15
1.5	1.78E-15
1.505	1.89E-15
1.51	1.89E-15
1.515	1.89E-15
1.52	1.94E-15
1.525	1.94E-15
1.53	1.94E-15
1.535	1.94E-15
1.54	2.00E-15
1.545	2.00E-15
1.55	2.11E-15
1.555	2.11E-15
1.56	2.16E-15
1.565	2.16E-15
1.57	2.22E-15
1.575	2.22E-15
1.58	2.28E-15
1.585	2.28E-15
1.59	2.39E-15
1.595	2.39E-15
1.6	2.44E-15
1.605	2.44E-15
1.61	2.50E-15
1.615	2.50E-15
1.62	2.50E-15
1.625	2.50E-15
1.63	2.61E-15
1.635	2.61E-15
1.64	2.61E-15
1.645	2.61E-15
1.65	2.66E-15
1.655	2.66E-15
1.66	2.66E-15
1.665	2.72E-15
1.67	2.78E-15
1.675	2.78E-15
1.68	2.78E-15
1.685	2.89E-15
1.69	2.89E-15
1.695	2.94E-15
1.7	2.94E-15
1.705	3.00E-15
1.71	3.05E-15
1.715	3.05E-15
1.72	3.16E-15
1.725	3.16E-15
1.73	3.22E-15
1.735	3.28E-15
1.74	3.28E-15
1.745	3.39E-15
1.75	3.44E-15
1.755	3.44E-15
1.76	3.50E-15
1.765	3.50E-15
1.77	3.55E-15
1.775	3.55E-15
1.78	3.55E-15
1.785	3.66E-15
1.79	3.66E-15
1.795	3.72E-15
1.8	3.72E-15
1.805	3.72E-15
1.81	3.72E-15
1.815	3.77E-15
1.82	3.77E-15
1.825	3.77E-15
1.83	3.77E-15
1.835	3.89E-15
1.84	3.89E-15
1.845	3.89E-15
1.85	3.89E-15
1.855	3.94E-15
1.86	3.94E-15
1.865	3.94E-15
1.87	3.94E-15
1.875	3.94E-15
1.88	3.94E-15
1.885	3.94E-15
1.89	3.94E-15
1.895	3.94E-15
1.9	3.94E-15
1.905	3.94E-15
1.91	3.94E-15
1.915	3.94E-15
1.92	3.94E-15
1.925	3.94E-15
1.93	3.89E-15
1.935	3.89E-15
1.94	3.94E-15
1.945	3.94E-15
1.95	3.94E-15
1.955	3.94E-15
1.96	3.94E-15
1.965	3.94E-15
1.97	3.94E-15
1.975	4.00E-15
1.98	4.00E-15
1.985	4.00E-15
1.99	4.00E-15
1.995	4.00E-15
2	4.00E-15
2.005	4.05E-15
2.01	4.05E-15
2.015	4.05E-15
2.02	4.05E-15
2.025	4.05E-15
2.03	4.05E-15
2.035	4.05E-15
2.04	4.05E-15
2.045	4.05E-15
2.05	4.05E-15
2.055	4.05E-15
2.06	4.05E-15
2.065	4.05E-15
2.07	4.05E-15
2.075	4.05E-15
2.08	4.05E-15
2.085	4.05E-15
2.09	4.05E-15
2.095	4.05E-15
2.1	4.05E-15
2.105	4.05E-15
2.11	4.05E-15
2.115	4.05E-15
2.12	4.05E-15
2.125	4.05E-15
2.13	4.05E-15
2.135	4.00E-15
2.14	4.00E-15
2.145	4.05E-15
2.15	4.00E-15
2.155	4.00E-15
2.16	4.00E-15
2.165	4.00E-15
2.17	4.00E-15
2.175	4.00E-15
2.18	3.94E-15
2.185	3.94E-15
2.19	3.94E-15
2.195	3.94E-15
2.2	3.94E-15
2.205	3.94E-15
2.21	3.94E-15
2.215	3.94E-15
2.22	3.94E-15
2.225	3.94E-15
2.23	3.94E-15
2.235	3.94E-15
2.24	3.94E-15
2.245	3.94E-15
2.25	3.94E-15
2.255	3.94E-15
2.26	3.94E-15
2.265	3.94E-15
2.27	3.89E-15
2.275	3.89E-15
2.28	3.89E-15
2.285	3.89E-15
2.29	3.77E-15
2.295	3.77E-15
2.3	3.77E-15
2.305	3.77E-15
2.31	3.72E-15
2.315	3.72E-15
2.32	3.66E-15
2.325	3.66E-15
2.33	3.66E-15
2.335	3.55E-15
2.34	3.55E-15
2.345	3.50E-15
2.35	3.50E-15
2.355	3.44E-15
2.36	3.44E-15
2.365	3.39E-15
2.37	3.39E-15
2.375	3.39E-15
2.38	3.28E-15
2.385	3.28E-15
2.39	3.28E-15
2.395	3.22E-15
2.4	3.22E-15
2.405	3.22E-15
2.41	3.22E-15
2.415	3.22E-15
2.42	3.16E-15
2.425	3.05E-15
2.43	3.05E-15
2.435	3.05E-15
2.44	3.00E-15
2.445	3.00E-15
2.45	3.00E-15
2.455	2.94E-15
2.46	2.94E-15
2.465	2.94E-15
2.47	2.94E-15
2.475	2.94E-15
2.48	2.89E-15
2.485	2.89E-15
2.49	2.89E-15
2.495	2.89E-15
2.5	2.89E-15
2.505	2.89E-15
2.51	2.89E-15
2.515	2.89E-15
2.52	2.89E-15
2.525	2.89E-15
2.53	2.89E-15
2.535	2.94E-15
2.54	2.94E-15
2.545	2.94E-15
2.55	2.94E-15
2.555	2.94E-15
2.56	3.00E-15
2.565	3.00E-15
2.57	3.00E-15
2.575	3.05E-15
2.58	3.05E-15
2.585	3.05E-15
2.59	3.16E-15
2.595	3.16E-15
2.6	3.16E-15
2.605	3.16E-15
2.61	3.22E-15
2.615	3.22E-15
2.62	3.28E-15
2.625	3.28E-15
2.63	3.28E-15
2.635	3.28E-15
2.64	3.39E-15
2.645	3.39E-15
2.65	3.39E-15
2.655	3.39E-15
2.66	3.44E-15
2.665	3.44E-15
2.67	3.44E-15
2.675	3.50E-15
2.68	3.50E-15
2.685	3.50E-15
2.69	3.55E-15
2.695	3.55E-15
2.7	3.55E-15
2.705	3.66E-15
2.71	3.66E-15
2.715	3.66E-15
2.72	3.72E-15
2.725	3.72E-15
2.73	3.72E-15
2.735	3.72E-15
2.74	3.72E-15
2.745	3.77E-15
2.75	3.77E-15
2.755	3.77E-15
2.76	3.77E-15
2.765	3.77E-15
2.77	3.77E-15
2.775	3.89E-15
2.78	3.89E-15
2.785	3.89E-15
2.79	3.89E-15
2.795	3.89E-15
2.8	3.89E-15
2.805	3.94E-15
2.81	3.94E-15
2.815	3.94E-15
2.82	3.94E-15
2.825	3.94E-15
2.83	3.94E-15
2.835	3.94E-15
2.84	4.00E-15
2.845	4.00E-15
2.85	4.05E-15
2.855	4.05E-15
2.86	4.05E-15
2.865	4.16E-15
2.87	4.16E-15
2.875	4.16E-15
2.88	4.22E-15
2.885	4.22E-15
2.89	4.27E-15
2.895	4.39E-15
2.9	4.39E-15
2.905	4.44E-15
2.91	4.44E-15
2.915	4.50E-15
2.92	4.50E-15
2.925	4.55E-15
2.93	4.55E-15
2.935	4.55E-15
2.94	4.55E-15
2.945	4.55E-15
2.95	4.66E-15
2.955	4.66E-15
2.96	4.66E-15
2.965	4.66E-15
2.97	4.66E-15
2.975	4.66E-15
2.98	4.66E-15
2.985	4.66E-15
2.99	4.66E-15
2.995	4.55E-15
3	4.55E-15
3.005	4.55E-15
3.01	4.55E-15
3.015	4.55E-15
3.02	4.55E-15
3.025	4.50E-15
3.03	4.50E-15
3.035	4.55E-15
3.04	4.55E-15
3.045	4.50E-15
3.05	4.50E-15
3.055	4.50E-15
3.06	4.55E-15
3.065	4.55E-15
3.07	4.55E-15
3.075	4.55E-15
3.08	4.55E-15
3.085	4.55E-15
3.09	4.55E-15
3.095	4.55E-15
3.1	4.55E-15
3.105	4.55E-15
3.11	4.55E-15
3.115	4.55E-15
3.12	4.55E-15
3.125	4.50E-15
3.13	4.50E-15
3.135	4.50E-15
3.14	4.50E-15
3.145	4.50E-15
3.15	4.50E-15
3.155	4.44E-15
3.16	4.44E-15
3.165	4.44E-15
3.17	4.44E-15
3.175	4.44E-15
3.18	4.39E-15
3.185	4.44E-15
3.19	4.44E-15
3.195	4.44E-15
3.2	4.44E-15
3.205	4.44E-15
3.21	4.44E-15
3.215	4.44E-15
3.22	4.44E-15
3.225	4.44E-15
3.23	4.44E-15
3.235	4.44E-15
3.24	4.50E-15
3.245	4.50E-15
3.25	4.55E-15
3.255	4.55E-15
3.26	4.55E-15
3.265	4.55E-15
3.27	4.66E-15
3.275	4.66E-15
3.28	4.66E-15
3.285	4.66E-15
3.29	4.66E-15
3.295	4.66E-15
3.3	4.66E-15
3.305	4.72E-15
3.31	4.66E-15
3.315	4.72E-15
3.32	4.72E-15
3.325	4.72E-15
3.33	4.72E-15
3.335	4.72E-15
3.34	4.72E-15
3.345	4.72E-15
3.35	4.77E-15
3.355	4.77E-15
3.36	4.72E-15
3.365	4.77E-15
3.37	4.77E-15
3.375	4.77E-15
3.38	4.72E-15
3.385	4.72E-15
3.39	4.72E-15
3.395	4.72E-15
3.4	4.72E-15
3.405	4.72E-15
3.41	4.72E-15
3.415	4.72E-15
3.42	4.72E-15
3.425	4.66E-15
3.43	4.66E-15
3.435	4.55E-15
3.44	4.55E-15
3.445	4.55E-15
3.45	4.50E-15
3.455	4.55E-15
3.46	4.50E-15
3.465	4.55E-15
3.47	4.55E-15
3.475	4.50E-15
3.48	4.50E-15
3.485	4.50E-15
3.49	4.50E-15
3.495	4.50E-15
3.5	4.55E-15
3.505	4.55E-15
3.51	4.55E-15
3.515	4.55E-15
3.52	4.66E-15
3.525	4.66E-15
3.53	4.66E-15
3.535	4.66E-15
3.54	4.72E-15
3.545	4.72E-15
3.55	4.72E-15
3.555	4.77E-15
3.56	4.83E-15
3.565	4.83E-15
3.57	4.83E-15
3.575	4.83E-15
3.58	4.83E-15
3.585	4.94E-15
3.59	4.94E-15
3.595	4.94E-15
3.6	4.94E-15
3.605	4.94E-15
3.61	4.94E-15
3.615	4.94E-15
3.62	5.00E-15
3.625	5.00E-15
3.63	5.00E-15
3.635	5.00E-15
3.64	5.00E-15
3.645	5.05E-15
3.65	5.00E-15
3.655	5.05E-15
3.66	5.05E-15
3.665	5.05E-15
3.67	5.05E-15
3.675	5.05E-15
3.68	5.05E-15
3.685	5.05E-15
3.69	5.16E-15
3.695	5.16E-15
3.7	5.16E-15
3.705	5.16E-15
3.71	5.16E-15
3.715	5.16E-15
3.72	5.16E-15
3.725	5.16E-15
3.73	5.16E-15
3.735	5.16E-15
3.74	5.16E-15
3.745	5.16E-15
3.75	5.16E-15
3.755	5.16E-15
3.76	5.16E-15
3.765	5.22E-15
3.77	5.22E-15
3.775	5.22E-15
3.78	5.22E-15
3.785	5.22E-15
3.79	5.27E-15
3.795	5.27E-15
3.8	5.27E-15
3.805	5.27E-15
3.81	5.33E-15
3.815	5.33E-15
3.82	5.33E-15
3.825	5.33E-15
3.83	5.44E-15
3.835	5.44E-15
3.84	5.44E-15
3.845	5.50E-15
3.85	5.50E-15
3.855	5.50E-15
3.86	5.50E-15
3.865	5.55E-15
3.87	5.55E-15
3.875	5.55E-15
3.88	5.55E-15
3.885	5.66E-15
3.89	5.66E-15
3.895	5.72E-15
3.9	5.72E-15
3.905	5.72E-15
3.91	5.77E-15
3.915	5.77E-15
3.92	5.83E-15
3.925	5.83E-15
3.93	5.83E-15
3.935	5.94E-15
3.94	5.94E-15
3.945	5.94E-15
3.95	5.94E-15
3.955	6.00E-15
3.96	6.00E-15
3.965	6.00E-15
3.97	6.00E-15
3.975	6.00E-15
3.98	6.00E-15
3.985	6.00E-15
3.99	6.00E-15
3.995	6.00E-15
4	6.00E-15
4.005	6.00E-15
4.01	6.00E-15
4.015	5.94E-15
4.02	5.94E-15
4.025	5.94E-15
4.03	5.94E-15
4.035	5.94E-15
4.04	5.83E-15
4.045	5.83E-15
4.05	5.83E-15
4.055	5.83E-15
4.06	5.83E-15
4.065	5.83E-15
4.07	5.77E-15
4.075	5.83E-15
4.08	5.77E-15
4.085	5.83E-15
4.09	5.83E-15
4.095	5.83E-15
4.1	5.83E-15
4.105	5.83E-15
4.11	5.94E-15
4.115	5.94E-15
4.12	5.94E-15
4.125	6.00E-15
4.13	6.00E-15
4.135	6.05E-15
4.14	6.05E-15
4.145	6.16E-15
4.15	6.16E-15
4.155	6.22E-15
4.16	6.22E-15
4.165	6.27E-15
4.17	6.27E-15
4.175	6.33E-15
4.18	6.33E-15
4.185	6.33E-15
4.19	6.33E-15
4.195	6.44E-15
4.2	6.44E-15
4.205	6.44E-15
4.21	6.44E-15
4.215	6.44E-15
4.22	6.44E-15
4.225	6.33E-15
4.23	6.33E-15
4.235	6.27E-15
4.24	6.27E-15
4.245	6.27E-15
4.25	6.22E-15
4.255	6.22E-15
4.26	6.16E-15
4.265	6.05E-15
4.27	6.05E-15
4.275	6.00E-15
4.28	5.94E-15
4.285	5.94E-15
4.29	5.83E-15
4.295	5.83E-15
4.3	5.83E-15
4.305	5.83E-15
4.31	5.77E-15
4.315	5.83E-15
4.32	5.83E-15
4.325	5.83E-15
4.33	5.94E-15
4.335	5.94E-15
4.34	5.94E-15
4.345	6.00E-15
4.35	6.05E-15
4.355	6.16E-15
4.36	6.22E-15
4.365	6.22E-15
4.37	6.33E-15
4.375	6.44E-15
4.38	6.49E-15
4.385	6.55E-15
4.39	6.61E-15
4.395	6.61E-15
4.4	6.77E-15
4.405	6.77E-15
4.41	6.83E-15
4.415	6.94E-15
4.42	6.99E-15
4.425	7.05E-15
4.43	7.05E-15
4.435	7.11E-15
4.44	7.22E-15
4.445	7.22E-15
4.45	7.27E-15
4.455	7.27E-15
4.46	7.33E-15
4.465	7.33E-15
4.47	7.44E-15
4.475	7.44E-15
4.48	7.49E-15
4.485	7.49E-15
4.49	7.55E-15
4.495	7.55E-15
4.5	7.55E-15
4.505	7.61E-15
4.51	7.61E-15
4.515	7.61E-15
4.52	7.72E-15
4.525	7.72E-15
4.53	7.72E-15
4.535	7.72E-15
4.54	7.77E-15
4.545	7.77E-15
4.55	7.77E-15
4.555	7.77E-15
4.56	7.83E-15
4.565	7.83E-15
4.57	7.83E-15
4.575	7.83E-15
4.58	7.83E-15
4.585	7.94E-15
4.59	7.94E-15
4.595	7.99E-15
4.6	7.99E-15
4.605	8.05E-15
4.61	8.05E-15
4.615	8.10E-15
4.62	8.22E-15
4.625	8.22E-15
4.63	8.27E-15
4.635	8.33E-15
4.64	8.33E-15
4.645	8.38E-15
4.65	8.49E-15
4.655	8.55E-15
4.66	8.55E-15
4.665	8.60E-15
4.67	8.60E-15
4.675	8.72E-15
4.68	8.72E-15
4.685	8.72E-15
4.69	8.77E-15
4.695	8.77E-15
4.7	8.83E-15
4.705	8.83E-15
4.71	8.83E-15
4.715	8.88E-15
4.72	8.88E-15
4.725	8.88E-15
4.73	8.99E-15
4.735	8.99E-15
4.74	8.99E-15
4.745	9.05E-15
4.75	9.05E-15
4.755	9.05E-15
4.76	9.10E-15
4.765	9.10E-15
4.77	9.21E-15
4.775	9.21E-15
4.78	9.27E-15
4.785	9.33E-15
4.79	9.33E-15
4.795	9.38E-15
4.8	9.38E-15
4.805	9.49E-15
4.81	9.55E-15
4.815	9.55E-15
4.82	9.55E-15
4.825	9.60E-15
4.83	9.60E-15
4.835	9.71E-15
4.84	9.71E-15
4.845	9.77E-15
4.85	9.77E-15
4.855	9.83E-15
4.86	9.83E-15
4.865	9.83E-15
4.87	9.83E-15
4.875	9.88E-15
4.88	9.88E-15
4.885	9.88E-15
4.89	9.88E-15
4.895	9.88E-15
4.9	9.88E-15
4.905	9.99E-15
4.91	9.99E-15
4.915	9.99E-15
4.92	1.00E-14
4.925	1.00E-14
4.93	1.00E-14
4.935	1.00E-14
4.94	1.00E-14
4.945	1.00E-14
4.95	1.00E-14
4.955	1.00E-14
4.96	1.00E-14
4.965	1.00E-14
4.97	1.01E-14
4.975	1.01E-14
4.98	1.01E-14
4.985	1.00E-14
4.99	1.01E-14
4.995	1.01E-14
5	1.01E-14
5.005	1.01E-14
5.01	1.01E-14
5.015	1.01E-14
5.02	1.01E-14
5.025	1.01E-14
5.03	1.01E-14
5.035	1.01E-14
5.04	1.01E-14
5.045	1.02E-14
5.05	1.01E-14
5.055	1.02E-14
5.06	1.02E-14
5.065	1.02E-14
5.07	1.02E-14
5.075	1.02E-14
5.08	1.02E-14
5.085	1.02E-14
5.09	1.02E-14
5.095	1.02E-14
5.1	1.03E-14
5.105	1.03E-14
5.11	1.03E-14
5.115	1.03E-14
5.12	1.03E-14
5.125	1.02E-14
5.13	1.02E-14
5.135	1.02E-14
5.14	1.02E-14
5.145	1.02E-14
5.15	1.02E-14
5.155	1.02E-14
5.16	1.01E-14
5.165	1.01E-14
5.17	1.01E-14
5.175	1.01E-14
5.18	1.00E-14
5.185	1.00E-14
5.19	1.00E-14
5.195	1.00E-14
5.2	1.00E-14
5.205	1.00E-14
5.21	1.00E-14
5.215	9.99E-15
5.22	9.99E-15
5.225	9.99E-15
5.23	9.99E-15
5.235	9.99E-15
5.24	9.99E-15
5.245	9.99E-15
5.25	9.99E-15
5.255	9.99E-15
5.26	9.99E-15
5.265	9.99E-15
5.27	9.99E-15
5.275	1.00E-14
5.28	1.00E-14
5.285	1.00E-14
5.29	1.00E-14
5.295	1.00E-14
5.3	1.00E-14
5.305	1.01E-14
5.31	1.01E-14
5.315	1.01E-14
5.32	1.01E-14
5.325	1.01E-14
5.33	1.01E-14
5.335	1.02E-14
5.34	1.02E-14
5.345	1.02E-14
5.35	1.02E-14
5.355	1.02E-14
5.36	1.02E-14
5.365	1.02E-14
5.37	1.02E-14
5.375	1.02E-14
5.38	1.03E-14
5.385	1.03E-14
5.39	1.03E-14
5.395	1.03E-14
5.4	1.02E-14
5.405	1.02E-14
5.41	1.02E-14
5.415	1.02E-14
5.42	1.02E-14
5.425	1.02E-14
5.43	1.02E-14
5.435	1.01E-14
5.44	1.01E-14
5.445	1.01E-14
5.45	1.01E-14
5.455	1.01E-14
5.46	1.01E-14
5.465	1.01E-14
5.47	1.00E-14
5.475	1.00E-14
5.48	1.00E-14
5.485	1.00E-14
5.49	1.00E-14
5.495	1.00E-14
5.5	9.99E-15
5.505	1.00E-14
5.51	9.99E-15
5.515	9.99E-15
5.52	9.99E-15
5.525	9.99E-15
5.53	9.99E-15
5.535	9.99E-15
5.54	9.99E-15
5.545	9.99E-15
5.55	9.99E-15
5.555	9.99E-15
5.56	9.99E-15
5.565	9.99E-15
5.57	9.88E-15
5.575	9.88E-15
5.58	9.88E-15
5.585	9.88E-15
5.59	9.88E-15
5.595	9.88E-15
5.6	9.88E-15
5.605	9.88E-15
5.61	9.88E-15
5.615	9.88E-15
5.62	9.88E-15
5.625	9.88E-15
5.63	9.88E-15
5.635	9.88E-15
5.64	9.99E-15
5.645	9.99E-15
5.65	9.99E-15
5.655	9.99E-15
5.66	1.00E-14
5.665	1.00E-14
5.67	1.00E-14
5.675	1.01E-14
5.68	1.01E-14
5.685	1.01E-14
5.69	1.01E-14
5.695	1.02E-14
5.7	1.02E-14
5.705	1.03E-14
5.71	1.03E-14
5.715	1.03E-14
5.72	1.03E-14
5.725	1.03E-14
5.73	1.04E-14
5.735	1.04E-14
5.74	1.04E-14
5.745	1.05E-14
5.75	1.05E-14
5.755	1.05E-14
5.76	1.05E-14
5.765	1.05E-14
5.77	1.05E-14
5.775	1.05E-14
5.78	1.05E-14
5.785	1.05E-14
5.79	1.05E-14
5.795	1.05E-14
5.8	1.05E-14
5.805	1.06E-14
5.81	1.06E-14
5.815	1.06E-14
5.82	1.06E-14
5.825	1.06E-14
5.83	1.06E-14
5.835	1.07E-14
5.84	1.07E-14
5.845	1.07E-14
5.85	1.08E-14
5.855	1.08E-14
5.86	1.08E-14
5.865	1.08E-14
5.87	1.08E-14
5.875	1.08E-14
5.88	1.08E-14
5.885	1.09E-14
5.89	1.09E-14
5.895	1.10E-14
5.9	1.10E-14
5.905	1.10E-14
5.91	1.10E-14
5.915	1.10E-14
5.92	1.10E-14
5.925	1.10E-14
5.93	1.10E-14
5.935	1.10E-14
5.94	1.11E-14
5.945	1.11E-14
5.95	1.11E-14
5.955	1.11E-14
5.96	1.11E-14
5.965	1.11E-14
5.97	1.11E-14
5.975	1.11E-14
5.98	1.11E-14
5.985	1.11E-14
5.99	1.11E-14
5.995	1.11E-14
6	1.11E-14
6.005	1.11E-14
6.01	1.11E-14
6.015	1.11E-14
6.02	1.11E-14
6.025	1.11E-14
6.03	1.11E-14
6.035	1.12E-14
6.04	1.12E-14
6.045	1.12E-14
6.05	1.12E-14
6.055	1.12E-14
6.06	1.12E-14
6.065	1.13E-14
6.07	1.13E-14
6.075	1.13E-14
6.08	1.13E-14
6.085	1.13E-14
6.09	1.13E-14
6.095	1.14E-14
6.1	1.14E-14
6.105	1.14E-14
6.11	1.15E-14
6.115	1.15E-14
6.12	1.15E-14
6.125	1.15E-14
6.13	1.15E-14
6.135	1.16E-14
6.14	1.16E-14
6.145	1.16E-14
6.15	1.17E-14
6.155	1.17E-14
6.16	1.17E-14
6.165	1.17E-14
6.17	1.17E-14
6.175	1.17E-14
6.18	1.17E-14
6.185	1.18E-14
6.19	1.18E-14
6.195	1.18E-14
6.2	1.18E-14
6.205	1.18E-14
6.21	1.17E-14
6.215	1.18E-14
6.22	1.17E-14
6.225	1.18E-14
6.23	1.18E-14
6.235	1.17E-14
6.24	1.17E-14
6.245	1.17E-14
6.25	1.17E-14
6.255	1.17E-14
6.26	1.17E-14
6.265	1.17E-14
6.27	1.17E-14
6.275	1.17E-14
6.28	1.17E-14
6.285	1.17E-14
6.29	1.17E-14
6.295	1.17E-14
6.3	1.17E-14
6.305	1.17E-14
6.31	1.17E-14
6.315	1.16E-14
6.32	1.16E-14
6.325	1.16E-14
6.33	1.16E-14
6.335	1.16E-14
6.34	1.16E-14
6.345	1.15E-14
6.35	1.15E-14
6.355	1.15E-14
6.36	1.15E-14
6.365	1.15E-14
6.37	1.15E-14
6.375	1.14E-14
6.38	1.14E-14
6.385	1.14E-14
6.39	1.14E-14
6.395	1.14E-14
6.4	1.14E-14
6.405	1.13E-14
6.41	1.13E-14
6.415	1.13E-14
6.42	1.13E-14
6.425	1.13E-14
6.43	1.13E-14
6.435	1.13E-14
6.44	1.13E-14
6.445	1.13E-14
6.45	1.13E-14
6.455	1.13E-14
6.46	1.13E-14
6.465	1.13E-14
6.47	1.13E-14
6.475	1.13E-14
6.48	1.13E-14
6.485	1.13E-14
6.49	1.13E-14
6.495	1.13E-14
6.5	1.13E-14
6.505	1.13E-14
6.51	1.13E-14
6.515	1.13E-14
6.52	1.13E-14
6.525	1.13E-14
6.53	1.13E-14
6.535	1.13E-14
6.54	1.14E-14
6.545	1.14E-14
6.55	1.14E-14
6.555	1.14E-14
6.56	1.14E-14
6.565	1.14E-14
6.57	1.14E-14
6.575	1.14E-14
6.58	1.15E-14
6.585	1.14E-14
6.59	1.14E-14
6.595	1.15E-14
6.6	1.15E-14
6.605	1.15E-14
6.61	1.15E-14
6.615	1.15E-14
6.62	1.15E-14
6.625	1.15E-14
6.63	1.15E-14
6.635	1.15E-14
6.64	1.15E-14
6.645	1.15E-14
6.65	1.15E-14
6.655	1.15E-14
6.66	1.15E-14
6.665	1.15E-14
6.67	1.15E-14
6.675	1.16E-14
6.68	1.16E-14
6.685	1.16E-14
6.69	1.16E-14
6.695	1.16E-14
6.7	1.16E-14
6.705	1.16E-14
6.71	1.16E-14
6.715	1.16E-14
6.72	1.16E-14
6.725	1.16E-14
6.73	1.16E-14
6.735	1.16E-14
6.74	1.16E-14
6.745	1.16E-14
6.75	1.16E-14
6.755	1.16E-14
6.76	1.16E-14
6.765	1.16E-14
6.77	1.16E-14
6.775	1.16E-14
6.78	1.16E-14
6.785	1.16E-14
6.79	1.16E-14
6.795	1.15E-14
6.8	1.15E-14
6.805	1.15E-14
6.81	1.15E-14
6.815	1.15E-14
6.82	1.15E-14
6.825	1.15E-14
6.83	1.15E-14
6.835	1.15E-14
6.84	1.15E-14
6.845	1.15E-14
6.85	1.15E-14
6.855	1.15E-14
6.86	1.15E-14
6.865	1.15E-14
6.87	1.14E-14
6.875	1.15E-14
6.88	1.15E-14
6.885	1.14E-14
6.89	1.15E-14
6.895	1.15E-14
6.9	1.14E-14
6.905	1.14E-14
6.91	1.14E-14
6.915	1.14E-14
6.92	1.14E-14
6.925	1.14E-14
6.93	1.14E-14
6.935	1.14E-14
6.94	1.13E-14
6.945	1.13E-14
6.95	1.13E-14
6.955	1.13E-14
6.96	1.13E-14
6.965	1.13E-14
6.97	1.13E-14
6.975	1.13E-14
6.98	1.13E-14
6.985	1.13E-14
6.99	1.13E-14
6.995	1.12E-14
7	1.12E-14
7.005	1.12E-14
7.01	1.12E-14
7.015	1.12E-14
7.02	1.12E-14
7.025	1.12E-14
7.03	1.12E-14
7.035	1.12E-14
7.04	1.12E-14
7.045	1.12E-14
7.05	1.12E-14
7.055	1.12E-14
7.06	1.12E-14
7.065	1.12E-14
7.07	1.13E-14
7.075	1.13E-14
7.08	1.13E-14
7.085	1.13E-14
7.09	1.13E-14
7.095	1.13E-14
7.1	1.13E-14
7.105	1.13E-14
7.11	1.13E-14
7.115	1.13E-14
7.12	1.14E-14
7.125	1.14E-14
7.13	1.14E-14
7.135	1.14E-14
7.14	1.15E-14
7.145	1.15E-14
7.15	1.15E-14
7.155	1.15E-14
7.16	1.15E-14
7.165	1.15E-14
7.17	1.15E-14
7.175	1.15E-14
7.18	1.15E-14
7.185	1.15E-14
7.19	1.15E-14
7.195	1.15E-14
7.2	1.15E-14
7.205	1.16E-14
7.21	1.16E-14
7.215	1.16E-14
7.22	1.16E-14
7.225	1.16E-14
7.23	1.17E-14
7.235	1.17E-14
7.24	1.17E-14
7.245	1.17E-14
7.25	1.17E-14
7.255	1.18E-14
7.26	1.18E-14
7.265	1.18E-14
7.27	1.18E-14
7.275	1.18E-14
7.28	1.18E-14
7.285	1.18E-14
7.29	1.18E-14
7.295	1.18E-14
7.3	1.18E-14
7.305	1.18E-14
7.31	1.18E-14
7.315	1.18E-14
7.32	1.18E-14
7.325	1.18E-14
7.33	1.18E-14
7.335	1.18E-14
7.34	1.18E-14
7.345	1.18E-14
7.35	1.18E-14
7.355	1.18E-14
7.36	1.18E-14
7.365	1.18E-14
7.37	1.18E-14
7.375	1.18E-14
7.38	1.18E-14
7.385	1.18E-14
7.39	1.18E-14
7.395	1.18E-14
7.4	1.18E-14
7.405	1.18E-14
7.41	1.18E-14
7.415	1.18E-14
7.42	1.18E-14
7.425	1.18E-14
7.43	1.18E-14
7.435	1.18E-14
7.44	1.18E-14
7.445	1.18E-14
7.45	1.18E-14
7.455	1.18E-14
7.46	1.18E-14
7.465	1.18E-14
7.47	1.18E-14
7.475	1.18E-14
7.48	1.19E-14
7.485	1.19E-14
7.49	1.19E-14
7.495	1.19E-14
7.5	1.19E-14
7.505	1.19E-14
7.51	1.19E-14
7.515	1.19E-14
7.52	1.19E-14
7.525	1.19E-14
7.53	1.19E-14
7.535	1.19E-14
7.54	1.19E-14
7.545	1.19E-14
7.55	1.19E-14
7.555	1.19E-14
7.56	1.20E-14
7.565	1.20E-14
7.57	1.20E-14
7.575	1.20E-14
7.58	1.20E-14
7.585	1.20E-14
7.59	1.20E-14
7.595	1.21E-14
7.6	1.21E-14
7.605	1.21E-14
7.61	1.22E-14
7.615	1.22E-14
7.62	1.22E-14
7.625	1.23E-14
7.63	1.23E-14
7.635	1.23E-14
7.64	1.23E-14
7.645	1.24E-14
7.65	1.24E-14
7.655	1.25E-14
7.66	1.25E-14
7.665	1.26E-14
7.67	1.27E-14
7.675	1.27E-14
7.68	1.28E-14
7.685	1.28E-14
7.69	1.28E-14
7.695	1.28E-14
7.7	1.29E-14
7.705	1.29E-14
7.71	1.29E-14
7.715	1.30E-14
7.72	1.30E-14
7.725	1.30E-14
7.73	1.30E-14
7.735	1.31E-14
7.74	1.31E-14
7.745	1.31E-14
7.75	1.31E-14
7.755	1.31E-14
7.76	1.31E-14
7.765	1.31E-14
7.77	1.31E-14
7.775	1.31E-14
7.78	1.31E-14
7.785	1.31E-14
7.79	1.31E-14
7.795	1.30E-14
7.8	1.30E-14
7.805	1.30E-14
7.81	1.30E-14
7.815	1.29E-14
7.82	1.29E-14
7.825	1.29E-14
7.83	1.29E-14
7.835	1.29E-14
7.84	1.28E-14
7.845	1.28E-14
7.85	1.28E-14
7.855	1.28E-14
7.86	1.28E-14
7.865	1.28E-14
7.87	1.28E-14
7.875	1.28E-14
7.88	1.28E-14
7.885	1.28E-14
7.89	1.28E-14
7.895	1.28E-14
7.9	1.28E-14
7.905	1.29E-14
7.91	1.29E-14
7.915	1.29E-14
7.92	1.29E-14
7.925	1.29E-14
7.93	1.29E-14
7.935	1.30E-14
7.94	1.31E-14
7.945	1.31E-14
7.95	1.32E-14
7.955	1.33E-14
7.96	1.33E-14
7.965	1.33E-14
7.97	1.33E-14
7.975	1.34E-14
7.98	1.34E-14
7.985	1.34E-14
7.99	1.34E-14
7.995	1.35E-14
8	1.35E-14
8.005	1.36E-14
8.01	1.37E-14
8.015	1.37E-14
8.02	1.37E-14
8.025	1.37E-14
8.03	1.38E-14
8.035	1.38E-14
8.04	1.39E-14
8.045	1.39E-14
8.05	1.39E-14
8.055	1.39E-14
8.06	1.39E-14
8.065	1.40E-14
8.07	1.40E-14
8.075	1.40E-14
8.08	1.40E-14
8.085	1.41E-14
8.09	1.41E-14
8.095	1.41E-14
8.1	1.41E-14
8.105	1.41E-14
8.11	1.41E-14
8.115	1.41E-14
8.12	1.41E-14
8.125	1.42E-14
8.13	1.41E-14
8.135	1.42E-14
8.14	1.41E-14
8.145	1.41E-14
8.15	1.42E-14
8.155	1.41E-14
8.16	1.41E-14
8.165	1.41E-14
8.17	1.41E-14
8.175	1.41E-14
8.18	1.41E-14
8.185	1.41E-14
8.19	1.41E-14
8.195	1.42E-14
8.2	1.41E-14
8.205	1.42E-14
8.21	1.42E-14
8.215	1.42E-14
8.22	1.42E-14
8.225	1.42E-14
8.23	1.42E-14
8.235	1.43E-14
8.24	1.43E-14
8.245	1.44E-14
8.25	1.44E-14
8.255	1.44E-14
8.26	1.44E-14
8.265	1.45E-14
8.27	1.46E-14
8.275	1.46E-14
8.28	1.47E-14
8.285	1.47E-14
8.29	1.47E-14
8.295	1.48E-14
8.3	1.48E-14
8.305	1.49E-14
8.31	1.49E-14
8.315	1.50E-14
8.32	1.50E-14
8.325	1.51E-14
8.33	1.52E-14
8.335	1.52E-14
8.34	1.52E-14
8.345	1.53E-14
8.35	1.53E-14
8.355	1.54E-14
8.36	1.54E-14
8.365	1.54E-14
8.37	1.55E-14
8.375	1.55E-14
8.38	1.56E-14
8.385	1.56E-14
8.39	1.56E-14
8.395	1.57E-14
8.4	1.57E-14
8.405	1.57E-14
8.41	1.57E-14
8.415	1.57E-14
8.42	1.58E-14
8.425	1.58E-14
8.43	1.58E-14
8.435	1.58E-14
8.44	1.59E-14
8.445	1.59E-14
8.45	1.59E-14
8.455	1.59E-14
8.46	1.59E-14
8.465	1.59E-14
8.47	1.59E-14
8.475	1.60E-14
8.48	1.60E-14
8.485	1.60E-14
8.49	1.61E-14
8.495	1.61E-14
8.5	1.61E-14
8.505	1.62E-14
8.51	1.62E-14
8.515	1.62E-14
8.52	1.62E-14
8.525	1.62E-14
8.53	1.63E-14
8.535	1.63E-14
8.54	1.64E-14
8.545	1.64E-14
8.55	1.64E-14
8.555	1.65E-14
8.56	1.65E-14
8.565	1.66E-14
8.57	1.66E-14
8.575	1.67E-14
8.58	1.67E-14
8.585	1.67E-14
8.59	1.67E-14
8.595	1.68E-14
8.6	1.68E-14
8.605	1.69E-14
8.61	1.69E-14
8.615	1.69E-14
8.62	1.69E-14
8.625	1.69E-14
8.63	1.70E-14
8.635	1.71E-14
8.64	1.71E-14
8.645	1.72E-14
8.65	1.72E-14
8.655	1.72E-14
8.66	1.72E-14
8.665	1.73E-14
8.67	1.73E-14
8.675	1.74E-14
8.68	1.74E-14
8.685	1.74E-14
8.69	1.75E-14
8.695	1.75E-14
8.7	1.76E-14
8.705	1.76E-14
8.71	1.77E-14
8.715	1.77E-14
8.72	1.77E-14
8.725	1.77E-14
8.73	1.77E-14
8.735	1.78E-14
8.74	1.78E-14
8.745	1.78E-14
8.75	1.78E-14
8.755	1.79E-14
8.76	1.78E-14
8.765	1.79E-14
8.77	1.79E-14
8.775	1.79E-14
8.78	1.79E-14
8.785	1.79E-14
8.79	1.79E-14
8.795	1.79E-14
8.8	1.79E-14
8.805	1.79E-14
8.81	1.79E-14
8.815	1.79E-14
8.82	1.79E-14
8.825	1.79E-14
8.83	1.79E-14
8.835	1.79E-14
8.84	1.79E-14
8.845	1.79E-14
8.85	1.79E-14
8.855	1.78E-14
8.86	1.78E-14
8.865	1.78E-14
8.87	1.78E-14
8.875	1.79E-14
8.88	1.79E-14
8.885	1.79E-14
8.89	1.79E-14
8.895	1.79E-14
8.9	1.79E-14
8.905	1.79E-14
8.91	1.79E-14
8.915	1.80E-14
8.92	1.80E-14
8.925	1.80E-14
8.93	1.81E-14
8.935	1.81E-14
8.94	1.82E-14
8.945	1.82E-14
8.95	1.82E-14
8.955	1.83E-14
8.96	1.84E-14
8.965	1.84E-14
8.97	1.84E-14
8.975	1.85E-14
8.98	1.86E-14
8.985	1.87E-14
8.99	1.87E-14
8.995	1.88E-14
9	1.89E-14
9.005	1.89E-14
9.01	1.89E-14
9.015	1.90E-14
9.02	1.92E-14
9.025	1.92E-14
9.03	1.92E-14
9.035	1.93E-14
9.04	1.94E-14
9.045	1.94E-14
9.05	1.95E-14
9.055	1.95E-14
9.06	1.95E-14
9.065	1.97E-14
9.07	1.97E-14
9.075	1.97E-14
9.08	1.98E-14
9.085	1.98E-14
9.09	1.99E-14
9.095	1.99E-14
9.1	1.99E-14
9.105	2.00E-14
9.11	2.00E-14
9.115	2.02E-14
9.12	2.02E-14
9.125	2.02E-14
9.13	2.02E-14
9.135	2.03E-14
9.14	2.04E-14
9.145	2.04E-14
9.15	2.04E-14
9.155	2.05E-14
9.16	2.05E-14
9.165	2.05E-14
9.17	2.07E-14
9.175	2.07E-14
9.18	2.07E-14
9.185	2.07E-14
9.19	2.08E-14
9.195	2.08E-14
9.2	2.08E-14
9.205	2.08E-14
9.21	2.09E-14
9.215	2.09E-14
9.22	2.10E-14
9.225	2.10E-14
9.23	2.10E-14
9.235	2.10E-14
9.24	2.10E-14
9.245	2.10E-14
9.25	2.11E-14
9.255	2.11E-14
9.26	2.10E-14
9.265	2.11E-14
9.27	2.11E-14
9.275	2.11E-14
9.28	2.12E-14
9.285	2.12E-14
9.29	2.12E-14
9.295	2.12E-14
9.3	2.13E-14
9.305	2.13E-14
9.31	2.13E-14
9.315	2.13E-14
9.32	2.13E-14
9.325	2.13E-14
9.33	2.14E-14
9.335	2.14E-14
9.34	2.14E-14
9.345	2.15E-14
9.35	2.15E-14
9.355	2.15E-14
9.36	2.15E-14
9.365	2.15E-14
9.37	2.16E-14
9.375	2.16E-14
9.38	2.17E-14
9.385	2.17E-14
9.39	2.17E-14
9.395	2.18E-14
9.4	2.18E-14
9.405	2.18E-14
9.41	2.18E-14
9.415	2.18E-14
9.42	2.18E-14
9.425	2.18E-14
9.43	2.18E-14
9.435	2.18E-14
9.44	2.19E-14
9.445	2.19E-14
9.45	2.19E-14
9.455	2.19E-14
9.46	2.19E-14
9.465	2.19E-14
9.47	2.20E-14
9.475	2.20E-14
9.48	2.20E-14
9.485	2.20E-14
9.49	2.20E-14
9.495	2.20E-14
9.5	2.20E-14
9.505	2.20E-14
9.51	2.20E-14
9.515	2.20E-14
9.52	2.20E-14
9.525	2.20E-14
9.53	2.20E-14
9.535	2.20E-14
9.54	2.20E-14
9.545	2.20E-14
9.55	2.20E-14
9.555	2.20E-14
9.56	2.20E-14
9.565	2.20E-14
9.57	2.20E-14
9.575	2.20E-14
9.58	2.20E-14
9.585	2.20E-14
9.59	2.20E-14
9.595	2.20E-14
9.6	2.20E-14
9.605	2.20E-14
9.61	2.20E-14
9.615	2.20E-14
9.62	2.19E-14
9.625	2.19E-14
9.63	2.19E-14
9.635	2.19E-14
9.64	2.19E-14
9.645	2.18E-14
9.65	2.18E-14
9.655	2.18E-14
9.66	2.18E-14
9.665	2.18E-14
9.67	2.18E-14
9.675	2.18E-14
9.68	2.18E-14
9.685	2.17E-14
9.69	2.17E-14
9.695	2.17E-14
9.7	2.16E-14
9.705	2.16E-14
9.71	2.16E-14
9.715	2.16E-14
9.72	2.15E-14
9.725	2.15E-14
9.73	2.15E-14
9.735	2.15E-14
9.74	2.15E-14
9.745	2.15E-14
9.75	2.15E-14
9.755	2.15E-14
9.76	2.14E-14
9.765	2.14E-14
9.77	2.14E-14
9.775	2.14E-14
9.78	2.14E-14
9.785	2.13E-14
9.79	2.13E-14
9.795	2.13E-14
9.8	2.13E-14
9.805	2.13E-14
9.81	2.13E-14
9.815	2.13E-14
9.82	2.13E-14
9.825	2.13E-14
9.83	2.13E-14
9.835	2.13E-14
9.84	2.13E-14
9.845	2.13E-14
9.85	2.13E-14
9.855	2.13E-14
9.86	2.13E-14
9.865	2.13E-14
9.87	2.12E-14
9.875	2.12E-14
9.88	2.12E-14
9.885	2.12E-14
9.89	2.12E-14
9.895	2.11E-14
9.9	2.11E-14
9.905	2.11E-14
9.91	2.11E-14
9.915	2.11E-14
9.92	2.11E-14
9.925	2.10E-14
9.93	2.10E-14
9.935	2.10E-14
9.94	2.10E-14
9.945	2.10E-14
9.95	2.10E-14
9.955	2.10E-14
9.96	2.10E-14
9.965	2.10E-14
9.97	2.09E-14
9.975	2.09E-14
9.98	2.09E-14
9.985	2.09E-14
9.99	2.09E-14
9.995	2.09E-14
10	2.08E-14
};
\end{semilogyaxis}
\end{tikzpicture}

\caption{Error of the energy for the case $r=20$ for the KdV equation: $|H(V\bz _n) - H(V\bz_0)|$ are plotted.
}
\label{fig:kdv1:EnergyError}
\end{figure}

%% file: figkdv1sol.tex
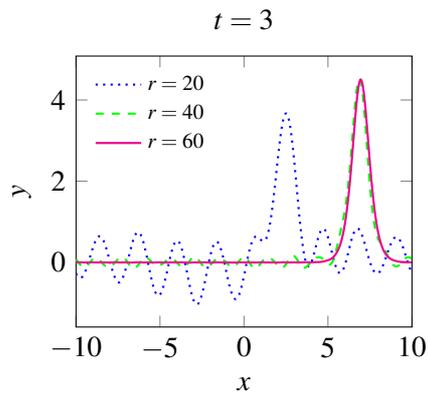
\begin{figure}[htbp]
\centering

\begin{tikzpicture}
\tikzstyle{every node}=[]
\begin{axis}[width=6cm,
xmax=10,xmin=-10,
xlabel={$x$},ylabel={$y$},
ylabel near ticks,
legend entries={$r=20$,$r=40$,$r=60$},
legend style={font = \scriptsize,legend pos=north west,legend cell align=left,draw=none,fill=none},
title = {$t=3$},
	]
\addplot[thick,dotted,color=blue
] table {
-10	-0.272
-9.96	-0.303
-9.92	-0.33
-9.88	-0.351
-9.84	-0.367
-9.8	-0.377
-9.76	-0.382
-9.72	-0.38
-9.68	-0.373
-9.64	-0.36
-9.6	-0.341
-9.56	-0.317
-9.52	-0.287
-9.48	-0.253
-9.44	-0.214
-9.4	-0.171
-9.36	-0.124
-9.32	-0.075
-9.28	-0.023
-9.24	0.031
-9.2	0.087
-9.16	0.142
-9.12	0.198
-9.08	0.253
-9.04	0.307
-9	0.359
-8.96	0.408
-8.92	0.454
-8.88	0.495
-8.84	0.533
-8.8	0.566
-8.76	0.593
-8.72	0.615
-8.68	0.631
-8.64	0.641
-8.6	0.644
-8.56	0.641
-8.52	0.632
-8.48	0.616
-8.44	0.594
-8.4	0.566
-8.36	0.533
-8.32	0.494
-8.28	0.45
-8.24	0.401
-8.2	0.349
-8.16	0.294
-8.12	0.235
-8.08	0.175
-8.04	0.114
-8	0.052
-7.96	-0.01
-7.92	-0.071
-7.88	-0.13
-7.84	-0.187
-7.8	-0.241
-7.76	-0.291
-7.72	-0.337
-7.68	-0.377
-7.64	-0.413
-7.6	-0.442
-7.56	-0.465
-7.52	-0.482
-7.48	-0.491
-7.44	-0.494
-7.4	-0.489
-7.36	-0.478
-7.32	-0.459
-7.28	-0.433
-7.24	-0.401
-7.2	-0.363
-7.16	-0.319
-7.12	-0.269
-7.08	-0.215
-7.04	-0.157
-7	-0.095
-6.96	-0.03
-6.92	0.036
-6.88	0.104
-6.84	0.172
-6.8	0.24
-6.76	0.306
-6.72	0.371
-6.68	0.432
-6.64	0.49
-6.6	0.543
-6.56	0.592
-6.52	0.634
-6.48	0.67
-6.44	0.7
-6.4	0.722
-6.36	0.737
-6.32	0.744
-6.28	0.743
-6.24	0.734
-6.2	0.716
-6.16	0.691
-6.12	0.658
-6.08	0.618
-6.04	0.57
-6	0.516
-5.96	0.456
-5.92	0.391
-5.88	0.321
-5.84	0.247
-5.8	0.169
-5.76	0.09
-5.72	0.008
-5.68	-0.073
-5.64	-0.155
-5.6	-0.235
-5.56	-0.313
-5.52	-0.388
-5.48	-0.459
-5.44	-0.526
-5.4	-0.587
-5.36	-0.642
-5.32	-0.691
-5.28	-0.732
-5.24	-0.766
-5.2	-0.792
-5.16	-0.809
-5.12	-0.818
-5.08	-0.819
-5.04	-0.811
-5	-0.794
-4.96	-0.769
-4.92	-0.737
-4.88	-0.697
-4.84	-0.65
-4.8	-0.597
-4.76	-0.538
-4.72	-0.474
-4.68	-0.406
-4.64	-0.335
-4.6	-0.261
-4.56	-0.186
-4.52	-0.11
-4.48	-0.034
-4.44	0.04
-4.4	0.112
-4.36	0.18
-4.32	0.245
-4.28	0.305
-4.24	0.359
-4.2	0.407
-4.16	0.448
-4.12	0.482
-4.08	0.508
-4.04	0.525
-4	0.534
-3.96	0.534
-3.92	0.525
-3.88	0.508
-3.84	0.482
-3.8	0.448
-3.76	0.406
-3.72	0.356
-3.68	0.3
-3.64	0.237
-3.6	0.169
-3.56	0.095
-3.52	0.018
-3.48	-0.062
-3.44	-0.145
-3.4	-0.229
-3.36	-0.313
-3.32	-0.396
-3.28	-0.478
-3.24	-0.557
-3.2	-0.633
-3.16	-0.704
-3.12	-0.77
-3.08	-0.83
-3.04	-0.884
-3	-0.93
-2.96	-0.968
-2.92	-0.998
-2.88	-1.019
-2.84	-1.032
-2.8	-1.035
-2.76	-1.03
-2.72	-1.015
-2.68	-0.992
-2.64	-0.96
-2.6	-0.919
-2.56	-0.872
-2.52	-0.817
-2.48	-0.755
-2.44	-0.688
-2.4	-0.616
-2.36	-0.54
-2.32	-0.461
-2.28	-0.379
-2.24	-0.296
-2.2	-0.213
-2.16	-0.13
-2.12	-0.049
-2.08	0.029
-2.04	0.104
-2	0.175
-1.96	0.241
-1.92	0.301
-1.88	0.354
-1.84	0.401
-1.8	0.439
-1.76	0.47
-1.72	0.492
-1.68	0.505
-1.64	0.51
-1.6	0.506
-1.56	0.493
-1.52	0.472
-1.48	0.443
-1.44	0.406
-1.4	0.361
-1.36	0.31
-1.32	0.253
-1.28	0.19
-1.24	0.123
-1.2	0.052
-1.16	-0.023
-1.12	-0.099
-1.08	-0.176
-1.04	-0.254
-1	-0.33
-0.96	-0.406
-0.92	-0.478
-0.88	-0.547
-0.84	-0.612
-0.8	-0.672
-0.76	-0.727
-0.72	-0.775
-0.68	-0.816
-0.64	-0.849
-0.6	-0.875
-0.56	-0.893
-0.52	-0.902
-0.48	-0.903
-0.44	-0.896
-0.4	-0.88
-0.36	-0.856
-0.32	-0.823
-0.28	-0.783
-0.24	-0.736
-0.2	-0.682
-0.16	-0.621
-0.12	-0.555
-0.08	-0.485
-0.04	-0.41
0	-0.332
0.04	-0.252
0.08	-0.17
0.12	-0.088
0.16	-0.007
0.2	0.073
0.24	0.15
0.28	0.224
0.32	0.293
0.36	0.357
0.4	0.416
0.44	0.468
0.48	0.513
0.52	0.552
0.56	0.583
0.6	0.607
0.64	0.623
0.68	0.634
0.72	0.637
0.76	0.636
0.8	0.629
0.84	0.618
0.88	0.605
0.92	0.589
0.96	0.573
1	0.556
1.04	0.542
1.08	0.53
1.12	0.523
1.16	0.521
1.2	0.525
1.24	0.536
1.28	0.556
1.32	0.586
1.36	0.625
1.4	0.674
1.44	0.734
1.48	0.806
1.52	0.888
1.56	0.981
1.6	1.085
1.64	1.198
1.68	1.321
1.72	1.452
1.76	1.591
1.8	1.735
1.84	1.885
1.88	2.038
1.92	2.193
1.96	2.348
2	2.502
2.04	2.652
2.08	2.798
2.12	2.938
2.16	3.069
2.2	3.191
2.24	3.302
2.28	3.4
2.32	3.484
2.36	3.554
2.4	3.608
2.44	3.645
2.48	3.665
2.52	3.667
2.56	3.652
2.6	3.618
2.64	3.568
2.68	3.5
2.72	3.415
2.76	3.315
2.8	3.201
2.84	3.072
2.88	2.932
2.92	2.78
2.96	2.619
3	2.451
3.04	2.277
3.08	2.098
3.12	1.917
3.16	1.736
3.2	1.555
3.24	1.378
3.28	1.205
3.32	1.038
3.36	0.879
3.4	0.729
3.44	0.59
3.48	0.461
3.52	0.345
3.56	0.241
3.6	0.151
3.64	0.075
3.68	0.013
3.72	-0.035
3.76	-0.069
3.8	-0.09
3.84	-0.097
3.88	-0.093
3.92	-0.077
3.96	-0.05
4	-0.015
4.04	0.03
4.08	0.081
4.12	0.138
4.16	0.199
4.2	0.263
4.24	0.329
4.28	0.395
4.32	0.461
4.36	0.524
4.4	0.583
4.44	0.638
4.48	0.688
4.52	0.731
4.56	0.767
4.6	0.796
4.64	0.816
4.68	0.827
4.72	0.83
4.76	0.825
4.8	0.81
4.84	0.788
4.88	0.757
4.92	0.719
4.96	0.674
5	0.622
5.04	0.566
5.08	0.505
5.12	0.441
5.16	0.374
5.2	0.307
5.24	0.238
5.28	0.171
5.32	0.105
5.36	0.042
5.4	-0.017
5.44	-0.071
5.48	-0.12
5.52	-0.162
5.56	-0.198
5.6	-0.226
5.64	-0.247
5.68	-0.259
5.72	-0.264
5.76	-0.259
5.8	-0.247
5.84	-0.227
5.88	-0.198
5.92	-0.163
5.96	-0.121
6	-0.072
6.04	-0.019
6.08	0.039
6.12	0.101
6.16	0.165
6.2	0.232
6.24	0.299
6.28	0.365
6.32	0.431
6.36	0.494
6.4	0.554
6.44	0.61
6.48	0.661
6.52	0.707
6.56	0.746
6.6	0.779
6.64	0.804
6.68	0.822
6.72	0.832
6.76	0.835
6.8	0.83
6.84	0.817
6.88	0.797
6.92	0.771
6.96	0.738
7	0.699
7.04	0.655
7.08	0.607
7.12	0.555
7.16	0.5
7.2	0.442
7.24	0.383
7.28	0.324
7.32	0.264
7.36	0.205
7.4	0.147
7.44	0.091
7.48	0.037
7.52	-0.013
7.56	-0.061
7.6	-0.104
7.64	-0.143
7.68	-0.177
7.72	-0.207
7.76	-0.231
7.8	-0.25
7.84	-0.263
7.88	-0.271
7.92	-0.273
7.96	-0.27
8	-0.261
8.04	-0.247
8.08	-0.229
8.12	-0.205
8.16	-0.176
8.2	-0.144
8.24	-0.108
8.28	-0.068
8.32	-0.026
8.36	0.019
8.4	0.065
8.44	0.113
8.48	0.162
8.52	0.21
8.56	0.259
8.6	0.306
8.64	0.352
8.68	0.395
8.72	0.437
8.76	0.475
8.8	0.509
8.84	0.54
8.88	0.566
8.92	0.588
8.96	0.605
9	0.616
9.04	0.623
9.08	0.623
9.12	0.619
9.16	0.609
9.2	0.594
9.24	0.573
9.28	0.548
9.32	0.518
9.36	0.483
9.4	0.444
9.44	0.402
9.48	0.357
9.52	0.309
9.56	0.258
9.6	0.206
9.64	0.153
9.68	0.1
9.72	0.047
9.76	-0.006
9.8	-0.057
9.84	-0.106
9.88	-0.152
9.92	-0.196
9.96	-0.236
};
\addplot[thick,dashed,color=green
] table {
-10	0.077
-9.96	0.058
-9.92	0.037
-9.88	0.016
-9.84	-0.007
-9.8	-0.029
-9.76	-0.049
-9.72	-0.069
-9.68	-0.085
-9.64	-0.099
-9.6	-0.109
-9.56	-0.115
-9.52	-0.118
-9.48	-0.116
-9.44	-0.11
-9.4	-0.101
-9.36	-0.088
-9.32	-0.072
-9.28	-0.054
-9.24	-0.035
-9.2	-0.014
-9.16	0.007
-9.12	0.027
-9.08	0.046
-9.04	0.064
-9	0.079
-8.96	0.091
-8.92	0.1
-8.88	0.105
-8.84	0.108
-8.8	0.106
-8.76	0.101
-8.72	0.094
-8.68	0.083
-8.64	0.07
-8.6	0.054
-8.56	0.038
-8.52	0.02
-8.48	0.001
-8.44	-0.017
-8.4	-0.035
-8.36	-0.052
-8.32	-0.067
-8.28	-0.081
-8.24	-0.093
-8.2	-0.103
-8.16	-0.11
-8.12	-0.114
-8.08	-0.116
-8.04	-0.115
-8	-0.111
-7.96	-0.105
-7.92	-0.096
-7.88	-0.084
-7.84	-0.071
-7.8	-0.055
-7.76	-0.038
-7.72	-0.02
-7.68	-0.001
-7.64	0.019
-7.6	0.038
-7.56	0.057
-7.52	0.075
-7.48	0.091
-7.44	0.105
-7.4	0.117
-7.36	0.125
-7.32	0.13
-7.28	0.132
-7.24	0.13
-7.2	0.125
-7.16	0.115
-7.12	0.102
-7.08	0.086
-7.04	0.066
-7	0.045
-6.96	0.022
-6.92	-0.003
-6.88	-0.027
-6.84	-0.051
-6.8	-0.074
-6.76	-0.095
-6.72	-0.114
-6.68	-0.129
-6.64	-0.14
-6.6	-0.147
-6.56	-0.149
-6.52	-0.147
-6.48	-0.14
-6.44	-0.13
-6.4	-0.115
-6.36	-0.097
-6.32	-0.076
-6.28	-0.054
-6.24	-0.03
-6.2	-0.005
-6.16	0.019
-6.12	0.042
-6.08	0.063
-6.04	0.081
-6	0.096
-5.96	0.108
-5.92	0.116
-5.88	0.121
-5.84	0.121
-5.8	0.118
-5.76	0.111
-5.72	0.101
-5.68	0.088
-5.64	0.073
-5.6	0.057
-5.56	0.04
-5.52	0.022
-5.48	0.004
-5.44	-0.013
-5.4	-0.029
-5.36	-0.043
-5.32	-0.055
-5.28	-0.065
-5.24	-0.073
-5.2	-0.078
-5.16	-0.081
-5.12	-0.082
-5.08	-0.08
-5.04	-0.075
-5	-0.069
-4.96	-0.061
-4.92	-0.052
-4.88	-0.041
-4.84	-0.029
-4.8	-0.017
-4.76	-0.004
-4.72	0.009
-4.68	0.023
-4.64	0.035
-4.6	0.047
-4.56	0.058
-4.52	0.067
-4.48	0.075
-4.44	0.082
-4.4	0.086
-4.36	0.087
-4.32	0.087
-4.28	0.084
-4.24	0.078
-4.2	0.07
-4.16	0.06
-4.12	0.047
-4.08	0.032
-4.04	0.016
-4	-0.001
-3.96	-0.019
-3.92	-0.037
-3.88	-0.055
-3.84	-0.071
-3.8	-0.086
-3.76	-0.098
-3.72	-0.108
-3.68	-0.114
-3.64	-0.116
-3.6	-0.115
-3.56	-0.11
-3.52	-0.101
-3.48	-0.089
-3.44	-0.073
-3.4	-0.055
-3.36	-0.034
-3.32	-0.012
-3.28	0.011
-3.24	0.034
-3.2	0.056
-3.16	0.077
-3.12	0.095
-3.08	0.11
-3.04	0.122
-3	0.13
-2.96	0.133
-2.92	0.132
-2.88	0.127
-2.84	0.117
-2.8	0.104
-2.76	0.087
-2.72	0.068
-2.68	0.047
-2.64	0.024
-2.6	0.002
-2.56	-0.021
-2.52	-0.043
-2.48	-0.063
-2.44	-0.08
-2.4	-0.095
-2.36	-0.106
-2.32	-0.114
-2.28	-0.119
-2.24	-0.119
-2.2	-0.117
-2.16	-0.111
-2.12	-0.102
-2.08	-0.09
-2.04	-0.076
-2	-0.061
-1.96	-0.044
-1.92	-0.027
-1.88	-0.009
-1.84	0.008
-1.8	0.024
-1.76	0.04
-1.72	0.054
-1.68	0.066
-1.64	0.076
-1.6	0.085
-1.56	0.091
-1.52	0.095
-1.48	0.097
-1.44	0.097
-1.4	0.094
-1.36	0.089
-1.32	0.082
-1.28	0.073
-1.24	0.063
-1.2	0.05
-1.16	0.037
-1.12	0.022
-1.08	0.007
-1.04	-0.009
-1	-0.026
-0.96	-0.041
-0.92	-0.056
-0.88	-0.07
-0.84	-0.082
-0.8	-0.092
-0.76	-0.1
-0.72	-0.105
-0.68	-0.106
-0.64	-0.104
-0.6	-0.099
-0.56	-0.091
-0.52	-0.079
-0.48	-0.064
-0.44	-0.046
-0.4	-0.026
-0.36	-0.005
-0.32	0.017
-0.28	0.039
-0.24	0.061
-0.2	0.081
-0.16	0.099
-0.12	0.113
-0.08	0.124
-0.04	0.131
0	0.133
0.04	0.131
0.08	0.123
0.12	0.111
0.16	0.095
0.2	0.075
0.24	0.052
0.28	0.026
0.32	0
0.36	-0.026
0.4	-0.052
0.44	-0.076
0.48	-0.098
0.52	-0.116
0.56	-0.13
0.6	-0.139
0.64	-0.144
0.68	-0.144
0.72	-0.14
0.76	-0.131
0.8	-0.119
0.84	-0.103
0.88	-0.085
0.92	-0.064
0.96	-0.043
1	-0.021
1.04	0.001
1.08	0.021
1.12	0.041
1.16	0.058
1.2	0.074
1.24	0.087
1.28	0.098
1.32	0.106
1.36	0.111
1.4	0.114
1.44	0.114
1.48	0.112
1.52	0.107
1.56	0.101
1.6	0.092
1.64	0.082
1.68	0.069
1.72	0.055
1.76	0.04
1.8	0.023
1.84	0.005
1.88	-0.013
1.92	-0.032
1.96	-0.05
2	-0.068
2.04	-0.085
2.08	-0.1
2.12	-0.113
2.16	-0.122
2.2	-0.129
2.24	-0.131
2.28	-0.13
2.32	-0.124
2.36	-0.113
2.4	-0.098
2.44	-0.08
2.48	-0.058
2.52	-0.033
2.56	-0.006
2.6	0.022
2.64	0.049
2.68	0.076
2.72	0.101
2.76	0.123
2.8	0.141
2.84	0.154
2.88	0.162
2.92	0.164
2.96	0.161
3	0.152
3.04	0.138
3.08	0.119
3.12	0.097
3.16	0.071
3.2	0.044
3.24	0.015
3.28	-0.013
3.32	-0.04
3.36	-0.064
3.4	-0.086
3.44	-0.104
3.48	-0.118
3.52	-0.127
3.56	-0.132
3.6	-0.132
3.64	-0.128
3.68	-0.12
3.72	-0.108
3.76	-0.094
3.8	-0.078
3.84	-0.06
3.88	-0.041
3.92	-0.021
3.96	-0.002
4	0.017
4.04	0.034
4.08	0.051
4.12	0.066
4.16	0.08
4.2	0.092
4.24	0.103
4.28	0.112
4.32	0.119
4.36	0.125
4.4	0.128
4.44	0.13
4.48	0.129
4.52	0.127
4.56	0.121
4.6	0.114
4.64	0.104
4.68	0.092
4.72	0.077
4.76	0.061
4.8	0.043
4.84	0.025
4.88	0.006
4.92	-0.011
4.96	-0.027
5	-0.04
5.04	-0.049
5.08	-0.053
5.12	-0.051
5.16	-0.043
5.2	-0.027
5.24	-0.004
5.28	0.027
5.32	0.066
5.36	0.111
5.4	0.164
5.44	0.222
5.48	0.286
5.52	0.353
5.56	0.424
5.6	0.496
5.64	0.57
5.68	0.644
5.72	0.719
5.76	0.793
5.8	0.868
5.84	0.944
5.88	1.022
5.92	1.103
5.96	1.188
6	1.279
6.04	1.379
6.08	1.488
6.12	1.608
6.16	1.74
6.2	1.885
6.24	2.043
6.28	2.215
6.32	2.399
6.36	2.593
6.4	2.795
6.44	3.003
6.48	3.214
6.52	3.422
6.56	3.624
6.6	3.815
6.64	3.992
6.68	4.149
6.72	4.283
6.76	4.39
6.8	4.467
6.84	4.513
6.88	4.525
6.92	4.504
6.96	4.449
7	4.362
7.04	4.246
7.08	4.101
7.12	3.933
7.16	3.745
7.2	3.54
7.24	3.323
7.28	3.098
7.32	2.869
7.36	2.641
7.4	2.416
7.44	2.198
7.48	1.991
7.52	1.795
7.56	1.614
7.6	1.447
7.64	1.297
7.68	1.163
7.72	1.044
7.76	0.941
7.8	0.852
7.84	0.777
7.88	0.713
7.92	0.659
7.96	0.614
8	0.576
8.04	0.543
8.08	0.513
8.12	0.486
8.16	0.46
8.2	0.435
8.24	0.409
8.28	0.382
8.32	0.354
8.36	0.325
8.4	0.294
8.44	0.261
8.48	0.228
8.52	0.194
8.56	0.16
8.6	0.127
8.64	0.094
8.68	0.063
8.72	0.033
8.76	0.006
8.8	-0.019
8.84	-0.04
8.88	-0.059
8.92	-0.074
8.96	-0.085
9	-0.093
9.04	-0.097
9.08	-0.097
9.12	-0.094
9.16	-0.088
9.2	-0.078
9.24	-0.066
9.28	-0.052
9.32	-0.035
9.36	-0.017
9.4	0.002
9.44	0.021
9.48	0.041
9.52	0.059
9.56	0.077
9.6	0.092
9.64	0.106
9.68	0.116
9.72	0.124
9.76	0.128
9.8	0.128
9.84	0.125
9.88	0.118
9.92	0.108
9.96	0.094
};
\addplot[thick,color=magenta
] table {
-10	0.001
-9.96	0.001
-9.92	0.002
-9.88	0.002
-9.84	0.002
-9.8	0.003
-9.76	0.003
-9.72	0.003
-9.68	0.002
-9.64	0.002
-9.6	0.001
-9.56	0.001
-9.52	0
-9.48	0
-9.44	0
-9.4	-0.001
-9.36	0
-9.32	0
-9.28	0
-9.24	0.001
-9.2	0.001
-9.16	0.002
-9.12	0.002
-9.08	0.002
-9.04	0.002
-9	0.002
-8.96	0.002
-8.92	0.001
-8.88	0.001
-8.84	0
-8.8	0
-8.76	-0.001
-8.72	-0.001
-8.68	-0.001
-8.64	-0.001
-8.6	-0.001
-8.56	0
-8.52	0
-8.48	0
-8.44	0.001
-8.4	0.001
-8.36	0.001
-8.32	0
-8.28	0
-8.24	-0.001
-8.2	-0.001
-8.16	-0.001
-8.12	-0.002
-8.08	-0.002
-8.04	-0.002
-8	-0.002
-7.96	-0.001
-7.92	-0.001
-7.88	0
-7.84	0.001
-7.8	0.002
-7.76	0.002
-7.72	0.002
-7.68	0.003
-7.64	0.003
-7.6	0.002
-7.56	0.002
-7.52	0.001
-7.48	0.001
-7.44	0.001
-7.4	0
-7.36	0
-7.32	0
-7.28	0
-7.24	0
-7.2	0.001
-7.16	0.001
-7.12	0.001
-7.08	0.001
-7.04	0.001
-7	0
-6.96	0
-6.92	-0.001
-6.88	-0.002
-6.84	-0.002
-6.8	-0.003
-6.76	-0.003
-6.72	-0.003
-6.68	-0.003
-6.64	-0.002
-6.6	-0.002
-6.56	-0.001
-6.52	0
-6.48	0.001
-6.44	0.001
-6.4	0.002
-6.36	0.002
-6.32	0.002
-6.28	0.002
-6.24	0.002
-6.2	0.001
-6.16	0.001
-6.12	0
-6.08	0
-6.04	-0.001
-6	-0.001
-5.96	-0.001
-5.92	-0.001
-5.88	0
-5.84	0
-5.8	0
-5.76	0
-5.72	0
-5.68	0
-5.64	0
-5.6	0
-5.56	-0.001
-5.52	-0.001
-5.48	-0.002
-5.44	-0.002
-5.4	-0.002
-5.36	-0.002
-5.32	-0.002
-5.28	-0.001
-5.24	-0.001
-5.2	0
-5.16	0.001
-5.12	0.001
-5.08	0.002
-5.04	0.002
-5	0.002
-4.96	0.002
-4.92	0.001
-4.88	0
-4.84	0
-4.8	-0.001
-4.76	-0.002
-4.72	-0.002
-4.68	-0.003
-4.64	-0.003
-4.6	-0.003
-4.56	-0.003
-4.52	-0.002
-4.48	-0.001
-4.44	-0.001
-4.4	0
-4.36	0
-4.32	0
-4.28	0.001
-4.24	0
-4.2	0
-4.16	0
-4.12	0
-4.08	0
-4.04	0
-4	0
-3.96	0
-3.92	0.001
-3.88	0.001
-3.84	0.002
-3.8	0.002
-3.76	0.003
-3.72	0.003
-3.68	0.002
-3.64	0.002
-3.6	0.001
-3.56	0
-3.52	-0.001
-3.48	-0.002
-3.44	-0.003
-3.4	-0.004
-3.36	-0.004
-3.32	-0.005
-3.28	-0.005
-3.24	-0.005
-3.2	-0.004
-3.16	-0.003
-3.12	-0.003
-3.08	-0.002
-3.04	-0.001
-3	0
-2.96	0
-2.92	0
-2.88	0
-2.84	0
-2.8	0
-2.76	0
-2.72	0
-2.68	0
-2.64	0
-2.6	0
-2.56	0.001
-2.52	0.001
-2.48	0.002
-2.44	0.002
-2.4	0.002
-2.36	0.002
-2.32	0.002
-2.28	0.001
-2.24	0
-2.2	0
-2.16	-0.001
-2.12	-0.002
-2.08	-0.003
-2.04	-0.004
-2	-0.004
-1.96	-0.004
-1.92	-0.004
-1.88	-0.004
-1.84	-0.003
-1.8	-0.002
-1.76	-0.001
-1.72	0
-1.68	0
-1.64	0.001
-1.6	0.001
-1.56	0.002
-1.52	0.001
-1.48	0.001
-1.44	0.001
-1.4	0.001
-1.36	0
-1.32	0
-1.28	0
-1.24	0
-1.2	0.001
-1.16	0.001
-1.12	0.001
-1.08	0.002
-1.04	0.002
-1	0.002
-0.96	0.001
-0.92	0.001
-0.88	0
-0.84	-0.001
-0.8	-0.002
-0.76	-0.002
-0.72	-0.003
-0.68	-0.003
-0.64	-0.004
-0.6	-0.003
-0.56	-0.003
-0.52	-0.002
-0.48	-0.002
-0.44	-0.001
-0.4	0
-0.36	0.001
-0.32	0.002
-0.28	0.002
-0.24	0.002
-0.2	0.002
-0.16	0.001
-0.12	0.001
-0.08	0
-0.04	0
0	-0.001
0.04	-0.001
0.08	-0.001
0.12	-0.001
0.16	-0.001
0.2	0
0.24	0
0.28	0
0.32	0.001
0.36	0.001
0.4	0
0.44	0
0.48	-0.001
0.52	-0.001
0.56	-0.002
0.6	-0.002
0.64	-0.003
0.68	-0.003
0.72	-0.003
0.76	-0.002
0.8	-0.001
0.84	0
0.88	0.001
0.92	0.002
0.96	0.002
1	0.003
1.04	0.003
1.08	0.003
1.12	0.002
1.16	0.002
1.2	0.001
1.24	0
1.28	0
1.32	-0.001
1.36	-0.001
1.4	-0.001
1.44	-0.001
1.48	-0.001
1.52	0
1.56	0
1.6	0
1.64	0
1.68	0
1.72	0
1.76	0
1.8	-0.001
1.84	-0.002
1.88	-0.002
1.92	-0.002
1.96	-0.002
2	-0.002
2.04	-0.002
2.08	-0.001
2.12	0
2.16	0.001
2.2	0.001
2.24	0.002
2.28	0.003
2.32	0.003
2.36	0.003
2.4	0.002
2.44	0.002
2.48	0.001
2.52	0
2.56	-0.001
2.6	-0.001
2.64	-0.001
2.68	-0.001
2.72	-0.001
2.76	-0.001
2.8	0
2.84	0
2.88	0
2.92	0.001
2.96	0.001
3	0.001
3.04	0
3.08	0
3.12	-0.001
3.16	-0.001
3.2	-0.002
3.24	-0.002
3.28	-0.002
3.32	-0.001
3.36	-0.001
3.4	0
3.44	0.001
3.48	0.002
3.52	0.003
3.56	0.004
3.6	0.004
3.64	0.004
3.68	0.004
3.72	0.003
3.76	0.003
3.8	0.002
3.84	0.002
3.88	0.001
3.92	0.001
3.96	0.001
4	0.001
4.04	0.002
4.08	0.003
4.12	0.004
4.16	0.005
4.2	0.006
4.24	0.006
4.28	0.007
4.32	0.007
4.36	0.008
4.4	0.008
4.44	0.008
4.48	0.009
4.52	0.01
4.56	0.012
4.6	0.014
4.64	0.017
4.68	0.02
4.72	0.023
4.76	0.027
4.8	0.031
4.84	0.035
4.88	0.039
4.92	0.044
4.96	0.048
5	0.054
5.04	0.059
5.08	0.066
5.12	0.073
5.16	0.082
5.2	0.092
5.24	0.103
5.28	0.117
5.32	0.132
5.36	0.15
5.4	0.169
5.44	0.191
5.48	0.215
5.52	0.242
5.56	0.272
5.6	0.305
5.64	0.342
5.68	0.382
5.72	0.428
5.76	0.479
5.8	0.536
5.84	0.6
5.88	0.671
5.92	0.749
5.96	0.836
6	0.932
6.04	1.037
6.08	1.151
6.12	1.275
6.16	1.409
6.2	1.553
6.24	1.708
6.28	1.873
6.32	2.048
6.36	2.233
6.4	2.426
6.44	2.626
6.48	2.832
6.52	3.04
6.56	3.248
6.6	3.452
6.64	3.649
6.68	3.834
6.72	4.004
6.76	4.154
6.8	4.281
6.84	4.381
6.88	4.452
6.92	4.492
6.96	4.5
7	4.475
7.04	4.419
7.08	4.334
7.12	4.22
7.16	4.083
7.2	3.924
7.24	3.748
7.28	3.558
7.32	3.359
7.36	3.154
7.4	2.946
7.44	2.739
7.48	2.536
7.52	2.337
7.56	2.147
7.6	1.965
7.64	1.793
7.68	1.631
7.72	1.48
7.76	1.34
7.8	1.21
7.84	1.091
7.88	0.982
7.92	0.882
7.96	0.792
8	0.71
8.04	0.636
8.08	0.569
8.12	0.509
8.16	0.454
8.2	0.406
8.24	0.362
8.28	0.323
8.32	0.287
8.36	0.256
8.4	0.228
8.44	0.203
8.48	0.18
8.52	0.16
8.56	0.143
8.6	0.127
8.64	0.113
8.68	0.101
8.72	0.09
8.76	0.08
8.8	0.072
8.84	0.064
8.88	0.057
8.92	0.051
8.96	0.045
9	0.04
9.04	0.035
9.08	0.031
9.12	0.027
9.16	0.023
9.2	0.02
9.24	0.018
9.28	0.015
9.32	0.013
9.36	0.012
9.4	0.011
9.44	0.01
9.48	0.009
9.52	0.008
9.56	0.007
9.6	0.007
9.64	0.006
9.68	0.005
9.72	0.004
9.76	0.003
9.8	0.003
9.84	0.002
9.88	0.002
9.92	0.001
9.96	0.001
};
\end{axis}
\end{tikzpicture}

\caption{Numerical solutions $V\bz$ at $t=3$ for $r=20,40,60$ for the KdV equation.
}
\label{fig:kdv1:sol}
\end{figure}

%% file: figmkdv1svd.tex
\begin{figure}[htbp]
\centering

\begin{tikzpicture}
\tikzstyle{every node}=[]
\begin{semilogyaxis}[width=6cm,
xmax=500,xmin=0,
xlabel={order of modes},ylabel={singlura values of $Y$},
ylabel near ticks,
	]
\addplot[thick,
] table {
1	2.11E+02
2	1.66E+02
3	1.61E+02
4	1.24E+02
5	1.22E+02
6	8.29E+01
7	8.01E+01
8	5.15E+01
9	5.12E+01
10	3.20E+01
11	3.15E+01
12	1.96E+01
13	1.93E+01
14	1.20E+01
15	1.19E+01
16	7.33E+00
17	7.23E+00
18	4.45E+00
19	4.45E+00
20	2.73E+00
21	2.71E+00
22	1.67E+00
23	1.66E+00
24	1.02E+00
25	1.01E+00
26	6.23E-01
27	6.18E-01
28	3.79E-01
29	3.79E-01
30	2.33E-01
31	2.31E-01
32	1.42E-01
33	1.41E-01
34	8.68E-02
35	8.64E-02
36	5.31E-02
37	5.28E-02
38	3.24E-02
39	3.24E-02
40	1.99E-02
41	1.97E-02
42	1.21E-02
43	1.21E-02
44	7.42E-03
45	7.40E-03
46	4.54E-03
47	4.52E-03
48	2.78E-03
49	2.78E-03
50	1.71E-03
51	1.70E-03
52	1.05E-03
53	1.05E-03
54	6.53E-04
55	6.51E-04
56	4.17E-04
57	4.15E-04
58	2.72E-04
59	2.72E-04
60	1.96E-04
61	1.96E-04
62	1.51E-04
63	1.50E-04
64	1.27E-04
65	1.26E-04
66	1.08E-04
67	1.06E-04
68	1.06E-04
69	1.02E-04
70	9.83E-05
71	9.81E-05
72	9.79E-05
73	9.70E-05
74	9.62E-05
75	9.59E-05
76	9.57E-05
77	9.57E-05
78	9.54E-05
79	9.53E-05
80	9.50E-05
81	9.49E-05
82	9.48E-05
83	9.47E-05
84	9.46E-05
85	9.45E-05
86	9.44E-05
87	9.44E-05
88	9.36E-05
89	9.34E-05
90	9.26E-05
91	9.26E-05
92	8.99E-05
93	8.13E-05
94	6.37E-05
95	2.70E-05
96	7.30E-06
97	3.75E-06
98	3.66E-06
99	2.31E-06
100	1.37E-06
101	8.58E-07
102	5.23E-07
103	3.89E-07
104	2.51E-07
105	1.87E-07
106	1.66E-07
107	1.60E-07
108	1.33E-07
109	9.30E-08
110	7.27E-08
111	6.12E-08
112	4.36E-08
113	3.47E-08
114	2.57E-08
115	1.92E-08
116	1.41E-08
117	1.18E-08
118	1.17E-08
119	1.00E-08
120	7.45E-09
121	5.43E-09
122	4.05E-09
123	3.23E-09
124	2.79E-09
125	2.67E-09
126	2.35E-09
127	1.90E-09
128	1.43E-09
129	1.22E-09
130	8.99E-10
131	8.66E-10
132	8.62E-10
133	6.63E-10
134	5.21E-10
135	4.88E-10
136	3.69E-10
137	3.21E-10
138	2.88E-10
139	2.58E-10
140	2.28E-10
141	1.75E-10
142	1.46E-10
143	1.30E-10
144	1.23E-10
145	1.13E-10
146	8.93E-11
147	7.23E-11
148	6.13E-11
149	5.83E-11
150	5.22E-11
151	3.98E-11
152	3.60E-11
153	3.37E-11
154	3.17E-11
155	2.50E-11
156	2.12E-11
157	1.98E-11
158	1.63E-11
159	1.51E-11
160	1.34E-11
161	1.18E-11
162	1.09E-11
163	8.69E-12
164	8.33E-12
165	7.12E-12
166	6.71E-12
167	5.68E-12
168	5.16E-12
169	4.45E-12
170	4.06E-12
171	3.50E-12
172	3.16E-12
173	2.89E-12
174	2.45E-12
175	2.36E-12
176	2.23E-12
177	2.01E-12
178	1.85E-12
179	1.74E-12
180	1.63E-12
181	1.58E-12
182	1.46E-12
183	1.43E-12
184	1.34E-12
185	1.17E-12
186	1.14E-12
187	1.06E-12
188	9.65E-13
189	8.63E-13
190	7.93E-13
191	7.52E-13
192	7.08E-13
193	6.57E-13
194	6.15E-13
195	5.28E-13
196	4.43E-13
197	4.27E-13
198	3.61E-13
199	3.29E-13
200	2.81E-13
201	2.78E-13
202	2.23E-13
203	1.90E-13
204	1.81E-13
205	1.70E-13
206	1.44E-13
207	1.27E-13
208	1.24E-13
209	1.10E-13
210	9.64E-14
211	8.65E-14
212	8.15E-14
213	7.71E-14
214	7.39E-14
215	6.98E-14
216	5.81E-14
217	5.29E-14
218	5.13E-14
219	4.79E-14
220	4.48E-14
221	4.10E-14
222	4.09E-14
223	3.76E-14
224	3.47E-14
225	3.34E-14
226	3.06E-14
227	2.97E-14
228	2.97E-14
229	2.81E-14
230	2.50E-14
231	2.39E-14
232	2.12E-14
233	1.97E-14
234	1.88E-14
235	1.84E-14
236	1.71E-14
237	1.62E-14
238	1.43E-14
239	1.43E-14
240	1.43E-14
241	1.43E-14
242	1.43E-14
243	1.43E-14
244	1.43E-14
245	1.43E-14
246	1.43E-14
247	1.43E-14
248	1.43E-14
249	1.43E-14
250	1.43E-14
251	1.43E-14
252	1.43E-14
253	1.43E-14
254	1.43E-14
255	1.43E-14
256	1.43E-14
257	1.43E-14
258	1.43E-14
259	1.43E-14
260	1.43E-14
261	1.43E-14
262	1.43E-14
263	1.43E-14
264	1.43E-14
265	1.43E-14
266	1.43E-14
267	1.43E-14
268	1.43E-14
269	1.43E-14
270	1.43E-14
271	1.43E-14
272	1.43E-14
273	1.43E-14
274	1.43E-14
275	1.43E-14
276	1.43E-14
277	1.43E-14
278	1.43E-14
279	1.43E-14
280	1.43E-14
281	1.43E-14
282	1.43E-14
283	1.43E-14
284	1.43E-14
285	1.43E-14
286	1.43E-14
287	1.43E-14
288	1.43E-14
289	1.43E-14
290	1.43E-14
291	1.43E-14
292	1.43E-14
293	1.43E-14
294	1.43E-14
295	1.43E-14
296	1.43E-14
297	1.43E-14
298	1.43E-14
299	1.43E-14
300	1.43E-14
301	1.43E-14
302	1.43E-14
303	1.43E-14
304	1.43E-14
305	1.43E-14
306	1.43E-14
307	1.43E-14
308	1.43E-14
309	1.43E-14
310	1.43E-14
311	1.43E-14
312	1.43E-14
313	1.43E-14
314	1.43E-14
315	1.43E-14
316	1.43E-14
317	1.43E-14
318	1.43E-14
319	1.43E-14
320	1.43E-14
321	1.43E-14
322	1.43E-14
323	1.43E-14
324	1.43E-14
325	1.43E-14
326	1.43E-14
327	1.43E-14
328	1.43E-14
329	1.43E-14
330	1.43E-14
331	1.43E-14
332	1.43E-14
333	1.43E-14
334	1.43E-14
335	1.43E-14
336	1.43E-14
337	1.43E-14
338	1.43E-14
339	1.43E-14
340	1.43E-14
341	1.43E-14
342	1.43E-14
343	1.43E-14
344	1.43E-14
345	1.43E-14
346	1.43E-14
347	1.43E-14
348	1.43E-14
349	1.43E-14
350	1.43E-14
351	1.43E-14
352	1.43E-14
353	1.43E-14
354	1.43E-14
355	1.43E-14
356	1.43E-14
357	1.43E-14
358	1.43E-14
359	1.43E-14
360	1.43E-14
361	1.43E-14
362	1.43E-14
363	1.43E-14
364	1.43E-14
365	1.43E-14
366	1.43E-14
367	1.43E-14
368	1.43E-14
369	1.43E-14
370	1.43E-14
371	1.43E-14
372	1.43E-14
373	1.43E-14
374	1.43E-14
375	1.43E-14
376	1.43E-14
377	1.43E-14
378	1.43E-14
379	1.43E-14
380	1.43E-14
381	1.43E-14
382	1.43E-14
383	1.43E-14
384	1.43E-14
385	1.43E-14
386	1.43E-14
387	1.43E-14
388	1.43E-14
389	1.43E-14
390	1.43E-14
391	1.43E-14
392	1.43E-14
393	1.43E-14
394	1.43E-14
395	1.43E-14
396	1.43E-14
397	1.43E-14
398	1.43E-14
399	1.43E-14
400	1.43E-14
401	1.43E-14
402	1.43E-14
403	1.43E-14
404	1.43E-14
405	1.43E-14
406	1.43E-14
407	1.43E-14
408	1.43E-14
409	1.43E-14
410	1.43E-14
411	1.43E-14
412	1.43E-14
413	1.43E-14
414	1.43E-14
415	1.43E-14
416	1.43E-14
417	1.43E-14
418	1.43E-14
419	1.43E-14
420	1.43E-14
421	1.43E-14
422	1.43E-14
423	1.43E-14
424	1.43E-14
425	1.43E-14
426	1.43E-14
427	1.43E-14
428	1.43E-14
429	1.43E-14
430	1.43E-14
431	1.43E-14
432	1.43E-14
433	1.43E-14
434	1.43E-14
435	1.43E-14
436	1.43E-14
437	1.43E-14
438	1.43E-14
439	1.43E-14
440	1.43E-14
441	1.43E-14
442	1.43E-14
443	1.43E-14
444	1.43E-14
445	1.43E-14
446	1.43E-14
447	1.43E-14
448	1.43E-14
449	1.43E-14
450	1.43E-14
451	1.43E-14
452	1.43E-14
453	1.43E-14
454	1.43E-14
455	1.43E-14
456	1.43E-14
457	1.43E-14
458	1.43E-14
459	1.43E-14
460	1.43E-14
461	1.43E-14
462	1.43E-14
463	1.43E-14
464	1.43E-14
465	1.43E-14
466	1.43E-14
467	1.43E-14
468	1.43E-14
469	1.43E-14
470	1.43E-14
471	1.43E-14
472	1.43E-14
473	1.43E-14
474	1.43E-14
475	1.43E-14
476	1.43E-14
477	1.43E-14
478	1.43E-14
479	1.43E-14
480	1.43E-14
481	1.43E-14
482	1.43E-14
483	1.43E-14
484	1.43E-14
485	1.43E-14
486	1.43E-14
487	1.43E-14
488	1.43E-14
489	1.43E-14
490	1.43E-14
491	1.43E-14
492	1.43E-14
493	1.43E-14
494	1.26E-14
495	9.89E-15
496	8.59E-15
497	8.22E-15
498	7.23E-15
499	6.23E-15
500	3.49E-15
};
\end{semilogyaxis}
\end{tikzpicture} 
\quad 
\begin{tikzpicture}
\tikzstyle{every node}=[]
\begin{semilogyaxis}[width=6cm,
xmax=500,xmin=0,
xlabel={order of modes},ylabel={singlura values of $S$},
ylabel near ticks,
	]
\addplot[thick,
] table {
1	1.21E+08
2	2.73E+04
3	2.45E+04
4	2.18E+04
5	2.12E+04
6	1.81E+04
7	1.77E+04
8	1.42E+04
9	1.41E+04
10	1.08E+04
11	1.06E+04
12	7.84E+03
13	7.77E+03
14	5.58E+03
15	5.53E+03
16	3.90E+03
17	3.85E+03
18	2.66E+03
19	2.66E+03
20	1.81E+03
21	1.80E+03
22	1.22E+03
23	1.21E+03
24	8.10E+02
25	8.05E+02
26	5.36E+02
27	5.32E+02
28	3.51E+02
29	3.51E+02
30	2.31E+02
31	2.29E+02
32	1.50E+02
33	1.49E+02
34	9.73E+01
35	9.69E+01
36	6.30E+01
37	6.26E+01
38	4.05E+01
39	4.05E+01
40	2.61E+01
41	2.59E+01
42	1.67E+01
43	1.67E+01
44	1.07E+01
45	1.07E+01
46	6.83E+00
47	6.80E+00
48	4.35E+00
49	4.34E+00
50	2.77E+00
51	2.76E+00
52	1.76E+00
53	1.75E+00
54	1.12E+00
55	1.11E+00
56	7.07E-01
57	7.04E-01
58	4.47E-01
59	4.46E-01
60	2.83E-01
61	2.82E-01
62	1.79E-01
63	1.79E-01
64	1.13E-01
65	1.13E-01
66	7.21E-02
67	7.19E-02
68	4.64E-02
69	4.63E-02
70	3.10E-02
71	3.09E-02
72	2.70E-02
73	2.41E-02
74	2.24E-02
75	2.05E-02
76	1.86E-02
77	1.71E-02
78	1.60E-02
79	1.49E-02
80	1.34E-02
81	1.24E-02
82	1.07E-02
83	9.64E-03
84	8.51E-03
85	6.87E-03
86	6.29E-03
87	5.38E-03
88	4.26E-03
89	3.76E-03
90	3.31E-03
91	2.64E-03
92	2.09E-03
93	2.03E-03
94	1.60E-03
95	1.29E-03
96	1.09E-03
97	9.82E-04
98	7.57E-04
99	6.36E-04
100	5.41E-04
101	4.63E-04
102	3.51E-04
103	3.09E-04
104	2.54E-04
105	2.12E-04
106	1.59E-04
107	1.40E-04
108	1.12E-04
109	9.58E-05
110	7.16E-05
111	6.27E-05
112	5.19E-05
113	4.48E-05
114	4.11E-05
115	3.77E-05
116	3.25E-05
117	2.88E-05
118	2.57E-05
119	2.26E-05
120	1.99E-05
121	1.69E-05
122	1.36E-05
123	1.21E-05
124	1.13E-05
125	9.41E-06
126	7.25E-06
127	7.07E-06
128	6.03E-06
129	5.29E-06
130	4.51E-06
131	4.10E-06
132	3.68E-06
133	3.41E-06
134	2.99E-06
135	2.88E-06
136	2.64E-06
137	2.21E-06
138	1.96E-06
139	1.78E-06
140	1.54E-06
141	1.28E-06
142	1.23E-06
143	1.05E-06
144	9.22E-07
145	8.57E-07
146	7.86E-07
147	7.14E-07
148	5.59E-07
149	5.14E-07
150	4.60E-07
151	4.10E-07
152	3.42E-07
153	3.04E-07
154	2.85E-07
155	2.48E-07
156	2.28E-07
157	2.13E-07
158	1.95E-07
159	1.73E-07
160	1.66E-07
161	1.45E-07
162	1.38E-07
163	1.29E-07
164	1.21E-07
165	1.16E-07
166	1.12E-07
167	1.02E-07
168	9.62E-08
169	9.34E-08
170	8.92E-08
171	7.74E-08
172	7.39E-08
173	6.94E-08
174	6.56E-08
175	5.62E-08
176	4.85E-08
177	4.44E-08
178	4.12E-08
179	3.86E-08
180	3.79E-08
181	3.44E-08
182	3.10E-08
183	2.94E-08
184	2.62E-08
185	2.53E-08
186	2.18E-08
187	1.96E-08
188	1.69E-08
189	1.31E-08
190	1.25E-08
191	1.21E-08
192	1.21E-08
193	1.21E-08
194	1.21E-08
195	1.21E-08
196	1.21E-08
197	1.21E-08
198	1.21E-08
199	1.21E-08
200	1.21E-08
201	1.21E-08
202	1.21E-08
203	1.21E-08
204	1.21E-08
205	1.21E-08
206	1.21E-08
207	1.21E-08
208	1.21E-08
209	1.21E-08
210	1.21E-08
211	1.21E-08
212	1.21E-08
213	1.21E-08
214	1.21E-08
215	1.21E-08
216	1.21E-08
217	1.21E-08
218	1.21E-08
219	1.21E-08
220	1.21E-08
221	1.21E-08
222	1.21E-08
223	1.21E-08
224	1.21E-08
225	1.21E-08
226	1.21E-08
227	1.21E-08
228	1.21E-08
229	1.21E-08
230	1.21E-08
231	1.21E-08
232	1.21E-08
233	1.21E-08
234	1.21E-08
235	1.21E-08
236	1.21E-08
237	1.21E-08
238	1.21E-08
239	1.21E-08
240	1.21E-08
241	1.21E-08
242	1.21E-08
243	1.21E-08
244	1.21E-08
245	1.21E-08
246	1.21E-08
247	1.21E-08
248	1.21E-08
249	1.21E-08
250	1.21E-08
251	1.21E-08
252	1.21E-08
253	1.21E-08
254	1.21E-08
255	1.21E-08
256	1.21E-08
257	1.21E-08
258	1.21E-08
259	1.21E-08
260	1.21E-08
261	1.21E-08
262	1.21E-08
263	1.21E-08
264	1.21E-08
265	1.21E-08
266	1.21E-08
267	1.21E-08
268	1.21E-08
269	1.21E-08
270	1.21E-08
271	1.21E-08
272	1.21E-08
273	1.21E-08
274	1.21E-08
275	1.21E-08
276	1.21E-08
277	1.21E-08
278	1.21E-08
279	1.21E-08
280	1.21E-08
281	1.21E-08
282	1.21E-08
283	1.21E-08
284	1.21E-08
285	1.21E-08
286	1.21E-08
287	1.21E-08
288	1.21E-08
289	1.21E-08
290	1.21E-08
291	1.21E-08
292	1.21E-08
293	1.21E-08
294	1.21E-08
295	1.21E-08
296	1.21E-08
297	1.21E-08
298	1.21E-08
299	1.21E-08
300	1.21E-08
301	1.21E-08
302	1.21E-08
303	1.21E-08
304	1.21E-08
305	1.21E-08
306	1.21E-08
307	1.21E-08
308	1.21E-08
309	1.21E-08
310	1.21E-08
311	1.21E-08
312	1.21E-08
313	1.21E-08
314	1.21E-08
315	1.21E-08
316	1.21E-08
317	1.21E-08
318	1.21E-08
319	1.21E-08
320	1.21E-08
321	1.21E-08
322	1.21E-08
323	1.21E-08
324	1.21E-08
325	1.21E-08
326	1.21E-08
327	1.21E-08
328	1.21E-08
329	1.21E-08
330	1.21E-08
331	1.21E-08
332	1.21E-08
333	1.21E-08
334	1.21E-08
335	1.21E-08
336	1.21E-08
337	1.21E-08
338	1.21E-08
339	1.21E-08
340	1.21E-08
341	1.21E-08
342	1.21E-08
343	1.21E-08
344	1.21E-08
345	1.21E-08
346	1.21E-08
347	1.21E-08
348	1.21E-08
349	1.21E-08
350	1.21E-08
351	1.21E-08
352	1.21E-08
353	1.21E-08
354	1.21E-08
355	1.21E-08
356	1.21E-08
357	1.21E-08
358	1.21E-08
359	1.21E-08
360	1.21E-08
361	1.21E-08
362	1.21E-08
363	1.21E-08
364	1.21E-08
365	1.21E-08
366	1.21E-08
367	1.21E-08
368	1.21E-08
369	1.21E-08
370	1.21E-08
371	1.21E-08
372	1.21E-08
373	1.21E-08
374	1.21E-08
375	1.21E-08
376	1.21E-08
377	1.21E-08
378	1.21E-08
379	1.21E-08
380	1.21E-08
381	1.21E-08
382	1.21E-08
383	1.21E-08
384	1.21E-08
385	1.21E-08
386	1.21E-08
387	1.21E-08
388	1.21E-08
389	1.21E-08
390	1.21E-08
391	1.21E-08
392	1.21E-08
393	1.21E-08
394	1.21E-08
395	1.21E-08
396	1.21E-08
397	1.21E-08
398	1.21E-08
399	1.21E-08
400	1.21E-08
401	1.21E-08
402	1.21E-08
403	1.21E-08
404	1.21E-08
405	1.21E-08
406	1.21E-08
407	1.21E-08
408	1.21E-08
409	1.21E-08
410	1.21E-08
411	1.21E-08
412	1.21E-08
413	1.21E-08
414	1.21E-08
415	1.21E-08
416	1.21E-08
417	1.21E-08
418	1.21E-08
419	1.21E-08
420	1.21E-08
421	1.21E-08
422	1.21E-08
423	1.21E-08
424	1.21E-08
425	1.21E-08
426	1.21E-08
427	1.21E-08
428	1.21E-08
429	1.21E-08
430	1.21E-08
431	1.21E-08
432	1.21E-08
433	1.21E-08
434	1.21E-08
435	1.21E-08
436	1.21E-08
437	1.21E-08
438	1.21E-08
439	1.21E-08
440	1.21E-08
441	1.21E-08
442	1.21E-08
443	1.21E-08
444	1.21E-08
445	1.21E-08
446	1.21E-08
447	1.21E-08
448	1.21E-08
449	1.21E-08
450	1.21E-08
451	1.21E-08
452	1.21E-08
453	1.21E-08
454	1.21E-08
455	1.21E-08
456	1.21E-08
457	1.21E-08
458	1.21E-08
459	1.21E-08
460	1.21E-08
461	1.21E-08
462	1.21E-08
463	1.21E-08
464	1.21E-08
465	1.21E-08
466	1.21E-08
467	1.21E-08
468	1.21E-08
469	1.21E-08
470	1.21E-08
471	1.21E-08
472	1.21E-08
473	1.21E-08
474	1.21E-08
475	1.21E-08
476	1.21E-08
477	1.21E-08
478	1.21E-08
479	1.21E-08
480	1.21E-08
481	1.21E-08
482	1.21E-08
483	1.21E-08
484	1.21E-08
485	1.21E-08
486	1.21E-08
487	1.21E-08
488	1.21E-08
489	1.21E-08
490	1.21E-08
491	1.21E-08
492	1.21E-08
493	1.21E-08
494	1.21E-08
495	1.21E-08
496	1.21E-08
497	1.21E-08
498	1.21E-08
499	1.21E-08
500	1.21E-08
501	1.21E-08
502	1.21E-08
503	1.21E-08
504	1.21E-08
505	1.21E-08
506	1.21E-08
507	1.21E-08
508	1.21E-08
509	1.21E-08
510	1.21E-08
511	1.21E-08
512	1.21E-08
513	1.21E-08
514	1.21E-08
515	1.21E-08
516	1.21E-08
517	1.21E-08
518	1.21E-08
519	1.21E-08
520	1.21E-08
521	1.21E-08
522	1.21E-08
523	1.21E-08
524	1.21E-08
525	1.21E-08
526	1.21E-08
527	1.21E-08
528	1.21E-08
529	1.21E-08
530	1.21E-08
531	1.21E-08
532	1.21E-08
533	1.21E-08
534	1.21E-08
535	1.21E-08
536	1.21E-08
537	1.21E-08
538	1.21E-08
539	1.21E-08
540	1.21E-08
541	1.21E-08
542	1.21E-08
543	1.21E-08
544	1.21E-08
545	1.21E-08
546	1.21E-08
547	1.21E-08
548	1.21E-08
549	1.21E-08
550	1.21E-08
551	1.21E-08
552	1.21E-08
553	1.21E-08
554	1.21E-08
555	1.21E-08
556	1.21E-08
557	1.21E-08
558	1.21E-08
559	1.21E-08
560	1.21E-08
561	1.21E-08
562	1.21E-08
563	1.21E-08
564	1.21E-08
565	1.21E-08
566	1.21E-08
567	1.21E-08
568	1.21E-08
569	1.21E-08
570	1.21E-08
571	1.21E-08
572	1.21E-08
573	1.21E-08
574	1.21E-08
575	1.21E-08
576	1.21E-08
577	1.21E-08
578	1.21E-08
579	1.21E-08
580	1.21E-08
581	1.21E-08
582	1.21E-08
583	1.21E-08
584	1.21E-08
585	1.21E-08
586	1.21E-08
587	1.21E-08
588	1.21E-08
589	1.21E-08
590	1.21E-08
591	1.21E-08
592	1.21E-08
593	1.21E-08
594	1.21E-08
595	1.21E-08
596	1.21E-08
597	1.21E-08
598	1.21E-08
599	1.21E-08
600	1.21E-08
601	1.21E-08
602	1.21E-08
603	1.21E-08
604	1.21E-08
605	1.21E-08
606	1.21E-08
607	1.21E-08
608	1.21E-08
609	1.21E-08
610	1.21E-08
611	1.21E-08
612	1.21E-08
613	1.21E-08
614	1.21E-08
615	1.21E-08
616	1.21E-08
617	1.21E-08
618	1.21E-08
619	1.21E-08
620	1.21E-08
621	1.21E-08
622	1.21E-08
623	1.21E-08
624	1.21E-08
625	1.21E-08
626	1.21E-08
627	1.21E-08
628	1.21E-08
629	1.21E-08
630	1.21E-08
631	1.21E-08
632	1.21E-08
633	1.21E-08
634	1.21E-08
635	1.21E-08
636	1.21E-08
637	1.21E-08
638	1.21E-08
639	1.21E-08
640	1.21E-08
641	1.21E-08
642	1.21E-08
643	1.21E-08
644	1.21E-08
645	1.21E-08
646	1.21E-08
647	1.21E-08
648	1.21E-08
649	1.21E-08
650	1.21E-08
651	1.21E-08
652	1.21E-08
653	1.21E-08
654	1.21E-08
655	1.21E-08
656	1.21E-08
657	1.21E-08
658	1.21E-08
659	1.21E-08
660	1.21E-08
661	1.21E-08
662	1.21E-08
663	1.21E-08
664	1.21E-08
665	1.21E-08
666	1.21E-08
667	1.21E-08
668	1.21E-08
669	1.21E-08
670	1.21E-08
671	1.21E-08
672	1.21E-08
673	1.21E-08
674	1.21E-08
675	1.21E-08
676	1.21E-08
677	1.21E-08
678	1.21E-08
679	1.21E-08
680	1.21E-08
681	1.21E-08
682	1.21E-08
683	1.21E-08
684	1.21E-08
685	1.21E-08
686	1.21E-08
687	1.21E-08
688	1.21E-08
689	1.21E-08
690	1.21E-08
691	1.21E-08
692	1.21E-08
693	1.21E-08
694	1.21E-08
695	1.21E-08
696	1.21E-08
697	1.21E-08
698	1.21E-08
699	1.21E-08
700	1.21E-08
701	1.21E-08
702	1.21E-08
703	1.21E-08
704	1.21E-08
705	1.21E-08
706	1.21E-08
707	1.21E-08
708	1.21E-08
709	1.21E-08
710	1.21E-08
711	1.21E-08
712	1.21E-08
713	1.21E-08
714	1.21E-08
715	1.21E-08
716	1.21E-08
717	1.21E-08
718	1.21E-08
719	1.21E-08
720	1.21E-08
721	1.21E-08
722	1.21E-08
723	1.21E-08
724	1.21E-08
725	1.21E-08
726	1.21E-08
727	1.21E-08
728	1.21E-08
729	1.21E-08
730	1.21E-08
731	1.21E-08
732	1.21E-08
733	1.21E-08
734	1.21E-08
735	1.21E-08
736	1.21E-08
737	1.21E-08
738	1.21E-08
739	1.21E-08
740	1.21E-08
741	1.21E-08
742	1.21E-08
743	1.21E-08
744	1.21E-08
745	1.21E-08
746	1.21E-08
747	1.11E-08
748	7.59E-09
749	7.46E-09
750	4.86E-09
751	8.58E-10
};
\end{semilogyaxis}
\end{tikzpicture}

\caption{Singular values corresponding to the POD modes for (left) $Y$ and (right) $S=[\veco (S(\by_0)), \veco (S(\by_1)),\dots,\veco (S(\by_{750}))]$ for the mKdV equation.
}
\label{fig:mkdv1:svd}
\end{figure}
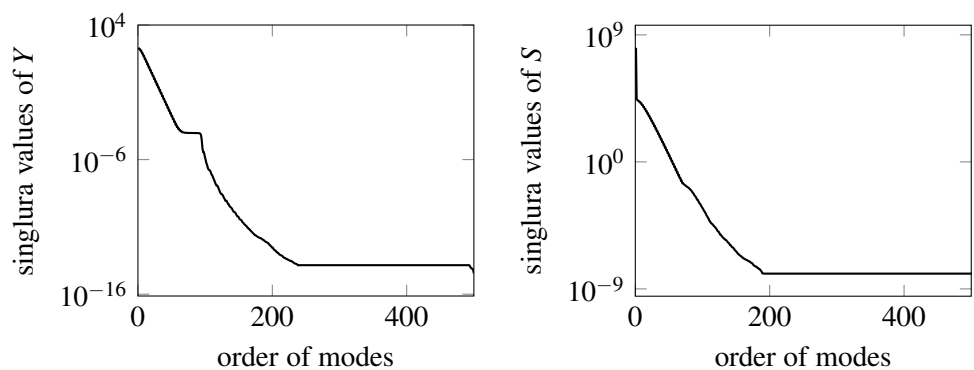

%% file: figmkdv1Eerror.tex
\begin{figure}[htbp]
\centering

\begin{tikzpicture}
\tikzstyle{every node}=[]
\begin{semilogyaxis}[width=6cm,
xmax=10,xmin=0,
xlabel={$t$},ylabel={error of $H(V\bz)$},
ylabel near ticks,
	]
\addplot[thick,
] table {
0	0
0.004	5.83E-16
0.008	1.14E-15
0.012	1.78E-15
0.016	2.39E-15
0.02	2.97E-15
0.024	3.72E-15
0.028	4.36E-15
0.032	5.08E-15
0.036	5.69E-15
0.04	6.41E-15
0.044	7.16E-15
0.048	7.91E-15
0.052	8.66E-15
0.056	9.30E-15
0.06	1.02E-14
0.064	1.09E-14
0.068	1.18E-14
0.072	1.26E-14
0.076	1.34E-14
0.08	1.43E-14
0.084	1.52E-14
0.088	1.61E-14
0.092	1.71E-14
0.096	1.81E-14
0.1	1.90E-14
0.104	1.99E-14
0.108	2.09E-14
0.112	2.18E-14
0.116	2.27E-14
0.12	2.35E-14
0.124	2.45E-14
0.128	2.52E-14
0.132	2.61E-14
0.136	2.69E-14
0.14	2.76E-14
0.144	2.83E-14
0.148	2.90E-14
0.152	2.95E-14
0.156	3.00E-14
0.16	3.08E-14
0.164	3.12E-14
0.168	3.16E-14
0.172	3.21E-14
0.176	3.24E-14
0.18	3.26E-14
0.184	3.30E-14
0.188	3.33E-14
0.192	3.33E-14
0.196	3.34E-14
0.2	3.34E-14
0.204	3.34E-14
0.208	3.33E-14
0.212	3.32E-14
0.216	3.29E-14
0.22	3.27E-14
0.224	3.23E-14
0.228	3.19E-14
0.232	3.15E-14
0.236	3.10E-14
0.24	3.04E-14
0.244	2.97E-14
0.248	2.91E-14
0.252	2.85E-14
0.256	2.77E-14
0.26	2.70E-14
0.264	2.61E-14
0.268	2.53E-14
0.272	2.43E-14
0.276	2.35E-14
0.28	2.25E-14
0.284	2.16E-14
0.288	2.05E-14
0.292	1.97E-14
0.296	1.87E-14
0.3	1.77E-14
0.304	1.66E-14
0.308	1.55E-14
0.312	1.47E-14
0.316	1.35E-14
0.32	1.26E-14
0.324	1.16E-14
0.328	1.05E-14
0.332	9.55E-15
0.336	8.69E-15
0.34	7.72E-15
0.344	6.72E-15
0.348	5.83E-15
0.352	4.80E-15
0.356	3.91E-15
0.36	3.14E-15
0.364	2.08E-15
0.368	1.19E-15
0.372	3.05E-16
0.376	4.72E-16
0.38	1.33E-15
0.384	2.25E-15
0.388	3.08E-15
0.392	3.94E-15
0.396	4.77E-15
0.4	5.52E-15
0.404	6.38E-15
0.408	7.24E-15
0.412	7.91E-15
0.416	8.77E-15
0.42	9.44E-15
0.424	1.02E-14
0.428	1.09E-14
0.432	1.16E-14
0.436	1.24E-14
0.44	1.30E-14
0.444	1.36E-14
0.448	1.43E-14
0.452	1.50E-14
0.456	1.55E-14
0.46	1.62E-14
0.464	1.67E-14
0.468	1.72E-14
0.472	1.78E-14
0.476	1.83E-14
0.48	1.88E-14
0.484	1.93E-14
0.488	1.99E-14
0.492	2.05E-14
0.496	2.10E-14
0.5	2.14E-14
0.504	2.19E-14
0.508	2.24E-14
0.512	2.29E-14
0.516	2.33E-14
0.52	2.39E-14
0.524	2.44E-14
0.528	2.49E-14
0.532	2.53E-14
0.536	2.58E-14
0.54	2.63E-14
0.544	2.69E-14
0.548	2.72E-14
0.552	2.78E-14
0.556	2.84E-14
0.56	2.88E-14
0.564	2.92E-14
0.568	2.98E-14
0.572	3.02E-14
0.576	3.09E-14
0.58	3.15E-14
0.584	3.20E-14
0.588	3.24E-14
0.592	3.31E-14
0.596	3.37E-14
0.6	3.43E-14
0.604	3.49E-14
0.608	3.55E-14
0.612	3.60E-14
0.616	3.66E-14
0.62	3.73E-14
0.624	3.78E-14
0.628	3.84E-14
0.632	3.89E-14
0.636	3.95E-14
0.64	4.01E-14
0.644	4.06E-14
0.648	4.12E-14
0.652	4.18E-14
0.656	4.23E-14
0.66	4.28E-14
0.664	4.33E-14
0.668	4.37E-14
0.672	4.43E-14
0.676	4.47E-14
0.68	4.51E-14
0.684	4.57E-14
0.688	4.61E-14
0.692	4.65E-14
0.696	4.69E-14
0.7	4.74E-14
0.704	4.78E-14
0.708	4.82E-14
0.712	4.85E-14
0.716	4.89E-14
0.72	4.91E-14
0.724	4.95E-14
0.728	4.98E-14
0.732	5.01E-14
0.736	5.04E-14
0.74	5.08E-14
0.744	5.12E-14
0.748	5.15E-14
0.752	5.18E-14
0.756	5.22E-14
0.76	5.26E-14
0.764	5.29E-14
0.768	5.33E-14
0.772	5.36E-14
0.776	5.40E-14
0.78	5.43E-14
0.784	5.47E-14
0.788	5.51E-14
0.792	5.54E-14
0.796	5.58E-14
0.8	5.62E-14
0.804	5.67E-14
0.808	5.71E-14
0.812	5.76E-14
0.816	5.82E-14
0.82	5.88E-14
0.824	5.95E-14
0.828	6.02E-14
0.832	6.10E-14
0.836	6.18E-14
0.84	6.27E-14
0.844	6.37E-14
0.848	6.47E-14
0.852	6.58E-14
0.856	6.69E-14
0.86	6.81E-14
0.864	6.93E-14
0.868	7.07E-14
0.872	7.20E-14
0.876	7.33E-14
0.88	7.47E-14
0.884	7.62E-14
0.888	7.77E-14
0.892	7.92E-14
0.896	8.08E-14
0.9	8.27E-14
0.904	8.43E-14
0.908	8.64E-14
0.912	8.85E-14
0.916	9.09E-14
0.92	9.35E-14
0.924	9.63E-14
0.928	9.94E-14
0.932	1.03E-13
0.936	1.06E-13
0.94	1.10E-13
0.944	1.15E-13
0.948	1.19E-13
0.952	1.24E-13
0.956	1.30E-13
0.96	1.35E-13
0.964	1.41E-13
0.968	1.47E-13
0.972	1.53E-13
0.976	1.59E-13
0.98	1.65E-13
0.984	1.72E-13
0.988	1.79E-13
0.992	1.87E-13
0.996	1.94E-13
1	2.03E-13
1.004	2.11E-13
1.008	2.21E-13
1.012	2.31E-13
1.016	2.42E-13
1.02	2.54E-13
1.024	2.66E-13
1.028	2.80E-13
1.032	2.95E-13
1.036	3.12E-13
1.04	3.29E-13
1.044	3.48E-13
1.048	3.69E-13
1.052	3.90E-13
1.056	4.13E-13
1.06	4.37E-13
1.064	4.63E-13
1.068	4.89E-13
1.072	5.17E-13
1.076	5.46E-13
1.08	5.76E-13
1.084	6.07E-13
1.088	6.40E-13
1.092	6.73E-13
1.096	7.09E-13
1.1	7.46E-13
1.104	7.84E-13
1.108	8.25E-13
1.112	8.67E-13
1.116	9.13E-13
1.12	9.60E-13
1.124	1.01E-12
1.128	1.06E-12
1.132	1.12E-12
1.136	1.18E-12
1.14	1.24E-12
1.144	1.31E-12
1.148	1.38E-12
1.152	1.46E-12
1.156	1.53E-12
1.16	1.61E-12
1.164	1.70E-12
1.168	1.78E-12
1.172	1.87E-12
1.176	1.96E-12
1.18	2.06E-12
1.184	2.15E-12
1.188	2.25E-12
1.192	2.35E-12
1.196	2.45E-12
1.2	2.56E-12
1.204	2.67E-12
1.208	2.78E-12
1.212	2.89E-12
1.216	3.01E-12
1.22	3.13E-12
1.224	3.25E-12
1.228	3.38E-12
1.232	3.51E-12
1.236	3.64E-12
1.24	3.77E-12
1.244	3.91E-12
1.248	4.05E-12
1.252	4.19E-12
1.256	4.33E-12
1.26	4.47E-12
1.264	4.61E-12
1.268	4.75E-12
1.272	4.89E-12
1.276	5.03E-12
1.28	5.16E-12
1.284	5.29E-12
1.288	5.42E-12
1.292	5.54E-12
1.296	5.66E-12
1.3	5.78E-12
1.304	5.90E-12
1.308	6.01E-12
1.312	6.12E-12
1.316	6.22E-12
1.32	6.33E-12
1.324	6.43E-12
1.328	6.52E-12
1.332	6.62E-12
1.336	6.71E-12
1.34	6.79E-12
1.344	6.88E-12
1.348	6.96E-12
1.352	7.03E-12
1.356	7.11E-12
1.36	7.17E-12
1.364	7.24E-12
1.368	7.30E-12
1.372	7.36E-12
1.376	7.42E-12
1.38	7.47E-12
1.384	7.52E-12
1.388	7.57E-12
1.392	7.62E-12
1.396	7.66E-12
1.4	7.70E-12
1.404	7.74E-12
1.408	7.78E-12
1.412	7.81E-12
1.416	7.85E-12
1.42	7.88E-12
1.424	7.91E-12
1.428	7.94E-12
1.432	7.97E-12
1.436	7.99E-12
1.44	8.02E-12
1.444	8.04E-12
1.448	8.07E-12
1.452	8.09E-12
1.456	8.11E-12
1.46	8.13E-12
1.464	8.15E-12
1.468	8.16E-12
1.472	8.18E-12
1.476	8.19E-12
1.48	8.21E-12
1.484	8.22E-12
1.488	8.23E-12
1.492	8.25E-12
1.496	8.26E-12
1.5	8.27E-12
1.504	8.28E-12
1.508	8.29E-12
1.512	8.29E-12
1.516	8.30E-12
1.52	8.31E-12
1.524	8.32E-12
1.528	8.32E-12
1.532	8.33E-12
1.536	8.34E-12
1.54	8.34E-12
1.544	8.35E-12
1.548	8.35E-12
1.552	8.36E-12
1.556	8.37E-12
1.56	8.37E-12
1.564	8.37E-12
1.568	8.38E-12
1.572	8.38E-12
1.576	8.39E-12
1.58	8.39E-12
1.584	8.39E-12
1.588	8.40E-12
1.592	8.40E-12
1.596	8.40E-12
1.6	8.41E-12
1.604	8.41E-12
1.608	8.41E-12
1.612	8.42E-12
1.616	8.42E-12
1.62	8.42E-12
1.624	8.42E-12
1.628	8.42E-12
1.632	8.43E-12
1.636	8.43E-12
1.64	8.43E-12
1.644	8.43E-12
1.648	8.43E-12
1.652	8.44E-12
1.656	8.44E-12
1.66	8.44E-12
1.664	8.44E-12
1.668	8.44E-12
1.672	8.44E-12
1.676	8.44E-12
1.68	8.45E-12
1.684	8.45E-12
1.688	8.45E-12
1.692	8.45E-12
1.696	8.45E-12
1.7	8.45E-12
1.704	8.45E-12
1.708	8.45E-12
1.712	8.45E-12
1.716	8.45E-12
1.72	8.45E-12
1.724	8.45E-12
1.728	8.45E-12
1.732	8.46E-12
1.736	8.46E-12
1.74	8.46E-12
1.744	8.46E-12
1.748	8.46E-12
1.752	8.46E-12
1.756	8.46E-12
1.76	8.46E-12
1.764	8.46E-12
1.768	8.46E-12
1.772	8.46E-12
1.776	8.46E-12
1.78	8.46E-12
1.784	8.46E-12
1.788	8.46E-12
1.792	8.46E-12
1.796	8.46E-12
1.8	8.46E-12
1.804	8.46E-12
1.808	8.46E-12
1.812	8.46E-12
1.816	8.46E-12
1.82	8.46E-12
1.824	8.46E-12
1.828	8.46E-12
1.832	8.46E-12
1.836	8.45E-12
1.84	8.45E-12
1.844	8.45E-12
1.848	8.45E-12
1.852	8.45E-12
1.856	8.45E-12
1.86	8.45E-12
1.864	8.45E-12
1.868	8.45E-12
1.872	8.45E-12
1.876	8.45E-12
1.88	8.45E-12
1.884	8.45E-12
1.888	8.45E-12
1.892	8.45E-12
1.896	8.45E-12
1.9	8.45E-12
1.904	8.45E-12
1.908	8.45E-12
1.912	8.45E-12
1.916	8.45E-12
1.92	8.44E-12
1.924	8.44E-12
1.928	8.44E-12
1.932	8.44E-12
1.936	8.44E-12
1.94	8.44E-12
1.944	8.44E-12
1.948	8.44E-12
1.952	8.44E-12
1.956	8.44E-12
1.96	8.44E-12
1.964	8.44E-12
1.968	8.44E-12
1.972	8.44E-12
1.976	8.44E-12
1.98	8.44E-12
1.984	8.44E-12
1.988	8.44E-12
1.992	8.44E-12
1.996	8.43E-12
2	8.43E-12
2.004	8.43E-12
2.008	8.43E-12
2.012	8.43E-12
2.016	8.43E-12
2.02	8.43E-12
2.024	8.43E-12
2.028	8.43E-12
2.032	8.43E-12
2.036	8.43E-12
2.04	8.43E-12
2.044	8.43E-12
2.048	8.43E-12
2.052	8.43E-12
2.056	8.43E-12
2.06	8.43E-12
2.064	8.43E-12
2.068	8.43E-12
2.072	8.43E-12
2.076	8.43E-12
2.08	8.43E-12
2.084	8.43E-12
2.088	8.43E-12
2.092	8.43E-12
2.096	8.43E-12
2.1	8.43E-12
2.104	8.43E-12
2.108	8.43E-12
2.112	8.43E-12
2.116	8.43E-12
2.12	8.43E-12
2.124	8.43E-12
2.128	8.43E-12
2.132	8.43E-12
2.136	8.43E-12
2.14	8.43E-12
2.144	8.43E-12
2.148	8.43E-12
2.152	8.43E-12
2.156	8.43E-12
2.16	8.43E-12
2.164	8.43E-12
2.168	8.43E-12
2.172	8.43E-12
2.176	8.43E-12
2.18	8.43E-12
2.184	8.43E-12
2.188	8.43E-12
2.192	8.43E-12
2.196	8.43E-12
2.2	8.43E-12
2.204	8.43E-12
2.208	8.43E-12
2.212	8.43E-12
2.216	8.43E-12
2.22	8.43E-12
2.224	8.43E-12
2.228	8.43E-12
2.232	8.43E-12
2.236	8.43E-12
2.24	8.43E-12
2.244	8.43E-12
2.248	8.43E-12
2.252	8.43E-12
2.256	8.44E-12
2.26	8.44E-12
2.264	8.44E-12
2.268	8.44E-12
2.272	8.44E-12
2.276	8.44E-12
2.28	8.44E-12
2.284	8.44E-12
2.288	8.44E-12
2.292	8.44E-12
2.296	8.44E-12
2.3	8.44E-12
2.304	8.44E-12
2.308	8.44E-12
2.312	8.44E-12
2.316	8.44E-12
2.32	8.44E-12
2.324	8.44E-12
2.328	8.44E-12
2.332	8.44E-12
2.336	8.44E-12
2.34	8.44E-12
2.344	8.44E-12
2.348	8.44E-12
2.352	8.44E-12
2.356	8.44E-12
2.36	8.44E-12
2.364	8.44E-12
2.368	8.44E-12
2.372	8.44E-12
2.376	8.44E-12
2.38	8.44E-12
2.384	8.44E-12
2.388	8.44E-12
2.392	8.44E-12
2.396	8.44E-12
2.4	8.44E-12
2.404	8.44E-12
2.408	8.44E-12
2.412	8.44E-12
2.416	8.44E-12
2.42	8.44E-12
2.424	8.44E-12
2.428	8.44E-12
2.432	8.44E-12
2.436	8.44E-12
2.44	8.44E-12
2.444	8.44E-12
2.448	8.44E-12
2.452	8.44E-12
2.456	8.44E-12
2.46	8.44E-12
2.464	8.44E-12
2.468	8.44E-12
2.472	8.44E-12
2.476	8.44E-12
2.48	8.44E-12
2.484	8.44E-12
2.488	8.44E-12
2.492	8.44E-12
2.496	8.44E-12
2.5	8.43E-12
2.504	8.43E-12
2.508	8.43E-12
2.512	8.43E-12
2.516	8.43E-12
2.52	8.43E-12
2.524	8.43E-12
2.528	8.43E-12
2.532	8.43E-12
2.536	8.43E-12
2.54	8.43E-12
2.544	8.43E-12
2.548	8.43E-12
2.552	8.43E-12
2.556	8.42E-12
2.56	8.42E-12
2.564	8.42E-12
2.568	8.42E-12
2.572	8.42E-12
2.576	8.42E-12
2.58	8.42E-12
2.584	8.42E-12
2.588	8.42E-12
2.592	8.42E-12
2.596	8.42E-12
2.6	8.42E-12
2.604	8.41E-12
2.608	8.41E-12
2.612	8.41E-12
2.616	8.41E-12
2.62	8.41E-12
2.624	8.41E-12
2.628	8.41E-12
2.632	8.41E-12
2.636	8.41E-12
2.64	8.41E-12
2.644	8.41E-12
2.648	8.40E-12
2.652	8.40E-12
2.656	8.40E-12
2.66	8.40E-12
2.664	8.40E-12
2.668	8.40E-12
2.672	8.40E-12
2.676	8.40E-12
2.68	8.40E-12
2.684	8.40E-12
2.688	8.40E-12
2.692	8.40E-12
2.696	8.40E-12
2.7	8.40E-12
2.704	8.40E-12
2.708	8.40E-12
2.712	8.40E-12
2.716	8.40E-12
2.72	8.40E-12
2.724	8.40E-12
2.728	8.40E-12
2.732	8.40E-12
2.736	8.40E-12
2.74	8.40E-12
2.744	8.40E-12
2.748	8.40E-12
2.752	8.40E-12
2.756	8.41E-12
2.76	8.41E-12
2.764	8.41E-12
2.768	8.41E-12
2.772	8.41E-12
2.776	8.41E-12
2.78	8.41E-12
2.784	8.41E-12
2.788	8.41E-12
2.792	8.41E-12
2.796	8.41E-12
2.8	8.42E-12
2.804	8.42E-12
2.808	8.42E-12
2.812	8.42E-12
2.816	8.42E-12
2.82	8.42E-12
2.824	8.42E-12
2.828	8.42E-12
2.832	8.42E-12
2.836	8.42E-12
2.84	8.43E-12
2.844	8.43E-12
2.848	8.43E-12
2.852	8.43E-12
2.856	8.43E-12
2.86	8.43E-12
2.864	8.43E-12
2.868	8.43E-12
2.872	8.43E-12
2.876	8.43E-12
2.88	8.43E-12
2.884	8.44E-12
2.888	8.44E-12
2.892	8.44E-12
2.896	8.44E-12
2.9	8.44E-12
2.904	8.44E-12
2.908	8.44E-12
2.912	8.44E-12
2.916	8.44E-12
2.92	8.44E-12
2.924	8.44E-12
2.928	8.44E-12
2.932	8.44E-12
2.936	8.45E-12
2.94	8.45E-12
2.944	8.45E-12
2.948	8.45E-12
2.952	8.45E-12
2.956	8.45E-12
2.96	8.45E-12
2.964	8.45E-12
2.968	8.45E-12
2.972	8.45E-12
2.976	8.45E-12
2.98	8.45E-12
2.984	8.45E-12
2.988	8.45E-12
2.992	8.45E-12
2.996	8.45E-12
3	8.45E-12
3.004	8.46E-12
3.008	8.46E-12
3.012	8.46E-12
3.016	8.46E-12
3.02	8.46E-12
3.024	8.46E-12
3.028	8.46E-12
3.032	8.46E-12
3.036	8.46E-12
3.04	8.46E-12
3.044	8.46E-12
3.048	8.46E-12
3.052	8.46E-12
3.056	8.46E-12
3.06	8.46E-12
3.064	8.46E-12
3.068	8.46E-12
3.072	8.46E-12
3.076	8.46E-12
3.08	8.46E-12
3.084	8.46E-12
3.088	8.47E-12
3.092	8.47E-12
3.096	8.47E-12
3.1	8.47E-12
3.104	8.47E-12
3.108	8.47E-12
3.112	8.47E-12
3.116	8.47E-12
3.12	8.47E-12
3.124	8.47E-12
3.128	8.47E-12
3.132	8.47E-12
3.136	8.47E-12
3.14	8.47E-12
3.144	8.47E-12
3.148	8.47E-12
3.152	8.47E-12
3.156	8.48E-12
3.16	8.48E-12
3.164	8.48E-12
3.168	8.48E-12
3.172	8.48E-12
3.176	8.48E-12
3.18	8.48E-12
3.184	8.48E-12
3.188	8.48E-12
3.192	8.48E-12
3.196	8.48E-12
3.2	8.48E-12
3.204	8.48E-12
3.208	8.48E-12
3.212	8.48E-12
3.216	8.48E-12
3.22	8.48E-12
3.224	8.48E-12
3.228	8.48E-12
3.232	8.48E-12
3.236	8.48E-12
3.24	8.48E-12
3.244	8.48E-12
3.248	8.48E-12
3.252	8.48E-12
3.256	8.49E-12
3.26	8.49E-12
3.264	8.49E-12
3.268	8.49E-12
3.272	8.49E-12
3.276	8.49E-12
3.28	8.49E-12
3.284	8.49E-12
3.288	8.49E-12
3.292	8.49E-12
3.296	8.49E-12
3.3	8.49E-12
3.304	8.49E-12
3.308	8.49E-12
3.312	8.49E-12
3.316	8.49E-12
3.32	8.49E-12
3.324	8.49E-12
3.328	8.49E-12
3.332	8.49E-12
3.336	8.49E-12
3.34	8.50E-12
3.344	8.50E-12
3.348	8.50E-12
3.352	8.50E-12
3.356	8.50E-12
3.36	8.50E-12
3.364	8.50E-12
3.368	8.50E-12
3.372	8.50E-12
3.376	8.51E-12
3.38	8.51E-12
3.384	8.51E-12
3.388	8.51E-12
3.392	8.51E-12
3.396	8.51E-12
3.4	8.52E-12
3.404	8.52E-12
3.408	8.52E-12
3.412	8.52E-12
3.416	8.52E-12
3.42	8.53E-12
3.424	8.53E-12
3.428	8.53E-12
3.432	8.54E-12
3.436	8.54E-12
3.44	8.54E-12
3.444	8.55E-12
3.448	8.55E-12
3.452	8.56E-12
3.456	8.56E-12
3.46	8.57E-12
3.464	8.57E-12
3.468	8.58E-12
3.472	8.58E-12
3.476	8.59E-12
3.48	8.60E-12
3.484	8.60E-12
3.488	8.61E-12
3.492	8.62E-12
3.496	8.63E-12
3.5	8.63E-12
3.504	8.64E-12
3.508	8.65E-12
3.512	8.66E-12
3.516	8.67E-12
3.52	8.68E-12
3.524	8.70E-12
3.528	8.71E-12
3.532	8.73E-12
3.536	8.74E-12
3.54	8.76E-12
3.544	8.78E-12
3.548	8.80E-12
3.552	8.82E-12
3.556	8.84E-12
3.56	8.87E-12
3.564	8.89E-12
3.568	8.92E-12
3.572	8.95E-12
3.576	8.97E-12
3.58	9.00E-12
3.584	9.04E-12
3.588	9.07E-12
3.592	9.10E-12
3.596	9.14E-12
3.6	9.17E-12
3.604	9.21E-12
3.608	9.25E-12
3.612	9.29E-12
3.616	9.34E-12
3.62	9.39E-12
3.624	9.44E-12
3.628	9.49E-12
3.632	9.55E-12
3.636	9.60E-12
3.64	9.67E-12
3.644	9.73E-12
3.648	9.80E-12
3.652	9.88E-12
3.656	9.95E-12
3.66	1.00E-11
3.664	1.01E-11
3.668	1.02E-11
3.672	1.03E-11
3.676	1.04E-11
3.68	1.05E-11
3.684	1.06E-11
3.688	1.07E-11
3.692	1.08E-11
3.696	1.09E-11
3.7	1.10E-11
3.704	1.11E-11
3.708	1.12E-11
3.712	1.13E-11
3.716	1.14E-11
3.72	1.15E-11
3.724	1.17E-11
3.728	1.18E-11
3.732	1.19E-11
3.736	1.21E-11
3.74	1.22E-11
3.744	1.23E-11
3.748	1.25E-11
3.752	1.26E-11
3.756	1.27E-11
3.76	1.29E-11
3.764	1.30E-11
3.768	1.32E-11
3.772	1.33E-11
3.776	1.34E-11
3.78	1.36E-11
3.784	1.37E-11
3.788	1.38E-11
3.792	1.40E-11
3.796	1.41E-11
3.8	1.42E-11
3.804	1.43E-11
3.808	1.44E-11
3.812	1.45E-11
3.816	1.46E-11
3.82	1.47E-11
3.824	1.48E-11
3.828	1.49E-11
3.832	1.50E-11
3.836	1.51E-11
3.84	1.52E-11
3.844	1.53E-11
3.848	1.54E-11
3.852	1.55E-11
3.856	1.55E-11
3.86	1.56E-11
3.864	1.57E-11
3.868	1.57E-11
3.872	1.58E-11
3.876	1.58E-11
3.88	1.59E-11
3.884	1.60E-11
3.888	1.60E-11
3.892	1.60E-11
3.896	1.61E-11
3.9	1.61E-11
3.904	1.62E-11
3.908	1.62E-11
3.912	1.62E-11
3.916	1.63E-11
3.92	1.63E-11
3.924	1.63E-11
3.928	1.64E-11
3.932	1.64E-11
3.936	1.64E-11
3.94	1.64E-11
3.944	1.65E-11
3.948	1.65E-11
3.952	1.65E-11
3.956	1.65E-11
3.96	1.66E-11
3.964	1.66E-11
3.968	1.66E-11
3.972	1.66E-11
3.976	1.66E-11
3.98	1.66E-11
3.984	1.67E-11
3.988	1.67E-11
3.992	1.67E-11
3.996	1.67E-11
4	1.67E-11
4.004	1.67E-11
4.008	1.67E-11
4.012	1.67E-11
4.016	1.67E-11
4.02	1.67E-11
4.024	1.67E-11
4.028	1.68E-11
4.032	1.68E-11
4.036	1.68E-11
4.04	1.68E-11
4.044	1.68E-11
4.048	1.68E-11
4.052	1.68E-11
4.056	1.68E-11
4.06	1.68E-11
4.064	1.68E-11
4.068	1.68E-11
4.072	1.68E-11
4.076	1.68E-11
4.08	1.68E-11
4.084	1.68E-11
4.088	1.68E-11
4.092	1.68E-11
4.096	1.68E-11
4.1	1.68E-11
4.104	1.68E-11
4.108	1.68E-11
4.112	1.68E-11
4.116	1.69E-11
4.12	1.69E-11
4.124	1.69E-11
4.128	1.69E-11
4.132	1.69E-11
4.136	1.69E-11
4.14	1.69E-11
4.144	1.69E-11
4.148	1.69E-11
4.152	1.69E-11
4.156	1.69E-11
4.16	1.69E-11
4.164	1.69E-11
4.168	1.69E-11
4.172	1.69E-11
4.176	1.69E-11
4.18	1.69E-11
4.184	1.69E-11
4.188	1.69E-11
4.192	1.69E-11
4.196	1.69E-11
4.2	1.69E-11
4.204	1.69E-11
4.208	1.69E-11
4.212	1.69E-11
4.216	1.69E-11
4.22	1.69E-11
4.224	1.69E-11
4.228	1.69E-11
4.232	1.69E-11
4.236	1.69E-11
4.24	1.69E-11
4.244	1.69E-11
4.248	1.69E-11
4.252	1.69E-11
4.256	1.69E-11
4.26	1.69E-11
4.264	1.69E-11
4.268	1.69E-11
4.272	1.69E-11
4.276	1.69E-11
4.28	1.69E-11
4.284	1.69E-11
4.288	1.69E-11
4.292	1.69E-11
4.296	1.69E-11
4.3	1.69E-11
4.304	1.69E-11
4.308	1.69E-11
4.312	1.69E-11
4.316	1.69E-11
4.32	1.69E-11
4.324	1.69E-11
4.328	1.69E-11
4.332	1.69E-11
4.336	1.69E-11
4.34	1.69E-11
4.344	1.69E-11
4.348	1.69E-11
4.352	1.69E-11
4.356	1.69E-11
4.36	1.69E-11
4.364	1.69E-11
4.368	1.69E-11
4.372	1.69E-11
4.376	1.69E-11
4.38	1.69E-11
4.384	1.69E-11
4.388	1.69E-11
4.392	1.69E-11
4.396	1.69E-11
4.4	1.69E-11
4.404	1.69E-11
4.408	1.69E-11
4.412	1.69E-11
4.416	1.69E-11
4.42	1.69E-11
4.424	1.69E-11
4.428	1.69E-11
4.432	1.69E-11
4.436	1.69E-11
4.44	1.69E-11
4.444	1.69E-11
4.448	1.69E-11
4.452	1.69E-11
4.456	1.69E-11
4.46	1.69E-11
4.464	1.69E-11
4.468	1.69E-11
4.472	1.69E-11
4.476	1.69E-11
4.48	1.69E-11
4.484	1.69E-11
4.488	1.69E-11
4.492	1.69E-11
4.496	1.69E-11
4.5	1.69E-11
4.504	1.69E-11
4.508	1.69E-11
4.512	1.69E-11
4.516	1.69E-11
4.52	1.69E-11
4.524	1.69E-11
4.528	1.69E-11
4.532	1.69E-11
4.536	1.69E-11
4.54	1.69E-11
4.544	1.69E-11
4.548	1.69E-11
4.552	1.69E-11
4.556	1.69E-11
4.56	1.69E-11
4.564	1.69E-11
4.568	1.69E-11
4.572	1.69E-11
4.576	1.69E-11
4.58	1.69E-11
4.584	1.69E-11
4.588	1.69E-11
4.592	1.69E-11
4.596	1.69E-11
4.6	1.69E-11
4.604	1.69E-11
4.608	1.69E-11
4.612	1.69E-11
4.616	1.69E-11
4.62	1.69E-11
4.624	1.69E-11
4.628	1.69E-11
4.632	1.69E-11
4.636	1.69E-11
4.64	1.69E-11
4.644	1.69E-11
4.648	1.69E-11
4.652	1.69E-11
4.656	1.69E-11
4.66	1.69E-11
4.664	1.69E-11
4.668	1.69E-11
4.672	1.69E-11
4.676	1.69E-11
4.68	1.69E-11
4.684	1.69E-11
4.688	1.69E-11
4.692	1.69E-11
4.696	1.69E-11
4.7	1.69E-11
4.704	1.69E-11
4.708	1.69E-11
4.712	1.69E-11
4.716	1.69E-11
4.72	1.69E-11
4.724	1.69E-11
4.728	1.69E-11
4.732	1.69E-11
4.736	1.69E-11
4.74	1.69E-11
4.744	1.69E-11
4.748	1.69E-11
4.752	1.69E-11
4.756	1.69E-11
4.76	1.69E-11
4.764	1.69E-11
4.768	1.69E-11
4.772	1.69E-11
4.776	1.69E-11
4.78	1.69E-11
4.784	1.69E-11
4.788	1.69E-11
4.792	1.69E-11
4.796	1.69E-11
4.8	1.69E-11
4.804	1.69E-11
4.808	1.69E-11
4.812	1.69E-11
4.816	1.69E-11
4.82	1.69E-11
4.824	1.69E-11
4.828	1.69E-11
4.832	1.69E-11
4.836	1.69E-11
4.84	1.69E-11
4.844	1.69E-11
4.848	1.69E-11
4.852	1.69E-11
4.856	1.69E-11
4.86	1.69E-11
4.864	1.69E-11
4.868	1.69E-11
4.872	1.69E-11
4.876	1.69E-11
4.88	1.69E-11
4.884	1.69E-11
4.888	1.69E-11
4.892	1.69E-11
4.896	1.69E-11
4.9	1.69E-11
4.904	1.69E-11
4.908	1.69E-11
4.912	1.69E-11
4.916	1.69E-11
4.92	1.69E-11
4.924	1.69E-11
4.928	1.69E-11
4.932	1.69E-11
4.936	1.69E-11
4.94	1.69E-11
4.944	1.69E-11
4.948	1.69E-11
4.952	1.69E-11
4.956	1.69E-11
4.96	1.69E-11
4.964	1.69E-11
4.968	1.69E-11
4.972	1.69E-11
4.976	1.69E-11
4.98	1.69E-11
4.984	1.69E-11
4.988	1.69E-11
4.992	1.69E-11
4.996	1.69E-11
5	1.69E-11
5.004	1.69E-11
5.008	1.69E-11
5.012	1.69E-11
5.016	1.69E-11
5.02	1.69E-11
5.024	1.69E-11
5.028	1.69E-11
5.032	1.69E-11
5.036	1.69E-11
5.04	1.69E-11
5.044	1.69E-11
5.048	1.69E-11
5.052	1.69E-11
5.056	1.69E-11
5.06	1.69E-11
5.064	1.69E-11
5.068	1.69E-11
5.072	1.69E-11
5.076	1.69E-11
5.08	1.69E-11
5.084	1.69E-11
5.088	1.69E-11
5.092	1.69E-11
5.096	1.68E-11
5.1	1.68E-11
5.104	1.68E-11
5.108	1.68E-11
5.112	1.68E-11
5.116	1.68E-11
5.12	1.68E-11
5.124	1.68E-11
5.128	1.68E-11
5.132	1.68E-11
5.136	1.68E-11
5.14	1.68E-11
5.144	1.68E-11
5.148	1.68E-11
5.152	1.68E-11
5.156	1.68E-11
5.16	1.68E-11
5.164	1.68E-11
5.168	1.68E-11
5.172	1.68E-11
5.176	1.68E-11
5.18	1.68E-11
5.184	1.68E-11
5.188	1.68E-11
5.192	1.68E-11
5.196	1.68E-11
5.2	1.68E-11
5.204	1.68E-11
5.208	1.68E-11
5.212	1.68E-11
5.216	1.68E-11
5.22	1.68E-11
5.224	1.68E-11
5.228	1.68E-11
5.232	1.68E-11
5.236	1.68E-11
5.24	1.68E-11
5.244	1.68E-11
5.248	1.68E-11
5.252	1.68E-11
5.256	1.68E-11
5.26	1.68E-11
5.264	1.68E-11
5.268	1.68E-11
5.272	1.68E-11
5.276	1.68E-11
5.28	1.68E-11
5.284	1.68E-11
5.288	1.68E-11
5.292	1.68E-11
5.296	1.68E-11
5.3	1.68E-11
5.304	1.68E-11
5.308	1.69E-11
5.312	1.69E-11
5.316	1.69E-11
5.32	1.69E-11
5.324	1.69E-11
5.328	1.69E-11
5.332	1.69E-11
5.336	1.69E-11
5.34	1.69E-11
5.344	1.69E-11
5.348	1.69E-11
5.352	1.69E-11
5.356	1.69E-11
5.36	1.69E-11
5.364	1.69E-11
5.368	1.69E-11
5.372	1.69E-11
5.376	1.69E-11
5.38	1.69E-11
5.384	1.69E-11
5.388	1.69E-11
5.392	1.69E-11
5.396	1.69E-11
5.4	1.69E-11
5.404	1.69E-11
5.408	1.69E-11
5.412	1.69E-11
5.416	1.69E-11
5.42	1.69E-11
5.424	1.69E-11
5.428	1.69E-11
5.432	1.69E-11
5.436	1.69E-11
5.44	1.69E-11
5.444	1.69E-11
5.448	1.69E-11
5.452	1.69E-11
5.456	1.69E-11
5.46	1.69E-11
5.464	1.69E-11
5.468	1.69E-11
5.472	1.69E-11
5.476	1.69E-11
5.48	1.69E-11
5.484	1.69E-11
5.488	1.69E-11
5.492	1.69E-11
5.496	1.69E-11
5.5	1.69E-11
5.504	1.69E-11
5.508	1.69E-11
5.512	1.69E-11
5.516	1.69E-11
5.52	1.69E-11
5.524	1.69E-11
5.528	1.69E-11
5.532	1.69E-11
5.536	1.69E-11
5.54	1.69E-11
5.544	1.69E-11
5.548	1.69E-11
5.552	1.69E-11
5.556	1.69E-11
5.56	1.69E-11
5.564	1.69E-11
5.568	1.69E-11
5.572	1.69E-11
5.576	1.69E-11
5.58	1.69E-11
5.584	1.69E-11
5.588	1.69E-11
5.592	1.69E-11
5.596	1.69E-11
5.6	1.69E-11
5.604	1.69E-11
5.608	1.69E-11
5.612	1.69E-11
5.616	1.69E-11
5.62	1.69E-11
5.624	1.69E-11
5.628	1.69E-11
5.632	1.69E-11
5.636	1.69E-11
5.64	1.69E-11
5.644	1.69E-11
5.648	1.69E-11
5.652	1.69E-11
5.656	1.69E-11
5.66	1.69E-11
5.664	1.69E-11
5.668	1.69E-11
5.672	1.69E-11
5.676	1.69E-11
5.68	1.69E-11
5.684	1.69E-11
5.688	1.69E-11
5.692	1.69E-11
5.696	1.69E-11
5.7	1.69E-11
5.704	1.69E-11
5.708	1.69E-11
5.712	1.69E-11
5.716	1.69E-11
5.72	1.69E-11
5.724	1.69E-11
5.728	1.69E-11
5.732	1.69E-11
5.736	1.69E-11
5.74	1.69E-11
5.744	1.69E-11
5.748	1.69E-11
5.752	1.69E-11
5.756	1.69E-11
5.76	1.69E-11
5.764	1.69E-11
5.768	1.69E-11
5.772	1.69E-11
5.776	1.69E-11
5.78	1.69E-11
5.784	1.69E-11
5.788	1.69E-11
5.792	1.69E-11
5.796	1.69E-11
5.8	1.69E-11
5.804	1.69E-11
5.808	1.69E-11
5.812	1.69E-11
5.816	1.69E-11
5.82	1.69E-11
5.824	1.69E-11
5.828	1.69E-11
5.832	1.69E-11
5.836	1.69E-11
5.84	1.69E-11
5.844	1.69E-11
5.848	1.69E-11
5.852	1.69E-11
5.856	1.69E-11
5.86	1.69E-11
5.864	1.69E-11
5.868	1.69E-11
5.872	1.69E-11
5.876	1.69E-11
5.88	1.69E-11
5.884	1.69E-11
5.888	1.69E-11
5.892	1.69E-11
5.896	1.69E-11
5.9	1.69E-11
5.904	1.70E-11
5.908	1.70E-11
5.912	1.70E-11
5.916	1.70E-11
5.92	1.70E-11
5.924	1.70E-11
5.928	1.70E-11
5.932	1.70E-11
5.936	1.70E-11
5.94	1.70E-11
5.944	1.70E-11
5.948	1.70E-11
5.952	1.70E-11
5.956	1.70E-11
5.96	1.70E-11
5.964	1.70E-11
5.968	1.70E-11
5.972	1.70E-11
5.976	1.70E-11
5.98	1.70E-11
5.984	1.70E-11
5.988	1.70E-11
5.992	1.71E-11
5.996	1.71E-11
6	1.71E-11
6.004	1.71E-11
6.008	1.71E-11
6.012	1.71E-11
6.016	1.71E-11
6.02	1.71E-11
6.024	1.71E-11
6.028	1.71E-11
6.032	1.72E-11
6.036	1.72E-11
6.04	1.72E-11
6.044	1.72E-11
6.048	1.72E-11
6.052	1.72E-11
6.056	1.73E-11
6.06	1.73E-11
6.064	1.73E-11
6.068	1.73E-11
6.072	1.74E-11
6.076	1.74E-11
6.08	1.74E-11
6.084	1.75E-11
6.088	1.75E-11
6.092	1.75E-11
6.096	1.76E-11
6.1	1.76E-11
6.104	1.76E-11
6.108	1.77E-11
6.112	1.77E-11
6.116	1.78E-11
6.12	1.78E-11
6.124	1.79E-11
6.128	1.79E-11
6.132	1.80E-11
6.136	1.80E-11
6.14	1.81E-11
6.144	1.82E-11
6.148	1.82E-11
6.152	1.83E-11
6.156	1.84E-11
6.16	1.85E-11
6.164	1.85E-11
6.168	1.86E-11
6.172	1.87E-11
6.176	1.88E-11
6.18	1.89E-11
6.184	1.90E-11
6.188	1.91E-11
6.192	1.92E-11
6.196	1.93E-11
6.2	1.94E-11
6.204	1.95E-11
6.208	1.96E-11
6.212	1.97E-11
6.216	1.98E-11
6.22	2.00E-11
6.224	2.01E-11
6.228	2.02E-11
6.232	2.03E-11
6.236	2.05E-11
6.24	2.06E-11
6.244	2.07E-11
6.248	2.09E-11
6.252	2.10E-11
6.256	2.12E-11
6.26	2.13E-11
6.264	2.14E-11
6.268	2.16E-11
6.272	2.17E-11
6.276	2.18E-11
6.28	2.20E-11
6.284	2.21E-11
6.288	2.22E-11
6.292	2.24E-11
6.296	2.25E-11
6.3	2.26E-11
6.304	2.27E-11
6.308	2.28E-11
6.312	2.30E-11
6.316	2.31E-11
6.32	2.32E-11
6.324	2.33E-11
6.328	2.34E-11
6.332	2.35E-11
6.336	2.35E-11
6.34	2.36E-11
6.344	2.37E-11
6.348	2.38E-11
6.352	2.39E-11
6.356	2.39E-11
6.36	2.40E-11
6.364	2.41E-11
6.368	2.42E-11
6.372	2.42E-11
6.376	2.43E-11
6.38	2.43E-11
6.384	2.44E-11
6.388	2.44E-11
6.392	2.45E-11
6.396	2.45E-11
6.4	2.46E-11
6.404	2.46E-11
6.408	2.46E-11
6.412	2.47E-11
6.416	2.47E-11
6.42	2.47E-11
6.424	2.48E-11
6.428	2.48E-11
6.432	2.48E-11
6.436	2.49E-11
6.44	2.49E-11
6.444	2.49E-11
6.448	2.49E-11
6.452	2.49E-11
6.456	2.50E-11
6.46	2.50E-11
6.464	2.50E-11
6.468	2.50E-11
6.472	2.50E-11
6.476	2.51E-11
6.48	2.51E-11
6.484	2.51E-11
6.488	2.51E-11
6.492	2.51E-11
6.496	2.51E-11
6.5	2.51E-11
6.504	2.51E-11
6.508	2.51E-11
6.512	2.52E-11
6.516	2.52E-11
6.52	2.52E-11
6.524	2.52E-11
6.528	2.52E-11
6.532	2.52E-11
6.536	2.52E-11
6.54	2.52E-11
6.544	2.52E-11
6.548	2.52E-11
6.552	2.52E-11
6.556	2.52E-11
6.56	2.52E-11
6.564	2.52E-11
6.568	2.52E-11
6.572	2.52E-11
6.576	2.53E-11
6.58	2.53E-11
6.584	2.53E-11
6.588	2.53E-11
6.592	2.53E-11
6.596	2.53E-11
6.6	2.53E-11
6.604	2.53E-11
6.608	2.53E-11
6.612	2.53E-11
6.616	2.53E-11
6.62	2.53E-11
6.624	2.53E-11
6.628	2.53E-11
6.632	2.53E-11
6.636	2.53E-11
6.64	2.53E-11
6.644	2.53E-11
6.648	2.53E-11
6.652	2.53E-11
6.656	2.53E-11
6.66	2.53E-11
6.664	2.53E-11
6.668	2.53E-11
6.672	2.53E-11
6.676	2.53E-11
6.68	2.53E-11
6.684	2.53E-11
6.688	2.53E-11
6.692	2.53E-11
6.696	2.53E-11
6.7	2.53E-11
6.704	2.53E-11
6.708	2.53E-11
6.712	2.53E-11
6.716	2.53E-11
6.72	2.53E-11
6.724	2.53E-11
6.728	2.53E-11
6.732	2.53E-11
6.736	2.53E-11
6.74	2.53E-11
6.744	2.53E-11
6.748	2.53E-11
6.752	2.53E-11
6.756	2.53E-11
6.76	2.53E-11
6.764	2.53E-11
6.768	2.53E-11
6.772	2.53E-11
6.776	2.53E-11
6.78	2.53E-11
6.784	2.53E-11
6.788	2.53E-11
6.792	2.53E-11
6.796	2.53E-11
6.8	2.53E-11
6.804	2.53E-11
6.808	2.53E-11
6.812	2.53E-11
6.816	2.53E-11
6.82	2.53E-11
6.824	2.53E-11
6.828	2.53E-11
6.832	2.53E-11
6.836	2.53E-11
6.84	2.53E-11
6.844	2.53E-11
6.848	2.53E-11
6.852	2.53E-11
6.856	2.53E-11
6.86	2.53E-11
6.864	2.53E-11
6.868	2.53E-11
6.872	2.53E-11
6.876	2.53E-11
6.88	2.53E-11
6.884	2.53E-11
6.888	2.53E-11
6.892	2.53E-11
6.896	2.53E-11
6.9	2.53E-11
6.904	2.53E-11
6.908	2.53E-11
6.912	2.53E-11
6.916	2.53E-11
6.92	2.53E-11
6.924	2.53E-11
6.928	2.53E-11
6.932	2.53E-11
6.936	2.53E-11
6.94	2.53E-11
6.944	2.53E-11
6.948	2.53E-11
6.952	2.53E-11
6.956	2.53E-11
6.96	2.53E-11
6.964	2.53E-11
6.968	2.53E-11
6.972	2.53E-11
6.976	2.53E-11
6.98	2.53E-11
6.984	2.53E-11
6.988	2.53E-11
6.992	2.53E-11
6.996	2.53E-11
7	2.53E-11
7.004	2.53E-11
7.008	2.53E-11
7.012	2.53E-11
7.016	2.53E-11
7.02	2.53E-11
7.024	2.53E-11
7.028	2.53E-11
7.032	2.53E-11
7.036	2.53E-11
7.04	2.53E-11
7.044	2.53E-11
7.048	2.53E-11
7.052	2.53E-11
7.056	2.53E-11
7.06	2.53E-11
7.064	2.53E-11
7.068	2.53E-11
7.072	2.53E-11
7.076	2.53E-11
7.08	2.53E-11
7.084	2.53E-11
7.088	2.53E-11
7.092	2.53E-11
7.096	2.53E-11
7.1	2.53E-11
7.104	2.53E-11
7.108	2.53E-11
7.112	2.53E-11
7.116	2.53E-11
7.12	2.53E-11
7.124	2.53E-11
7.128	2.53E-11
7.132	2.53E-11
7.136	2.53E-11
7.14	2.53E-11
7.144	2.53E-11
7.148	2.53E-11
7.152	2.53E-11
7.156	2.53E-11
7.16	2.53E-11
7.164	2.53E-11
7.168	2.53E-11
7.172	2.53E-11
7.176	2.53E-11
7.18	2.53E-11
7.184	2.53E-11
7.188	2.53E-11
7.192	2.53E-11
7.196	2.53E-11
7.2	2.53E-11
7.204	2.53E-11
7.208	2.53E-11
7.212	2.53E-11
7.216	2.53E-11
7.22	2.53E-11
7.224	2.53E-11
7.228	2.53E-11
7.232	2.53E-11
7.236	2.53E-11
7.24	2.53E-11
7.244	2.53E-11
7.248	2.53E-11
7.252	2.53E-11
7.256	2.53E-11
7.26	2.53E-11
7.264	2.53E-11
7.268	2.53E-11
7.272	2.53E-11
7.276	2.53E-11
7.28	2.53E-11
7.284	2.53E-11
7.288	2.53E-11
7.292	2.53E-11
7.296	2.53E-11
7.3	2.53E-11
7.304	2.53E-11
7.308	2.53E-11
7.312	2.53E-11
7.316	2.53E-11
7.32	2.53E-11
7.324	2.53E-11
7.328	2.53E-11
7.332	2.53E-11
7.336	2.53E-11
7.34	2.53E-11
7.344	2.53E-11
7.348	2.53E-11
7.352	2.53E-11
7.356	2.53E-11
7.36	2.53E-11
7.364	2.53E-11
7.368	2.53E-11
7.372	2.53E-11
7.376	2.53E-11
7.38	2.53E-11
7.384	2.53E-11
7.388	2.53E-11
7.392	2.53E-11
7.396	2.53E-11
7.4	2.53E-11
7.404	2.53E-11
7.408	2.53E-11
7.412	2.53E-11
7.416	2.53E-11
7.42	2.53E-11
7.424	2.53E-11
7.428	2.53E-11
7.432	2.53E-11
7.436	2.53E-11
7.44	2.53E-11
7.444	2.53E-11
7.448	2.53E-11
7.452	2.53E-11
7.456	2.53E-11
7.46	2.53E-11
7.464	2.53E-11
7.468	2.53E-11
7.472	2.53E-11
7.476	2.53E-11
7.48	2.53E-11
7.484	2.53E-11
7.488	2.53E-11
7.492	2.53E-11
7.496	2.53E-11
7.5	2.53E-11
7.504	2.53E-11
7.508	2.53E-11
7.512	2.53E-11
7.516	2.53E-11
7.52	2.53E-11
7.524	2.53E-11
7.528	2.53E-11
7.532	2.53E-11
7.536	2.53E-11
7.54	2.53E-11
7.544	2.53E-11
7.548	2.53E-11
7.552	2.53E-11
7.556	2.53E-11
7.56	2.53E-11
7.564	2.53E-11
7.568	2.53E-11
7.572	2.53E-11
7.576	2.53E-11
7.58	2.53E-11
7.584	2.53E-11
7.588	2.53E-11
7.592	2.53E-11
7.596	2.53E-11
7.6	2.53E-11
7.604	2.53E-11
7.608	2.53E-11
7.612	2.53E-11
7.616	2.53E-11
7.62	2.53E-11
7.624	2.53E-11
7.628	2.53E-11
7.632	2.53E-11
7.636	2.53E-11
7.64	2.53E-11
7.644	2.53E-11
7.648	2.53E-11
7.652	2.53E-11
7.656	2.53E-11
7.66	2.53E-11
7.664	2.53E-11
7.668	2.53E-11
7.672	2.53E-11
7.676	2.53E-11
7.68	2.53E-11
7.684	2.53E-11
7.688	2.53E-11
7.692	2.53E-11
7.696	2.53E-11
7.7	2.53E-11
7.704	2.53E-11
7.708	2.53E-11
7.712	2.53E-11
7.716	2.53E-11
7.72	2.53E-11
7.724	2.53E-11
7.728	2.53E-11
7.732	2.53E-11
7.736	2.53E-11
7.74	2.53E-11
7.744	2.53E-11
7.748	2.53E-11
7.752	2.53E-11
7.756	2.53E-11
7.76	2.53E-11
7.764	2.53E-11
7.768	2.53E-11
7.772	2.53E-11
7.776	2.53E-11
7.78	2.53E-11
7.784	2.53E-11
7.788	2.53E-11
7.792	2.53E-11
7.796	2.53E-11
7.8	2.53E-11
7.804	2.53E-11
7.808	2.53E-11
7.812	2.53E-11
7.816	2.53E-11
7.82	2.53E-11
7.824	2.53E-11
7.828	2.53E-11
7.832	2.53E-11
7.836	2.53E-11
7.84	2.53E-11
7.844	2.53E-11
7.848	2.53E-11
7.852	2.53E-11
7.856	2.53E-11
7.86	2.53E-11
7.864	2.53E-11
7.868	2.53E-11
7.872	2.53E-11
7.876	2.53E-11
7.88	2.53E-11
7.884	2.53E-11
7.888	2.53E-11
7.892	2.53E-11
7.896	2.53E-11
7.9	2.53E-11
7.904	2.53E-11
7.908	2.53E-11
7.912	2.53E-11
7.916	2.53E-11
7.92	2.53E-11
7.924	2.53E-11
7.928	2.53E-11
7.932	2.53E-11
7.936	2.53E-11
7.94	2.53E-11
7.944	2.53E-11
7.948	2.53E-11
7.952	2.53E-11
7.956	2.53E-11
7.96	2.53E-11
7.964	2.53E-11
7.968	2.53E-11
7.972	2.53E-11
7.976	2.53E-11
7.98	2.53E-11
7.984	2.53E-11
7.988	2.53E-11
7.992	2.53E-11
7.996	2.53E-11
8	2.53E-11
8.004	2.53E-11
8.008	2.53E-11
8.012	2.53E-11
8.016	2.53E-11
8.02	2.53E-11
8.024	2.53E-11
8.028	2.53E-11
8.032	2.53E-11
8.036	2.53E-11
8.04	2.53E-11
8.044	2.53E-11
8.048	2.53E-11
8.052	2.53E-11
8.056	2.53E-11
8.06	2.53E-11
8.064	2.53E-11
8.068	2.53E-11
8.072	2.53E-11
8.076	2.53E-11
8.08	2.53E-11
8.084	2.53E-11
8.088	2.53E-11
8.092	2.53E-11
8.096	2.53E-11
8.1	2.53E-11
8.104	2.53E-11
8.108	2.53E-11
8.112	2.53E-11
8.116	2.53E-11
8.12	2.53E-11
8.124	2.53E-11
8.128	2.53E-11
8.132	2.53E-11
8.136	2.53E-11
8.14	2.53E-11
8.144	2.53E-11
8.148	2.53E-11
8.152	2.53E-11
8.156	2.53E-11
8.16	2.53E-11
8.164	2.53E-11
8.168	2.53E-11
8.172	2.53E-11
8.176	2.53E-11
8.18	2.53E-11
8.184	2.53E-11
8.188	2.53E-11
8.192	2.53E-11
8.196	2.53E-11
8.2	2.53E-11
8.204	2.53E-11
8.208	2.53E-11
8.212	2.53E-11
8.216	2.53E-11
8.22	2.53E-11
8.224	2.53E-11
8.228	2.53E-11
8.232	2.53E-11
8.236	2.53E-11
8.24	2.53E-11
8.244	2.53E-11
8.248	2.54E-11
8.252	2.54E-11
8.256	2.54E-11
8.26	2.54E-11
8.264	2.54E-11
8.268	2.54E-11
8.272	2.54E-11
8.276	2.54E-11
8.28	2.54E-11
8.284	2.54E-11
8.288	2.54E-11
8.292	2.54E-11
8.296	2.54E-11
8.3	2.54E-11
8.304	2.54E-11
8.308	2.54E-11
8.312	2.54E-11
8.316	2.54E-11
8.32	2.54E-11
8.324	2.54E-11
8.328	2.54E-11
8.332	2.54E-11
8.336	2.54E-11
8.34	2.54E-11
8.344	2.54E-11
8.348	2.54E-11
8.352	2.54E-11
8.356	2.54E-11
8.36	2.54E-11
8.364	2.54E-11
8.368	2.54E-11
8.372	2.54E-11
8.376	2.54E-11
8.38	2.54E-11
8.384	2.54E-11
8.388	2.54E-11
8.392	2.54E-11
8.396	2.54E-11
8.4	2.54E-11
8.404	2.54E-11
8.408	2.54E-11
8.412	2.54E-11
8.416	2.54E-11
8.42	2.54E-11
8.424	2.54E-11
8.428	2.54E-11
8.432	2.54E-11
8.436	2.54E-11
8.44	2.54E-11
8.444	2.54E-11
8.448	2.54E-11
8.452	2.54E-11
8.456	2.54E-11
8.46	2.54E-11
8.464	2.54E-11
8.468	2.54E-11
8.472	2.54E-11
8.476	2.55E-11
8.48	2.55E-11
8.484	2.55E-11
8.488	2.55E-11
8.492	2.55E-11
8.496	2.55E-11
8.5	2.55E-11
8.504	2.55E-11
8.508	2.55E-11
8.512	2.55E-11
8.516	2.55E-11
8.52	2.55E-11
8.524	2.56E-11
8.528	2.56E-11
8.532	2.56E-11
8.536	2.56E-11
8.54	2.56E-11
8.544	2.56E-11
8.548	2.57E-11
8.552	2.57E-11
8.556	2.57E-11
8.56	2.57E-11
8.564	2.57E-11
8.568	2.58E-11
8.572	2.58E-11
8.576	2.58E-11
8.58	2.59E-11
8.584	2.59E-11
8.588	2.59E-11
8.592	2.60E-11
8.596	2.60E-11
8.6	2.60E-11
8.604	2.61E-11
8.608	2.61E-11
8.612	2.61E-11
8.616	2.62E-11
8.62	2.62E-11
8.624	2.63E-11
8.628	2.63E-11
8.632	2.64E-11
8.636	2.65E-11
8.64	2.65E-11
8.644	2.66E-11
8.648	2.66E-11
8.652	2.67E-11
8.656	2.68E-11
8.66	2.69E-11
8.664	2.70E-11
8.668	2.70E-11
8.672	2.71E-11
8.676	2.72E-11
8.68	2.73E-11
8.684	2.74E-11
8.688	2.75E-11
8.692	2.76E-11
8.696	2.77E-11
8.7	2.78E-11
8.704	2.79E-11
8.708	2.80E-11
8.712	2.81E-11
8.716	2.83E-11
8.72	2.84E-11
8.724	2.85E-11
8.728	2.86E-11
8.732	2.87E-11
8.736	2.89E-11
8.74	2.90E-11
8.744	2.91E-11
8.748	2.93E-11
8.752	2.94E-11
8.756	2.96E-11
8.76	2.97E-11
8.764	2.98E-11
8.768	3.00E-11
8.772	3.01E-11
8.776	3.03E-11
8.78	3.04E-11
8.784	3.05E-11
8.788	3.07E-11
8.792	3.08E-11
8.796	3.09E-11
8.8	3.10E-11
8.804	3.11E-11
8.808	3.13E-11
8.812	3.14E-11
8.816	3.15E-11
8.82	3.16E-11
8.824	3.17E-11
8.828	3.18E-11
8.832	3.19E-11
8.836	3.20E-11
8.84	3.21E-11
8.844	3.21E-11
8.848	3.22E-11
8.852	3.23E-11
8.856	3.24E-11
8.86	3.24E-11
8.864	3.25E-11
8.868	3.26E-11
8.872	3.26E-11
8.876	3.27E-11
8.88	3.27E-11
8.884	3.28E-11
8.888	3.28E-11
8.892	3.29E-11
8.896	3.29E-11
8.9	3.30E-11
8.904	3.30E-11
8.908	3.31E-11
8.912	3.31E-11
8.916	3.31E-11
8.92	3.32E-11
8.924	3.32E-11
8.928	3.32E-11
8.932	3.33E-11
8.936	3.33E-11
8.94	3.33E-11
8.944	3.33E-11
8.948	3.34E-11
8.952	3.34E-11
8.956	3.34E-11
8.96	3.34E-11
8.964	3.34E-11
8.968	3.35E-11
8.972	3.35E-11
8.976	3.35E-11
8.98	3.35E-11
8.984	3.35E-11
8.988	3.35E-11
8.992	3.35E-11
8.996	3.36E-11
9	3.36E-11
9.004	3.36E-11
9.008	3.36E-11
9.012	3.36E-11
9.016	3.36E-11
9.02	3.36E-11
9.024	3.36E-11
9.028	3.36E-11
9.032	3.36E-11
9.036	3.36E-11
9.04	3.36E-11
9.044	3.36E-11
9.048	3.37E-11
9.052	3.37E-11
9.056	3.37E-11
9.06	3.37E-11
9.064	3.37E-11
9.068	3.37E-11
9.072	3.37E-11
9.076	3.37E-11
9.08	3.37E-11
9.084	3.37E-11
9.088	3.37E-11
9.092	3.37E-11
9.096	3.37E-11
9.1	3.37E-11
9.104	3.37E-11
9.108	3.37E-11
9.112	3.37E-11
9.116	3.37E-11
9.12	3.37E-11
9.124	3.37E-11
9.128	3.37E-11
9.132	3.37E-11
9.136	3.37E-11
9.14	3.37E-11
9.144	3.37E-11
9.148	3.37E-11
9.152	3.37E-11
9.156	3.37E-11
9.16	3.37E-11
9.164	3.37E-11
9.168	3.37E-11
9.172	3.37E-11
9.176	3.37E-11
9.18	3.37E-11
9.184	3.37E-11
9.188	3.37E-11
9.192	3.37E-11
9.196	3.37E-11
9.2	3.37E-11
9.204	3.38E-11
9.208	3.38E-11
9.212	3.38E-11
9.216	3.38E-11
9.22	3.38E-11
9.224	3.38E-11
9.228	3.38E-11
9.232	3.38E-11
9.236	3.38E-11
9.24	3.38E-11
9.244	3.38E-11
9.248	3.38E-11
9.252	3.38E-11
9.256	3.38E-11
9.26	3.38E-11
9.264	3.38E-11
9.268	3.38E-11
9.272	3.38E-11
9.276	3.38E-11
9.28	3.38E-11
9.284	3.38E-11
9.288	3.38E-11
9.292	3.38E-11
9.296	3.38E-11
9.3	3.38E-11
9.304	3.38E-11
9.308	3.38E-11
9.312	3.38E-11
9.316	3.38E-11
9.32	3.38E-11
9.324	3.38E-11
9.328	3.38E-11
9.332	3.38E-11
9.336	3.38E-11
9.34	3.38E-11
9.344	3.38E-11
9.348	3.38E-11
9.352	3.38E-11
9.356	3.38E-11
9.36	3.38E-11
9.364	3.38E-11
9.368	3.38E-11
9.372	3.38E-11
9.376	3.37E-11
9.38	3.37E-11
9.384	3.37E-11
9.388	3.37E-11
9.392	3.37E-11
9.396	3.37E-11
9.4	3.37E-11
9.404	3.37E-11
9.408	3.37E-11
9.412	3.37E-11
9.416	3.37E-11
9.42	3.37E-11
9.424	3.37E-11
9.428	3.37E-11
9.432	3.37E-11
9.436	3.37E-11
9.44	3.37E-11
9.444	3.37E-11
9.448	3.37E-11
9.452	3.37E-11
9.456	3.37E-11
9.46	3.37E-11
9.464	3.37E-11
9.468	3.37E-11
9.472	3.37E-11
9.476	3.37E-11
9.48	3.37E-11
9.484	3.37E-11
9.488	3.37E-11
9.492	3.37E-11
9.496	3.37E-11
9.5	3.37E-11
9.504	3.37E-11
9.508	3.37E-11
9.512	3.37E-11
9.516	3.37E-11
9.52	3.37E-11
9.524	3.37E-11
9.528	3.37E-11
9.532	3.37E-11
9.536	3.37E-11
9.54	3.37E-11
9.544	3.37E-11
9.548	3.37E-11
9.552	3.37E-11
9.556	3.37E-11
9.56	3.37E-11
9.564	3.37E-11
9.568	3.37E-11
9.572	3.37E-11
9.576	3.37E-11
9.58	3.37E-11
9.584	3.37E-11
9.588	3.37E-11
9.592	3.37E-11
9.596	3.37E-11
9.6	3.37E-11
9.604	3.37E-11
9.608	3.37E-11
9.612	3.37E-11
9.616	3.37E-11
9.62	3.37E-11
9.624	3.37E-11
9.628	3.37E-11
9.632	3.37E-11
9.636	3.37E-11
9.64	3.37E-11
9.644	3.37E-11
9.648	3.37E-11
9.652	3.37E-11
9.656	3.37E-11
9.66	3.37E-11
9.664	3.37E-11
9.668	3.37E-11
9.672	3.37E-11
9.676	3.37E-11
9.68	3.37E-11
9.684	3.37E-11
9.688	3.37E-11
9.692	3.37E-11
9.696	3.37E-11
9.7	3.37E-11
9.704	3.37E-11
9.708	3.37E-11
9.712	3.37E-11
9.716	3.37E-11
9.72	3.37E-11
9.724	3.37E-11
9.728	3.37E-11
9.732	3.37E-11
9.736	3.37E-11
9.74	3.37E-11
9.744	3.37E-11
9.748	3.37E-11
9.752	3.37E-11
9.756	3.37E-11
9.76	3.37E-11
9.764	3.37E-11
9.768	3.37E-11
9.772	3.37E-11
9.776	3.37E-11
9.78	3.37E-11
9.784	3.37E-11
9.788	3.37E-11
9.792	3.37E-11
9.796	3.37E-11
9.8	3.37E-11
9.804	3.37E-11
9.808	3.37E-11
9.812	3.37E-11
9.816	3.37E-11
9.82	3.37E-11
9.824	3.37E-11
9.828	3.37E-11
9.832	3.37E-11
9.836	3.37E-11
9.84	3.37E-11
9.844	3.37E-11
9.848	3.37E-11
9.852	3.37E-11
9.856	3.37E-11
9.86	3.37E-11
9.864	3.37E-11
9.868	3.37E-11
9.872	3.37E-11
9.876	3.37E-11
9.88	3.37E-11
9.884	3.37E-11
9.888	3.37E-11
9.892	3.37E-11
9.896	3.37E-11
9.9	3.37E-11
9.904	3.37E-11
9.908	3.37E-11
9.912	3.37E-11
9.916	3.37E-11
9.92	3.37E-11
9.924	3.37E-11
9.928	3.37E-11
9.932	3.37E-11
9.936	3.37E-11
9.94	3.37E-11
9.944	3.37E-11
9.948	3.37E-11
9.952	3.37E-11
9.956	3.37E-11
9.96	3.37E-11
9.964	3.37E-11
9.968	3.37E-11
9.972	3.37E-11
9.976	3.37E-11
9.98	3.37E-11
9.984	3.37E-11
9.988	3.37E-11
9.992	3.37E-11
9.996	3.37E-11
10	3.37E-11
};
\end{semilogyaxis}
\end{tikzpicture}

\caption{Error of the energy for the case $r=m=50$ for the mKdV equation: $|H(V\bz _n) - H(V\bz_0)|$ are plotted.
}
\label{fig:mkdv1:EnergyError}
\end{figure}

%% file: figmkdv1sol.tex
\begin{figure}[htbp]
\centering

\begin{tikzpicture}
\tikzstyle{every node}=[]
\begin{axis}[width=6cm,
xmax=5,xmin=-5,
xlabel={$x$},ylabel={$y$},
ylabel near ticks,
legend entries={$r=40$,$r=50$},
legend style={font = \scriptsize,legend pos=north west,legend cell align=left,draw=none,fill=none},
title = {$t=3$},
	]
\addplot[thick,dashed,color=green
] table {
-5	1.029
-4.98	1.021
-4.96	1.003
-4.94	0.977
-4.92	0.943
-4.9	0.905
-4.88	0.865
-4.86	0.826
-4.84	0.789
-4.82	0.757
-4.8	0.732
-4.78	0.715
-4.76	0.707
-4.74	0.709
-4.72	0.719
-4.7	0.737
-4.68	0.762
-4.66	0.79
-4.64	0.821
-4.62	0.852
-4.6	0.88
-4.58	0.902
-4.56	0.916
-4.54	0.921
-4.52	0.916
-4.5	0.898
-4.48	0.869
-4.46	0.828
-4.44	0.777
-4.42	0.718
-4.4	0.651
-4.38	0.58
-4.36	0.507
-4.34	0.435
-4.32	0.366
-4.3	0.303
-4.28	0.247
-4.26	0.2
-4.24	0.162
-4.22	0.135
-4.2	0.117
-4.18	0.109
-4.16	0.108
-4.14	0.114
-4.12	0.123
-4.1	0.134
-4.08	0.145
-4.06	0.154
-4.04	0.158
-4.02	0.156
-4	0.147
-3.98	0.131
-3.96	0.107
-3.94	0.077
-3.92	0.04
-3.9	-0.002
-3.88	-0.046
-3.86	-0.092
-3.84	-0.137
-3.82	-0.178
-3.8	-0.216
-3.78	-0.246
-3.76	-0.27
-3.74	-0.285
-3.72	-0.292
-3.7	-0.291
-3.68	-0.282
-3.66	-0.268
-3.64	-0.249
-3.62	-0.227
-3.6	-0.204
-3.58	-0.183
-3.56	-0.164
-3.54	-0.149
-3.52	-0.14
-3.5	-0.138
-3.48	-0.142
-3.46	-0.153
-3.44	-0.169
-3.42	-0.19
-3.4	-0.214
-3.38	-0.24
-3.36	-0.264
-3.34	-0.286
-3.32	-0.303
-3.3	-0.314
-3.28	-0.316
-3.26	-0.309
-3.24	-0.293
-3.22	-0.266
-3.2	-0.231
-3.18	-0.186
-3.16	-0.135
-3.14	-0.078
-3.12	-0.019
-3.1	0.042
-3.08	0.102
-3.06	0.158
-3.04	0.208
-3.02	0.252
-3	0.287
-2.98	0.313
-2.96	0.331
-2.94	0.339
-2.92	0.339
-2.9	0.333
-2.88	0.321
-2.86	0.306
-2.84	0.29
-2.82	0.274
-2.8	0.26
-2.78	0.25
-2.76	0.244
-2.74	0.245
-2.72	0.251
-2.7	0.264
-2.68	0.281
-2.66	0.303
-2.64	0.327
-2.62	0.353
-2.6	0.377
-2.58	0.399
-2.56	0.417
-2.54	0.428
-2.52	0.433
-2.5	0.429
-2.48	0.417
-2.46	0.397
-2.44	0.369
-2.42	0.334
-2.4	0.294
-2.38	0.252
-2.36	0.208
-2.34	0.165
-2.32	0.125
-2.3	0.091
-2.28	0.063
-2.26	0.043
-2.24	0.033
-2.22	0.032
-2.2	0.04
-2.18	0.057
-2.16	0.082
-2.14	0.112
-2.12	0.147
-2.1	0.184
-2.08	0.222
-2.06	0.257
-2.04	0.289
-2.02	0.315
-2	0.335
-1.98	0.347
-1.96	0.352
-1.94	0.35
-1.92	0.341
-1.9	0.327
-1.88	0.309
-1.86	0.289
-1.84	0.268
-1.82	0.25
-1.8	0.234
-1.78	0.224
-1.76	0.22
-1.74	0.223
-1.72	0.233
-1.7	0.25
-1.68	0.274
-1.66	0.303
-1.64	0.336
-1.62	0.371
-1.6	0.407
-1.58	0.44
-1.56	0.47
-1.54	0.494
-1.52	0.51
-1.5	0.518
-1.48	0.517
-1.46	0.507
-1.44	0.488
-1.42	0.46
-1.4	0.426
-1.38	0.386
-1.36	0.343
-1.34	0.298
-1.32	0.254
-1.3	0.213
-1.28	0.175
-1.26	0.144
-1.24	0.119
-1.22	0.101
-1.2	0.091
-1.18	0.087
-1.16	0.09
-1.14	0.097
-1.12	0.107
-1.1	0.118
-1.08	0.129
-1.06	0.137
-1.04	0.14
-1.02	0.138
-1	0.129
-0.98	0.113
-0.96	0.088
-0.94	0.057
-0.92	0.018
-0.9	-0.025
-0.88	-0.072
-0.86	-0.12
-0.84	-0.168
-0.82	-0.213
-0.8	-0.253
-0.78	-0.287
-0.76	-0.312
-0.74	-0.327
-0.72	-0.332
-0.7	-0.327
-0.68	-0.312
-0.66	-0.287
-0.64	-0.255
-0.62	-0.216
-0.6	-0.172
-0.58	-0.126
-0.56	-0.08
-0.54	-0.035
-0.52	0.007
-0.5	0.045
-0.48	0.077
-0.46	0.103
-0.44	0.124
-0.42	0.139
-0.4	0.15
-0.38	0.158
-0.36	0.165
-0.34	0.173
-0.32	0.183
-0.3	0.198
-0.28	0.218
-0.26	0.246
-0.24	0.282
-0.22	0.326
-0.2	0.379
-0.18	0.44
-0.16	0.508
-0.14	0.58
-0.12	0.656
-0.1	0.733
-0.08	0.809
-0.06	0.881
-0.04	0.947
-0.02	1.006
0	1.056
0.02	1.097
0.04	1.127
0.06	1.147
0.08	1.158
0.1	1.16
0.12	1.156
0.14	1.148
0.16	1.137
0.18	1.127
0.2	1.118
0.22	1.113
0.24	1.114
0.26	1.122
0.28	1.137
0.3	1.159
0.32	1.188
0.34	1.224
0.36	1.263
0.38	1.305
0.4	1.348
0.42	1.389
0.44	1.426
0.46	1.458
0.48	1.483
0.5	1.499
0.52	1.507
0.54	1.507
0.56	1.498
0.58	1.482
0.6	1.461
0.62	1.436
0.64	1.41
0.66	1.384
0.68	1.361
0.7	1.343
0.72	1.331
0.74	1.327
0.76	1.331
0.78	1.343
0.8	1.363
0.82	1.39
0.84	1.423
0.86	1.458
0.88	1.495
0.9	1.53
0.92	1.562
0.94	1.589
0.96	1.608
0.98	1.618
1	1.618
1.02	1.608
1.04	1.589
1.06	1.56
1.08	1.524
1.1	1.482
1.12	1.436
1.14	1.389
1.16	1.342
1.18	1.298
1.2	1.258
1.22	1.225
1.24	1.199
1.26	1.18
1.28	1.169
1.3	1.164
1.32	1.165
1.34	1.17
1.36	1.177
1.38	1.184
1.4	1.189
1.42	1.19
1.44	1.185
1.46	1.173
1.48	1.152
1.5	1.124
1.52	1.086
1.54	1.041
1.56	0.99
1.58	0.933
1.6	0.874
1.62	0.814
1.64	0.756
1.66	0.701
1.68	0.651
1.7	0.609
1.72	0.574
1.74	0.548
1.76	0.53
1.78	0.52
1.8	0.517
1.82	0.519
1.84	0.525
1.86	0.532
1.88	0.539
1.9	0.543
1.92	0.543
1.94	0.536
1.96	0.522
1.98	0.501
2	0.472
2.02	0.435
2.04	0.392
2.06	0.345
2.08	0.293
2.1	0.241
2.12	0.19
2.14	0.141
2.16	0.098
2.18	0.061
2.2	0.032
2.22	0.012
2.24	0.001
2.26	0
2.28	0.006
2.3	0.02
2.32	0.039
2.34	0.062
2.36	0.087
2.38	0.112
2.4	0.134
2.42	0.152
2.44	0.165
2.46	0.172
2.48	0.171
2.5	0.163
2.52	0.148
2.54	0.127
2.56	0.101
2.58	0.072
2.6	0.041
2.62	0.011
2.64	-0.017
2.66	-0.04
2.68	-0.058
2.7	-0.07
2.72	-0.074
2.74	-0.07
2.76	-0.059
2.78	-0.041
2.8	-0.018
2.82	0.009
2.84	0.038
2.86	0.067
2.88	0.095
2.9	0.118
2.92	0.136
2.94	0.147
2.96	0.149
2.98	0.144
3	0.129
3.02	0.107
3.04	0.077
3.06	0.042
3.08	0.002
3.1	-0.04
3.12	-0.081
3.14	-0.121
3.16	-0.157
3.18	-0.187
3.2	-0.211
3.22	-0.226
3.24	-0.233
3.26	-0.231
3.28	-0.222
3.3	-0.207
3.32	-0.185
3.34	-0.161
3.36	-0.134
3.38	-0.108
3.4	-0.084
3.42	-0.064
3.44	-0.05
3.46	-0.042
3.48	-0.042
3.5	-0.049
3.52	-0.064
3.54	-0.084
3.56	-0.11
3.58	-0.14
3.6	-0.172
3.62	-0.203
3.64	-0.233
3.66	-0.258
3.68	-0.278
3.7	-0.291
3.72	-0.297
3.74	-0.294
3.76	-0.283
3.78	-0.265
3.8	-0.24
3.82	-0.211
3.84	-0.177
3.86	-0.143
3.88	-0.108
3.9	-0.077
3.92	-0.049
3.94	-0.027
3.96	-0.013
3.98	-0.005
4	-0.006
4.02	-0.013
4.04	-0.028
4.06	-0.047
4.08	-0.071
4.1	-0.096
4.12	-0.121
4.14	-0.144
4.16	-0.163
4.18	-0.176
4.2	-0.182
4.22	-0.179
4.24	-0.167
4.26	-0.147
4.28	-0.117
4.3	-0.081
4.32	-0.038
4.34	0.009
4.36	0.059
4.38	0.11
4.4	0.16
4.42	0.206
4.44	0.248
4.46	0.284
4.48	0.313
4.5	0.336
4.52	0.352
4.54	0.362
4.56	0.368
4.58	0.37
4.6	0.371
4.62	0.372
4.64	0.375
4.66	0.382
4.68	0.395
4.7	0.415
4.72	0.442
4.74	0.476
4.76	0.519
4.78	0.568
4.8	0.622
4.82	0.68
4.84	0.74
4.86	0.799
4.88	0.856
4.9	0.907
4.92	0.952
4.94	0.987
4.96	1.012
4.98	1.026
};
\addplot[thick,color=magenta
] table {
-5	0.01
-4.98	0.01
-4.96	0.009
-4.94	0.009
-4.92	0.009
-4.9	0.008
-4.88	0.008
-4.86	0.008
-4.84	0.007
-4.82	0.007
-4.8	0.007
-4.78	0.007
-4.76	0.006
-4.74	0.006
-4.72	0.006
-4.7	0.005
-4.68	0.005
-4.66	0.005
-4.64	0.005
-4.62	0.004
-4.6	0.004
-4.58	0.004
-4.56	0.004
-4.54	0.004
-4.52	0.004
-4.5	0.003
-4.48	0.003
-4.46	0.003
-4.44	0.003
-4.42	0.003
-4.4	0.003
-4.38	0.003
-4.36	0.003
-4.34	0.003
-4.32	0.002
-4.3	0.002
-4.28	0.002
-4.26	0.002
-4.24	0.002
-4.22	0.002
-4.2	0.002
-4.18	0.002
-4.16	0.002
-4.14	0.002
-4.12	0.002
-4.1	0.002
-4.08	0.002
-4.06	0.002
-4.04	0.002
-4.02	0.002
-4	0.002
-3.98	0.002
-3.96	0.001
-3.94	0.001
-3.92	0.001
-3.9	0.001
-3.88	0.001
-3.86	0.001
-3.84	0.001
-3.82	0.001
-3.8	0.001
-3.78	0.001
-3.76	0.001
-3.74	0.001
-3.72	0.001
-3.7	0.001
-3.68	0.001
-3.66	0.001
-3.64	0.001
-3.62	0.001
-3.6	0.001
-3.58	0.001
-3.56	0.001
-3.54	0.001
-3.52	0.001
-3.5	0.001
-3.48	0.001
-3.46	0.001
-3.44	0.001
-3.42	0
-3.4	0
-3.38	0
-3.36	0
-3.34	0
-3.32	0
-3.3	0
-3.28	0.001
-3.26	0.001
-3.24	0.001
-3.22	0.001
-3.2	0.001
-3.18	0.001
-3.16	0
-3.14	0
-3.12	0
-3.1	0
-3.08	0
-3.06	0
-3.04	0
-3.02	0
-3	0
-2.98	0
-2.96	0
-2.94	0
-2.92	0
-2.9	0
-2.88	0
-2.86	0
-2.84	0
-2.82	0
-2.8	0
-2.78	0
-2.76	0
-2.74	0
-2.72	0
-2.7	0
-2.68	0
-2.66	0
-2.64	0
-2.62	0
-2.6	0
-2.58	0
-2.56	0
-2.54	0
-2.52	0
-2.5	0.001
-2.48	0.001
-2.46	0.001
-2.44	0.001
-2.42	0.001
-2.4	0.001
-2.38	0.001
-2.36	0.001
-2.34	0.001
-2.32	0.001
-2.3	0.001
-2.28	0.001
-2.26	0.001
-2.24	0.001
-2.22	0.001
-2.2	0.001
-2.18	0.001
-2.16	0.001
-2.14	0.001
-2.12	0.001
-2.1	0.001
-2.08	0.001
-2.06	0.001
-2.04	0.001
-2.02	0.001
-2	0.002
-1.98	0.002
-1.96	0.002
-1.94	0.002
-1.92	0.002
-1.9	0.002
-1.88	0.002
-1.86	0.002
-1.84	0.002
-1.82	0.002
-1.8	0.002
-1.78	0.002
-1.76	0.002
-1.74	0.002
-1.72	0.002
-1.7	0.002
-1.68	0.002
-1.66	0.003
-1.64	0.003
-1.62	0.003
-1.6	0.003
-1.58	0.003
-1.56	0.003
-1.54	0.003
-1.52	0.003
-1.5	0.003
-1.48	0.004
-1.46	0.004
-1.44	0.004
-1.42	0.004
-1.4	0.004
-1.38	0.004
-1.36	0.004
-1.34	0.005
-1.32	0.005
-1.3	0.005
-1.28	0.006
-1.26	0.006
-1.24	0.006
-1.22	0.006
-1.2	0.007
-1.18	0.007
-1.16	0.007
-1.14	0.008
-1.12	0.008
-1.1	0.008
-1.08	0.008
-1.06	0.009
-1.04	0.009
-1.02	0.009
-1	0.01
-0.98	0.01
-0.96	0.011
-0.94	0.011
-0.92	0.012
-0.9	0.012
-0.88	0.013
-0.86	0.013
-0.84	0.014
-0.82	0.014
-0.8	0.015
-0.78	0.016
-0.76	0.016
-0.74	0.017
-0.72	0.018
-0.7	0.018
-0.68	0.019
-0.66	0.02
-0.64	0.02
-0.62	0.021
-0.6	0.022
-0.58	0.023
-0.56	0.024
-0.54	0.025
-0.52	0.026
-0.5	0.027
-0.48	0.028
-0.46	0.029
-0.44	0.031
-0.42	0.032
-0.4	0.033
-0.38	0.035
-0.36	0.036
-0.34	0.037
-0.32	0.039
-0.3	0.04
-0.28	0.042
-0.26	0.044
-0.24	0.045
-0.22	0.047
-0.2	0.049
-0.18	0.051
-0.16	0.053
-0.14	0.055
-0.12	0.058
-0.1	0.06
-0.08	0.063
-0.06	0.065
-0.04	0.068
-0.02	0.071
0	0.074
0.02	0.077
0.04	0.08
0.06	0.083
0.08	0.086
0.1	0.09
0.12	0.093
0.14	0.097
0.16	0.101
0.18	0.105
0.2	0.11
0.22	0.114
0.24	0.119
0.26	0.124
0.28	0.129
0.3	0.134
0.32	0.14
0.34	0.145
0.36	0.151
0.38	0.158
0.4	0.164
0.42	0.171
0.44	0.178
0.46	0.185
0.48	0.192
0.5	0.2
0.52	0.208
0.54	0.217
0.56	0.225
0.58	0.234
0.6	0.244
0.62	0.254
0.64	0.264
0.66	0.275
0.68	0.286
0.7	0.298
0.72	0.31
0.74	0.322
0.76	0.335
0.78	0.349
0.8	0.363
0.82	0.377
0.84	0.392
0.86	0.408
0.88	0.424
0.9	0.441
0.92	0.458
0.94	0.476
0.96	0.495
0.98	0.515
1	0.535
1.02	0.556
1.04	0.577
1.06	0.6
1.08	0.623
1.1	0.647
1.12	0.672
1.14	0.698
1.16	0.725
1.18	0.752
1.2	0.78
1.22	0.809
1.24	0.839
1.26	0.87
1.28	0.902
1.3	0.934
1.32	0.968
1.34	1.002
1.36	1.037
1.38	1.073
1.4	1.11
1.42	1.147
1.44	1.185
1.46	1.223
1.48	1.262
1.5	1.302
1.52	1.341
1.54	1.381
1.56	1.421
1.58	1.461
1.6	1.501
1.62	1.541
1.64	1.58
1.66	1.618
1.68	1.656
1.7	1.692
1.72	1.728
1.74	1.762
1.76	1.795
1.78	1.825
1.8	1.854
1.82	1.881
1.84	1.906
1.86	1.928
1.88	1.947
1.9	1.963
1.92	1.977
1.94	1.988
1.96	1.995
1.98	2
2	2.001
2.02	1.999
2.04	1.993
2.06	1.985
2.08	1.973
2.1	1.959
2.12	1.941
2.14	1.921
2.16	1.898
2.18	1.873
2.2	1.846
2.22	1.816
2.24	1.785
2.26	1.752
2.28	1.717
2.3	1.682
2.32	1.645
2.34	1.607
2.36	1.568
2.38	1.529
2.4	1.489
2.42	1.449
2.44	1.409
2.46	1.369
2.48	1.329
2.5	1.29
2.52	1.25
2.54	1.211
2.56	1.173
2.58	1.135
2.6	1.098
2.62	1.062
2.64	1.026
2.66	0.992
2.68	0.958
2.7	0.924
2.72	0.892
2.74	0.861
2.76	0.83
2.78	0.8
2.8	0.771
2.82	0.743
2.84	0.716
2.86	0.69
2.88	0.664
2.9	0.64
2.92	0.616
2.94	0.593
2.96	0.57
2.98	0.549
3	0.528
3.02	0.508
3.04	0.489
3.06	0.47
3.08	0.453
3.1	0.435
3.12	0.419
3.14	0.403
3.16	0.387
3.18	0.372
3.2	0.358
3.22	0.344
3.24	0.331
3.26	0.318
3.28	0.306
3.3	0.294
3.32	0.282
3.34	0.271
3.36	0.261
3.38	0.251
3.4	0.241
3.42	0.232
3.44	0.223
3.46	0.214
3.48	0.206
3.5	0.198
3.52	0.19
3.54	0.183
3.56	0.176
3.58	0.169
3.6	0.162
3.62	0.156
3.64	0.15
3.66	0.144
3.68	0.138
3.7	0.133
3.72	0.127
3.74	0.122
3.76	0.117
3.78	0.113
3.8	0.108
3.82	0.104
3.84	0.1
3.86	0.096
3.88	0.092
3.9	0.089
3.92	0.085
3.94	0.082
3.96	0.079
3.98	0.076
4	0.073
4.02	0.07
4.04	0.067
4.06	0.064
4.08	0.062
4.1	0.059
4.12	0.057
4.14	0.055
4.16	0.052
4.18	0.05
4.2	0.048
4.22	0.046
4.24	0.045
4.26	0.043
4.28	0.041
4.3	0.04
4.32	0.038
4.34	0.037
4.36	0.035
4.38	0.034
4.4	0.033
4.42	0.032
4.44	0.03
4.46	0.029
4.48	0.028
4.5	0.027
4.52	0.026
4.54	0.025
4.56	0.024
4.58	0.023
4.6	0.022
4.62	0.021
4.64	0.02
4.66	0.02
4.68	0.019
4.7	0.018
4.72	0.018
4.74	0.017
4.76	0.016
4.78	0.016
4.8	0.015
4.82	0.015
4.84	0.014
4.86	0.013
4.88	0.013
4.9	0.012
4.92	0.012
4.94	0.011
4.96	0.011
4.98	0.01
};
\end{axis}
\end{tikzpicture}

\caption{Numerical solutions $V\bz$ at $t=3$ for $r=m=40,50$ for the mKdV equation.
}
\label{fig:mkdv1:sol}
\end{figure}